\documentclass[fleqn,dvips]{article}
\usepackage[english]{babel}   
\usepackage{graphicx}
\usepackage{a4wide}
\usepackage{color}
\usepackage{xspace}  
\usepackage{amssymb}
\usepackage{array}
\usepackage{subfigure}
\usepackage{pictex}
\usepackage{multirow}
\def\e{\varepsilon}

\def\l{\lambda}

\newcommand{\R}{\mathbb{R}}

\newcommand{\N}{\mathbb{N}}
\newcommand{\C}{\mathbb{C}}

\newcommand{\pa}{\partial}

\newcommand{\beq}{\begin{eqnarray}}
\newcommand{\eeq}{\end{eqnarray}}
\newcommand{\beqs}{\begin{eqnarray*}}
\newcommand{\eeqs}{\end{eqnarray*}}
\newcommand{\bequ}{\begin{equation}}
\newcommand{\eequ}{\end{equation}}

\catcode`\@=11

\@addtoreset{equation}{section}

\begin{document}
\bibliographystyle{abbrv}
\title{
Pattern formation in a pseudo-parabolic equation
}
\author{ 
C.M. Cuesta\footnote{University of the Basque Country (UPV/EHU), Mathematics Department, Aptdo. 644, 48080 Bilbao, Spain, e-mail: carlotamaria.cuesta@ehu.es
}, 
J.R. King\footnote{University of Nottingham, Division of Applied Mathematics, School of Mathematical Sciences, 
University of Nottingham, NG7 2RD Nottingham, UK, 
e-mail: john.king@nottingham.ac.uk 
}, 
}

\date{}
\maketitle
\begin{abstract}
We address the propagation into an unstable state of a localised disturbance in the pseudo-parabolic equation
\[
\frac{\pa u}{\pa t} = \frac{\pa^2}{\pa x^2}\left(\phi(u) + \frac{\pa u}{\pa t}\right) \,, 
\]
where $\phi$ is a non-monotone function. We concentrate on the representative odd nonlinearities $\phi(u)=u^3 - u$ 
and $\phi(u)=-u e^{-u^2}$, and take the unstable state to be $u_u\equiv 0$ for most of the analysis. Three asymptotic regimes are distinguished as 
$t\to +\infty$, the first being a regime ahead of the propagating disturbance that is dominated by the linearised equation. The analysis of this 
leads to the determination of the speed of the leading edge of the propagating disturbance and implies that in the second, transition, regime the 
solution takes the form of a modulated travelling wave. In a third regime the solution approaches a nearly periodic steady state, where the period 
is obtained on matching with the modulated travelling wave. Detailed analysis of this pattern is also presented. The analysis is completed by 
contrasting the formal asymptotic description of the solution with numerical computations. It is assumed for the above analysis that the initial 
disturbance decays faster than an exponential rate; in this case a critical exponential decay rate at the leading edge of the front and propagation speed are found. 
We investigate the wave speed selection mechanism for exponentially decaying initial conditions. It is found that whenever the initial data behave as a real exponential (no matter how slow the rate of the decay) the speed selected is that selected by fast decaying initial conditions. 
However, for initial conditions with a complex exponential, thus allowing oscillatory perturbations, we find regimes of the decay rate and the wavelength for which the front 
propagates at a faster wave speed. This is investigated numerically and is worth emphasising since it gives a different scenario for wave speed behaviour than that exhibited by well-studied semilinear reaction-diffusion equations: there are initial conditions with exponential decay faster 
than the critical one for which the front propagates with a speed faster than the critical one. 
\end{abstract}
\section{Introduction}\label{section:1}
In this paper we study pattern formation initiated by a localised disturbance for the pseudo-parabolic equation
\begin{equation}\label{main:eq}
\frac{\pa u}{\pa t} = \frac{\pa^2}{\pa x^2}\left(\phi(u) + \frac{\pa u}{\pa t}\right) , \quad x \in \R, \quad t>0 \, , 
\end{equation}
subject to an initial condition
\begin{equation}\label{main:ic}
u(x,0)=u_0(x), \quad x \in \R \, 
\end{equation}
where the nonlinearity $\phi$ is a smooth non-monotone function. Such formulations arise in a number of physical and biological applications, as we outline below. We are interested in the dynamics around unstable states: we analyse, by means of matched asymptotics, 
front propagation into unstable states, i.e. the mechanisms by which, under an initial perturbation of 
an unstable state, stable patterns `win' over the unstable state, invading its domain. Before we go into this matter, let us recall some properties of (\ref{main:eq}).

Equation (\ref{main:eq}) typically appears as a so-called Sobolev regularisation 
(cf. \cite{ER2}) of the forward-backward diffusion equation
\begin{equation}\label{diff:eq}
\frac{\pa u}{\pa t }=\frac{\pa^2}{\pa x^2}(\phi(u)) , \quad x \in \R, \quad t>0 \, .
\end{equation}
Observe that in regions where $\phi'(u)<0$ equation (\ref{diff:eq}) is backward-parabolic and, 
thus, ill-posed; in particular, H\"ollig showed in \cite{hollig} that if 
$\phi$ is piecewise linear then
there exist initial conditions for which the Cauchy problem has infinitely many solutions.
 Uniqueness can be achieved by introducing a higher-order regularisation such as a fourth-order 
term, as in the Cahn-Hilliard equation, or a third-order term with mixed derivatives as in (\ref{main:eq}), cf. Latt\`es and Lions \cite{Lions}. 
The limit $\e\to 0$ of the Sobolev regularisation for a cubic nonlinearity such as $\phi(u)=u^3-u$ is studied rigorously 
in Plotnikov \cite{Plotnikov}; see also e.g. Evans and Portilheiro \cite{Evans} and Mascia, Terracina and Tesei \cite{terracina} and \cite{terracina2}, Gilding and Tesei \cite{GilTesei} and Lafitte and Mascia \cite{LfftM}. 

\begin{figure}[htb]
\centering
\mbox{
\subfigure[$\phi(u) = u^3 - u$.]
{
\includegraphics[width=0.295\textwidth,height=.28\textwidth]{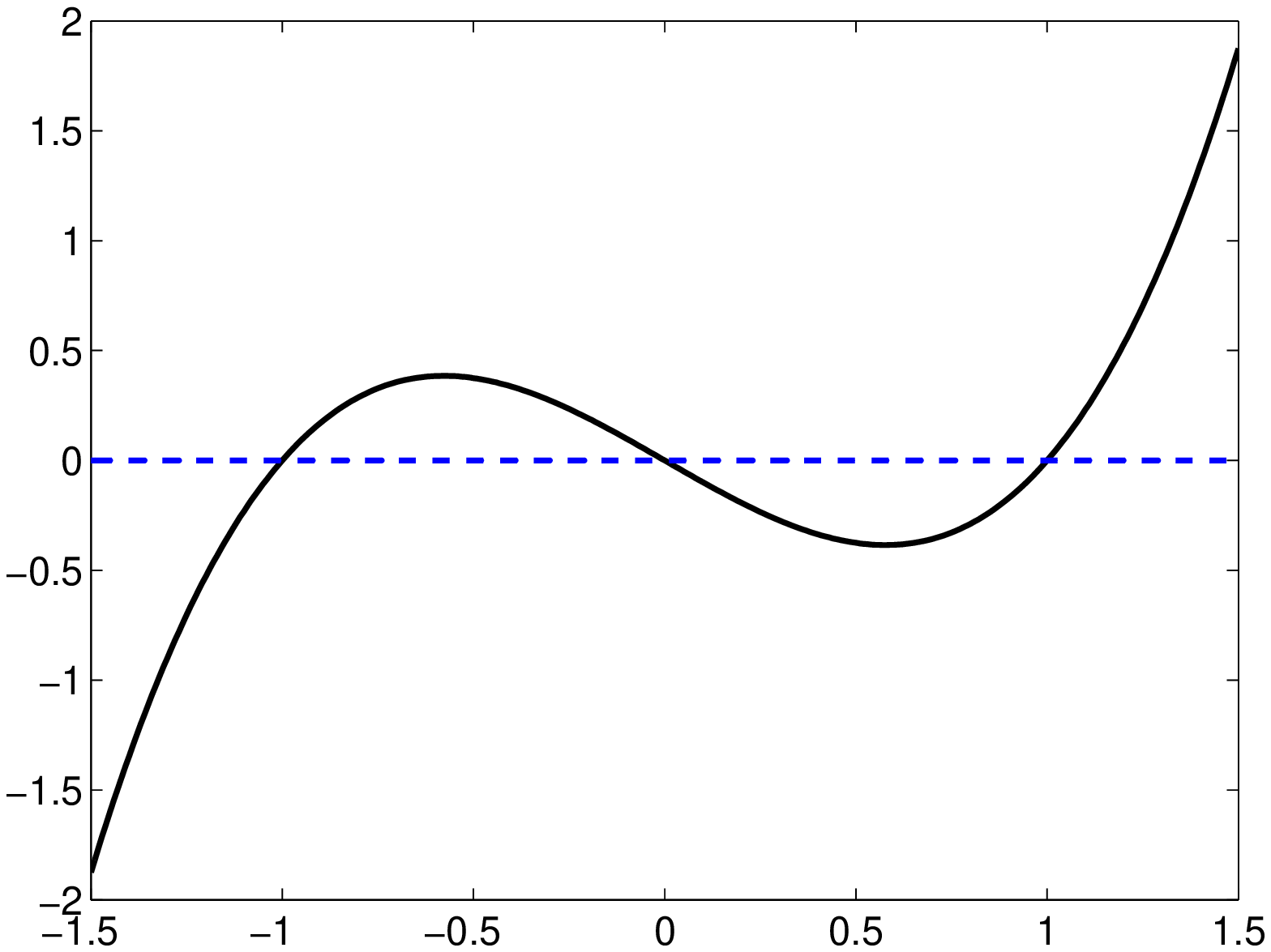}
\label{cubic}
}
}
\mbox{
\subfigure[$\phi(u) = u e^{-u}$.] 
{
\includegraphics[width=0.295\textwidth,height=.28\textwidth]{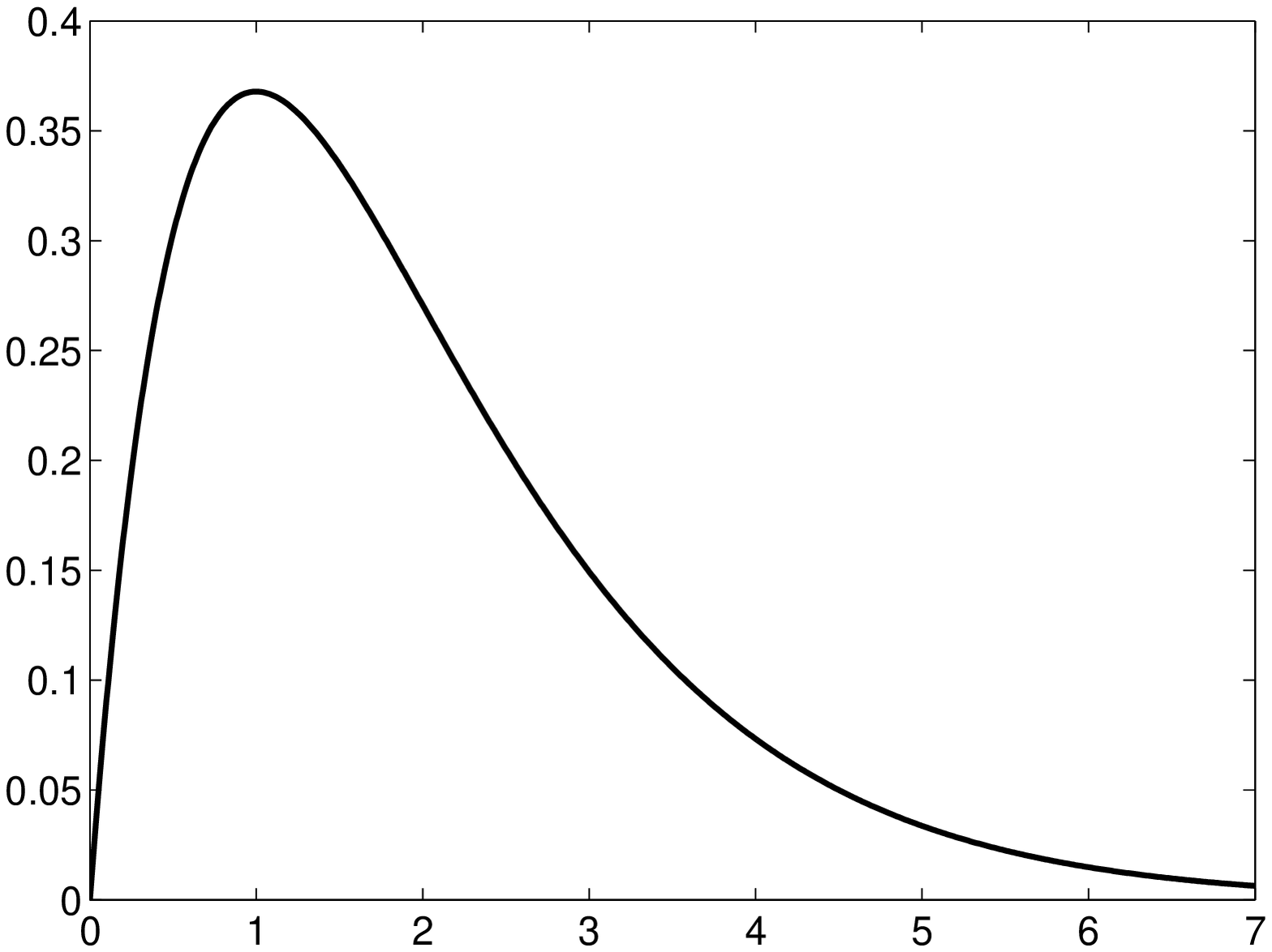}
\label{expo}
}
}
\mbox{
\subfigure[$\phi(u) = -u e^{-u^2}$.] 
{
\includegraphics[width=0.295\textwidth,height=.28\textwidth]{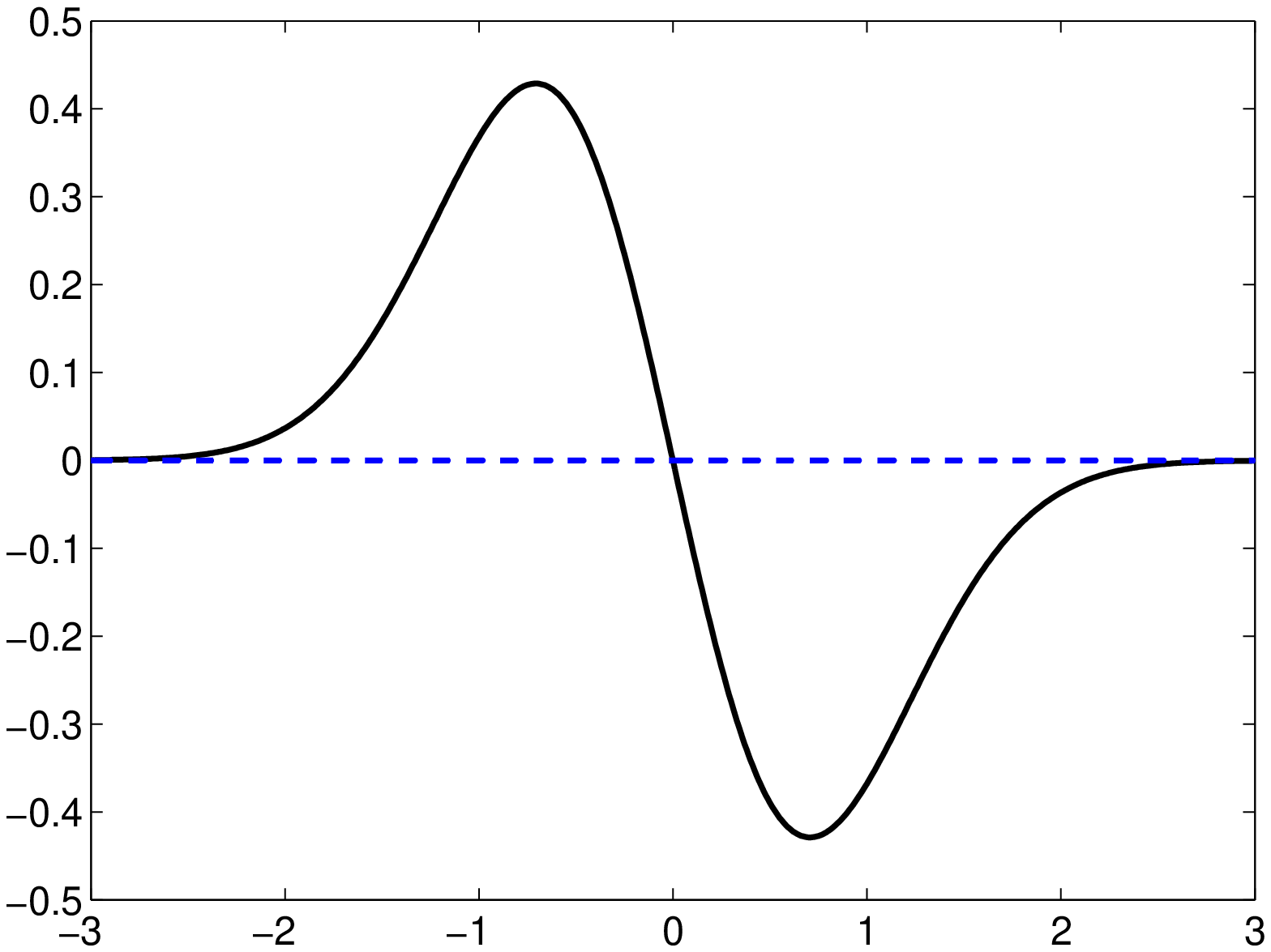}
\label{exposym}
}
}
\caption{Exemplar nonlinearities $\phi(u)$.}
\label{phi}
\end{figure}

Existence and regularity properties of (\ref{main:eq}) were derived in \cite{Cohen} and \cite{Padron} using 
different approaches.
It is well known that pseudo-parabolic equations (at least when the higher-order regularisation is linear) preserve 
the regularity in space of the initial data. For example, if the initial condition has a 
jump discontinuity at 
$x=x_0$ then the solution has a (time-dependent) jump discontinuity at $x=x_0$ for all $t>0$, satisfying
\[
\frac{d}{dt}\left[u\right]_-^+ = -\left[\phi(u)\right]_-^+\,.
\]
Global existence holds in the 
positively invariant regions, i.e. the regions in which $\phi'(u)>0$ (see e.g. \cite{Cohen}). 
We also observe that the zeroth moment (mass) and first moment are conserved, i.e.
\bequ\label{mass}
\int_{-\infty}^{\infty}\left(u(x,t) -u_0(x)\right)\,dx=0 \,, 
\ \int_{-\infty}^{\infty} x \left(u(x,t)-u_0(x)\right)\,dx=0 \quad\mbox{for} \ t>0 \,.
\eequ

The steady states of equations (\ref{main:eq}) and (\ref{diff:eq}) satisfy
\begin{equation}
\label{phi:const}
\phi(u)=A\, 
\end{equation}
for some constant $A$, but for non-monotonic $\phi$ this need not imply that $u$ is constant; 
clearly, any constant solution is a steady state, however. Linearisation shows that constant steady states 
$u_s$ such that $\phi'(u_s)>0$ are linearly stable, and we term this domain the {\it stable region}. 
Those satisfying $\phi'(u_s)<0$ are linearly unstable, this domain being the {\it unstable region}. 
More complicated stationary patterns, namely any piecewise combination of constant solutions satisfying 
(\ref{phi:const}), require non-trivial stability analysis. For $\phi$ of the form (\ref{phi:model1}), 
stability (to small perturbations) of (any) piecewise-constant steady states $u_s(x)$ satisfying (\ref{phi:const}) and $\phi'(u_s(x))>0$ a.e. 
was proved in \cite{Cohen}.


In what follows the constant $u_u$ will 
denote an unstable state, namely a constant such that $\phi'(u_u)<0$; in all cases the solution and the initial condition will be taken to satisfy
\bequ\label{far:field:uns}
u_0\,,\ u \to u_u \quad \mbox{as} \ |x| \to \infty\,. 
\eequ 
and the initial condition is a small localised perturbation to $u_u$. One expects that for $t$ sufficiently large and increasing, the perturbation to $u_u$ will grow and spread, invading the domain in both directions and leaving behind a pattern which approaches a steady state. In what follows we give some of the ingredients for analysing the associated propagating front and resulting pattern for equation (\ref{main:eq}). We recall that we shall apply matched-asymptotic methods in identifying three distinguished regimes of the solution. In the regime ahead of the front, the dominant balance as $t\to \infty$ is the equation linearised around the unstable state $u_u$, this regime applying back to the leading edge of the front where the perturbation grows to become non-negligible (so that linearisation is no longer appropriate). There is a second, transition, regime at the front, where the growing component of the perturbation is controlled by the nonlinearity (leading to a modulated travelling wave), and a third one where the pattern is established.

\begin{figure}
\centerline{
\includegraphics[width=0.65\textwidth,height=.34\textwidth]{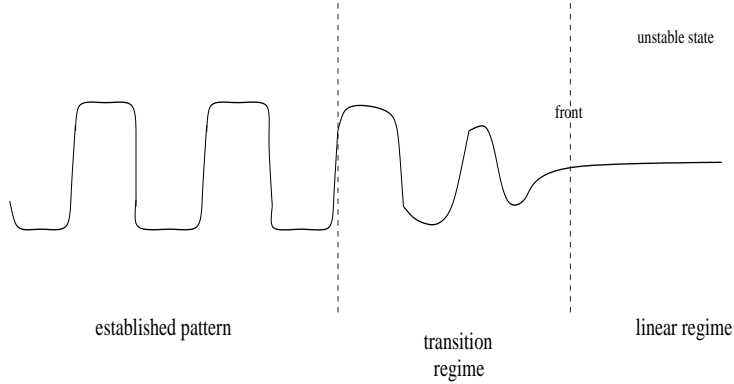}
}
\caption{Sketch of the asymptotic regimes for a front propagating to the right.}
\label{regimes}
\end{figure}

We shall concentrate on the example shown in Figure~\ref{cubic}, 
namely 
\bequ\label{phi:cubic}\phi(u) = u^3 - u
\eequ 
for a $\phi$ satisfying 
\begin{equation}\label{phi:model1}
\begin{array}{l}
\phi'(u) > 0  \quad \mbox{for} \ u \in (-\infty,u_M)\cup(u_m, \infty)\, ,   \\
\phi'(u) < 0 \quad \mbox{for} \ u \in (u_M,u_m)\, , \\ 
\lim_{u \to\pm \infty} \phi(u) = \pm \infty\, ,
\end{array}
\end{equation}
having the local maximum at $u_M$ and the local minimum at $u_m$. 
We expect that the perturbations will evolve to a stable piecewise constant solution taking distinct values $u_-$ and $u_+$ (with $\phi'(u_+)$, $\phi'(u_-)>0$), 
rather than to a constant solution, since (\ref{mass}) must be satisfied. 

It is instructive also to consider a second example, see Figure~\ref{exposym}, 
where $\phi$ satisfies (\ref{phi:model1}), but with the last condition replaced by 
\[
\lim_{u \to \pm \infty}\phi(u)=0\,,
\] 
namely 
\bequ\label{phi:exposym}
\phi(u) = -u e^{-u^2}\,.
\eequ 
For the nonlinearity (\ref{phi:exposym}) the case $u_u=0$ is particularly interesting, since there are {\it no} values $u_s$ in the stable region satisfying $\phi(u_s)=0$, and we expect a situation where the solution alternates between time-dependent values that approach $+\infty$ and $-\infty$ as 
$t\to +\infty$. Thus, this example illustrates a class in which $u$ grows unboundedly behind the advancing front 
and its analysis complements instructively that of (\ref{phi:cubic}). 
Observe that $\phi$ in (\ref{phi:cubic}) and (\ref{phi:exposym}) are odd 
functions of $u$; for these two cases we shall speak about the {\it symmetric} case when $u_u=0$. 
The current paper is restricted to these cases, although part of the analysis presented here is more general. 
We shall point out some additional complications that appear in the non-symmetric cases as appropriate. For example, we shall see later that, if we let $A:=\phi(u_-)=\phi(u_+)$, then in general $A \neq\phi(u_u)$, but the symmetric case (\ref{phi:cubic}) with $u_-=-1$, $u_u=0$ and $u_+=1$ has $\phi(u_-)=\phi(u_+)=\phi(0)=0$.


In Section~\ref{section:2} we analyse in some generality the regime ahead of the front and the wave speed selection mechanism for large times. 
We start by computing the speed of the front propagating into the unstable state. 
We distinguish three regimes (see Figure~\ref{regimes}): in the first, the 
equation linearised about the unstable state, namely, setting $u \sim u_u +  v$ with $|v|\ll 1$, 
\bequ\label{lin:equ}
\frac{\pa v}{\pa t} = -\Phi_u \frac{\pa^2 v}{\pa x^2} + \frac{\pa^3 v}{\pa x^2 \pa t}\, , \quad \mbox{with} \quad \Phi_u:=-\phi'(u_u)>0\,,
\eequ         
gives the dominant balance (in a sense made more precise below) 
and the position of the front can be determined from (\ref{lin:equ}) using an appropriate condition to detect the leading edge of 
the wave front. One approach to pursue such analysis is reviewed in \cite{frontsreview}, 
cf. \cite{kendall} for what amounts to an early application of such a procedure. In this approach one computes the Fourier transform 
of the linearised equation and solves the 
transformed equation, the inverse Fourier transform of this solution being approximated as $t\to +\infty$ 
by the steepest descent method. We adopt an alternative approach that is based on the Liouville-Green (JWKB) method, as we 
shall describe in Section~\ref{section:2}. The analysis allows in particular the (critical) front speed $\xi^*$ to be computed. In Section~\ref{section:2.2} we present this analysis requiring that the initial perturbation decays as $x \to \pm\infty$ at a rate greater than exponential. Under this condition it also shows that solutions decay exponentially with a characteristic (linearly selected) rate $\lambda^*$.
Then Section~\ref{sec:expdecay} is devoted to clarifying the front speed selection mechanism. 
For exponentially decaying perturbations
\bequ\label{exp:ic}
u_0(x)=u_u+ v_0(x) \quad \mbox{with} \quad v_0(x)\sim \e e^{-\lambda |x|}\quad \mbox{as} \ x\to +\infty \quad \mbox{and} \ 0<\e\ll 1\,
\eequ
additional exponential contributions, associated with the separation-of-variables solution
\bequ\label{slowdecay}
v(x,t) = e^{-\Phi_u \frac{\l^2}{1-\l^2}t -\l x} \quad \mbox{with}\quad \mbox{Re}(\l)>0\,,
\eequ
need to be taken into account. In principle, if these contributions 
are dominant they could lead to a front that propagates at a faster speed $\xi_f(\lambda)$ than the critical one, $\xi^*$. 
Such a mechanism is well understood for the Fisher equation and other systems for which the front solutions take the form of a 
travelling wave, cf. \cite{frontsreview}, \cite{CK}. For the Fisher equation with initial data of the form (\ref{slowdecay}), the selected front speed $\bar{\xi}(\lambda)$ is a smooth decreasing function with $\bar{\xi}(\lambda)\to +\infty$ as $\lambda\to 0^+$ and $\bar{\xi}(\lambda)=\xi^*$ for all $\lambda\geq \lambda^*$. We analyse this issue for (\ref{main:eq}) in Section~\ref{sec:expdecay}, by taking $\l\in \R$ and comparing the exponential contributions in the complex plane ($\xi$ complex). This requires a computation of the associated Stokes lines, since the various exponential 
contributions can be switched on or off across these Stokes lines (i.e. where the imaginary parts of every two exponents coincide and the exponential that is 
being switched is maximally subdominant to the other exponential). The analysis concludes that $\bar{\xi}(\l)=\xi^*$ for any $\l\in\R$, in marked contrast to the familiar Fisher case.

In Section~\ref{sec:3} we analyse the asymptotic regions. First, in Section~\ref{section:2.3new} we analyse the leading edge of a front propagating with speed $\xi^*$. Of the two further regimes that we distinguish, one (the transition regime) propagates 
with the front, while the third is that in which the pattern has been established. We analyse the former in Section~\ref{section:2:3}. In this transition regime the solution is asymptotically periodic in time in the coordinate system moving with the front speed (i.e. it takes the form of a modulated travelling wave). 
As we shall see, this suggests that the pattern that it lays down is periodic in space. The numerical examples presented in Section~\ref{section:3} are consistent with the conjecture that the spatial and temporal periods, as well as 
the wave speed, can be computed from the previous analysis of the linear regime. I.e. that the current problem 
is of the `pulled front' class according to the terminology adopted in e.g. \cite{frontsreview}\footnote{According to the terminology adopted in \cite{saarloos}, {\it pulled}  (linearly-selected) fronts are those 
which propagate at speed that can be compute from the linearised equation, while {\it pushed} (nonlinearly-selected) fronts propagate at a, in general, faster speed, which can be determined only by a nonlinear analysis of the associated (modulated) travelling wave. The mechanisms underlying pushed fronts are often 
subtle and, in general, specific to the equation; see \cite{aronson} for semi-linear parabolic equations and \cite{frontsreview} and \cite{CK} 
for a more general discussion. As we outline later, wave speeds greater than $\xi^*$ can also result for specific types of initial data; in contrast to pushed fronts, the wave speed in these cases can again be determined by linear arguments.} We complete the asymptotic analysis of the pattern in the 
symmetric cases in Section~\ref{section:3.2}. As mentioned above, the pattern alternates between two values. 
The transition between them is rather sharp, though necessarily continuous for continuous initial conditions. We analyse in this section the internal structure of these transitions. 

In Section~\ref{section:3} we check numerically the predictions for the front speed and the spatial period derived in the previous sections. For (\ref{phi:cubic}) we show that the resulting pattern $u_s(x)$ 
alternately takes (in a near-periodic fashion) the values $-1$ and $1$, 
implying that $\phi(u_s(x))=\phi(0)=0$. This relation is later confirmed numerically. 
For (\ref{phi:exposym}) we show that the periodic pattern laid down behind the front becomes unbounded 
as $t\to +\infty$, and we estimate this growth and compare it with the numerical results.

We complete the analysis in Section~\ref{section:new} by considering initial conditions of the form (\ref{slowdecay}) with $\l\in\C$. A continuity argument 
(into the complex $\l$-plane) gives the regions of $\l$ that lead to faster than critical wave speeds. The section is completed with numerical experiments 
for several values of $\l$. In particular, in this section we find that there are values of $\l\in\C$ with $\mbox{Re}(\l)>\l^*$ that have $\bar{\xi}(\l)>\xi^*$ and this is verified by the numerical results. 

We end this introduction by placing (\ref{main:eq}) into a wider context. 
Equation (\ref{main:eq}) seems to have been first considered by Novick-Cohen and Pego in \cite{Cohen}.
There the nonlinearity $\phi$ was taken to satisfy (\ref{phi:model1}) (cf. Figure~\ref{cubic}). Equation (\ref{main:eq}) can thus be viewed as the limiting case of the 
viscous Cahn-Hilliard equation, namely (in the notation of \cite{Cohen2})
\begin{equation}\label{cahn-hilliard}
\frac{\pa u}{\pa t} = \frac{\pa^2}{\pa x^2}\left(
\phi(u) -\delta^2 \frac{\pa^2 u}{\pa x^2} + \e^2 \frac{\pa u}{\pa t}
\right)\, ,
\end{equation} 
whereby the interfacial energy (see the final term in (\ref{free:energy})) is negligible (i.e. $\varepsilon=1$, $\delta =0$; as we shall see, the absence 
of penalisation of interfaces in this case has important implications for the dynamics). 
The third-order term was introduced in \cite{Cohen2} to account for viscous relaxation effects. 
The widely studied Cahn-Hilliard equation ($\e=0$ in (\ref{cahn-hilliard})), arises as a model for phase separation by 
spinoidal decomposition of a binary mixture, see \cite{cahn-hilliard}. In higher dimensions, 
the constant-mobility Cahn-Hillard equation reads
\bequ\label{close:cahn:hilliard}
\frac{\partial u}{\partial t} =\Delta \mu\,, 
\eequ
where the unknown $u$ represents the concentration of one of the two phases and 
$\mu=\phi(u)-\delta^2 \Delta u$ is the chemical 
potential, which is the functional derivative with respect to $u$ of the free energy ${\cal L}$ 
for a given volume $\Omega$, with
\bequ\label{free:energy}
{\cal L}(u)= \int_\Omega\left( \Phi(u)+\frac{1}{2} \delta^2 |\nabla u|^2\right)\, d{\bf x}\,,
\eequ
where $\phi(u)=\Phi'(u)$, the contribution $\Phi$ to the free energy per unit volume typically 
being taken to be a double-well potential. When $\Phi$ is convex, separation of the phases does not occur, 
whereas 
if $\Phi$ has two minima, corresponding to different concentration levels, separation occurs 
in which the minima are attained, the final state not being an homogeneous mixture. 
We observe that the viscous version of the Cahn-Hilliard equation can be seen as a 
Sobolev regularisation of (\ref{close:cahn:hilliard}), namely
\[
\frac{\pa}{\pa t}(u-\e^2 \Delta u)=\Delta \mu\,.
\]

Equation (\ref{cahn-hilliard}) exhibits metastable solutions characterised by (slowly-evolving) 
alternating regions of high and low concentration when $\Phi$ has minima of equal depth; see for example Reyna and Ward \cite{Ward}. 
Competition between these regions leads to phase coarsening, typically ending up in total separation of the phases; see \cite{nicolaenkoscheurer} and \cite{zheng}. 
Importantly, such behaviour does not occur 
when $\delta=0$ in (\ref{cahn-hilliard}). For a review on the Cahn-Hilliard equation and related 
models of phase separation we refer to \cite{fife}.

Equation (\ref{main:eq}) was also considered in \cite{Padron} 
as a model of aggregating populations. 
In this case, however, the nonlinearity $\phi$ was taken to be of the general form
\begin{equation}\label{phi:model2}
\begin{array}{l}
\phi'(u) > 0  \quad \mbox{for} \ u \in (u_-,u_M)\, , \\
\phi'(u) < 0 \quad \mbox{for} \ u \in (u_M,\infty)\,, \\ 
 \lim_{u\to \infty} \phi(u)  < +\infty \,
\end{array}
\end{equation}
(cf. Figure~\ref{expo}). 
Here $u$ stands for the population density. The nonlinearity reads $\phi(u) =u \varphi(u)$, where 
$\varphi$ is the migration rate, and is taken to be a decreasing function: as the population grows the 
tendency of individuals to migrate diminishes, so that $\phi'(u)<0$ for $u$ sufficiently large. 
The third-order term is introduced here as a regularisation of the ill-posed problem. Observe that 
in this model no growth interaction (births and deaths) has been included. 

If $\phi$ is of the form (\ref{phi:model2}) there can only be one $u$ in the stable region satisfying  
(\ref{phi:const}). In this case the solution might be expected to approach the only available constant stable 
solution. It is, however, not immediately clear how such a solution arranges itself in space as $t\to +\infty$, since 
the conditions (\ref{mass}) hold; we venture that the solution oscillates spatially between a stable value $u_s$ 
and values that tend to infinity as $t\to +\infty$, presumably approximating a function of the form $u(x)=u_s + \sum_{n=1}^N M_n \,\delta(x-x_n)$.

A related model, where $\phi$ is of the form (\ref{phi:model2}) but where the third-order term
is quasilinear, appears in \cite{Bar3} and \cite{Bar2} 
as a model of heat and mass transfer in turbulent shear flows. Variants of such model equations also appear in several other biological applications, 
see for instance \cite{Oliver} and \cite{lemon}; we shall not deal with such more complicated models here.

\section{Front speed selection}\label{section:2}
In this section we give the preliminary analysis of the linear regime and derive the linearly selected wave speed depending on the decay of the initial condition.

\subsection{The WKBJ approach}\label{section:2.2} 
We consider a general nonlinearity $\phi$. In the linear-dominated regime, linearisation around the constant unstable state $u_u$ gives the leading-order equation (\ref{lin:equ}). 
We use a JWKB method for solving it, thus we adopt the usual ansatz
\bequ\label{wkbj:ansatz}
v(x,t)\sim  a(x,t) e^{-{f(x,t)}}\,,
\eequ
whereby, ahead of the front, the function $a$ influences the amplitude of oscillations and $\mbox{Im}(f)$ determines  
their frequency, and $\mbox{Re}(f)$ records the decay of the solution. 
Substituting (\ref{wkbj:ansatz}) into (\ref{lin:equ}) gives at leading order
\bequ\label{wkbj:balance}
\frac{\pa f}{\pa t} =\Phi_u \left( \frac{\pa f}{\pa x }\right)^2 + \left(\frac{\pa f}{\pa x }\right)^2 \frac{\pa f}{\pa t}\quad \mbox{as} 
\ t\to +\infty \quad\mbox{with} \quad \frac{x}{t}=O(1)\,, 
\eequ
and an amplitude equation for $a$ follows from the next order balance:
\bequ\label{wkbj:amplitude}
\left( 1-\left(\frac{\pa f}{\pa x}\right)^2\right)\frac{\pa a}{ \pa t} - 2 \left(\Phi_u + \frac{\pa f}{ \pa t}\right)\frac{\pa f}{\pa x} \frac{\pa a}{\pa x}
=\left\{ \left(\Phi_u +\frac{\pa f}{\pa t} \right) \frac{\pa^2 f}{\pa x^2}+ 2\frac{\pa f}{\pa x}\frac{\pa^2 f}{\pa x\pa t}\right\}a \,,
\eequ
where $f$ satisfies (\ref{wkbj:balance}).

Setting $f(x,t)= t F(\xi)$ with $\xi=x/t$ in (\ref{wkbj:balance}), we obtain the Clairaut equation 
\bequ
F-\xi \frac{d F}{d\xi} =\Phi_u\frac{ \left(\frac{d F}{d\xi}\right)^2}{ 1-\left(\frac{d F}{ d\xi} \right)^2} \,,
\label{clairaut}
\eequ
the general solutions to which are 
\bequ\label{clairaut:general}
F(\xi)= \l\xi + \Phi_u\frac{\l^2}{1 -\l^2}\quad \mbox{for all} \ \l\neq 1 \,,
\eequ
while the singular solution is given parametrically in terms of $p=dF/d\xi$ by
\begin{eqnarray}\label{clairaut:singular}
F & = & -\Phi_u\left( \frac{ 2p^2}{(1-p^2)^2 } - \frac{p^2}{1-p^2}\right)\,, 
\\
\xi & = & - \Phi_u\frac{ 2 p}{(1-p^2)^2}\,. \label{saddle:eq}
\end{eqnarray}
The graph of the solution (\ref{clairaut:singular})-(\ref{saddle:eq}) is the envelope of the graphs of the family 
(\ref{clairaut:general}). Equations (\ref{clairaut:general}) and (\ref{clairaut:singular})-(\ref{saddle:eq}) of course correspond to solutions of
 the Charpit equations for (\ref{wkbj:balance}), which have $p=\pa f/\pa x$ and $\pa f/\pa t$ constant along the 
rays 
\bequ\label{rays}
x(t)=-\Phi_u \frac{2p }{(1-p^2)^2}\, t + x_0\,,
\eequ 
with $x(0)=x_0$. The general solution with Cauchy data $f(x,0)=f_0(x)$ is 
\[
f(x(t),t)= -\Phi_u\left( \frac{ 2p^2}{(1-p^2)^2 } - \frac{p^2}{1-p^2}\right)t + f_0(x_0)
\]
which for $f_0(x_0)=\lambda x_0$, $p=\lambda$ gives (\ref{clairaut:general}), and for $x_0=0$ with $f_0(0)=0$ (and $p$ arbitrary) gives (\ref{clairaut:singular})-(\ref{saddle:eq}).

\begin{figure}
\centering
\mbox{
\subfigure[]
{
\includegraphics[width=0.45\textwidth,height=.30\textwidth]{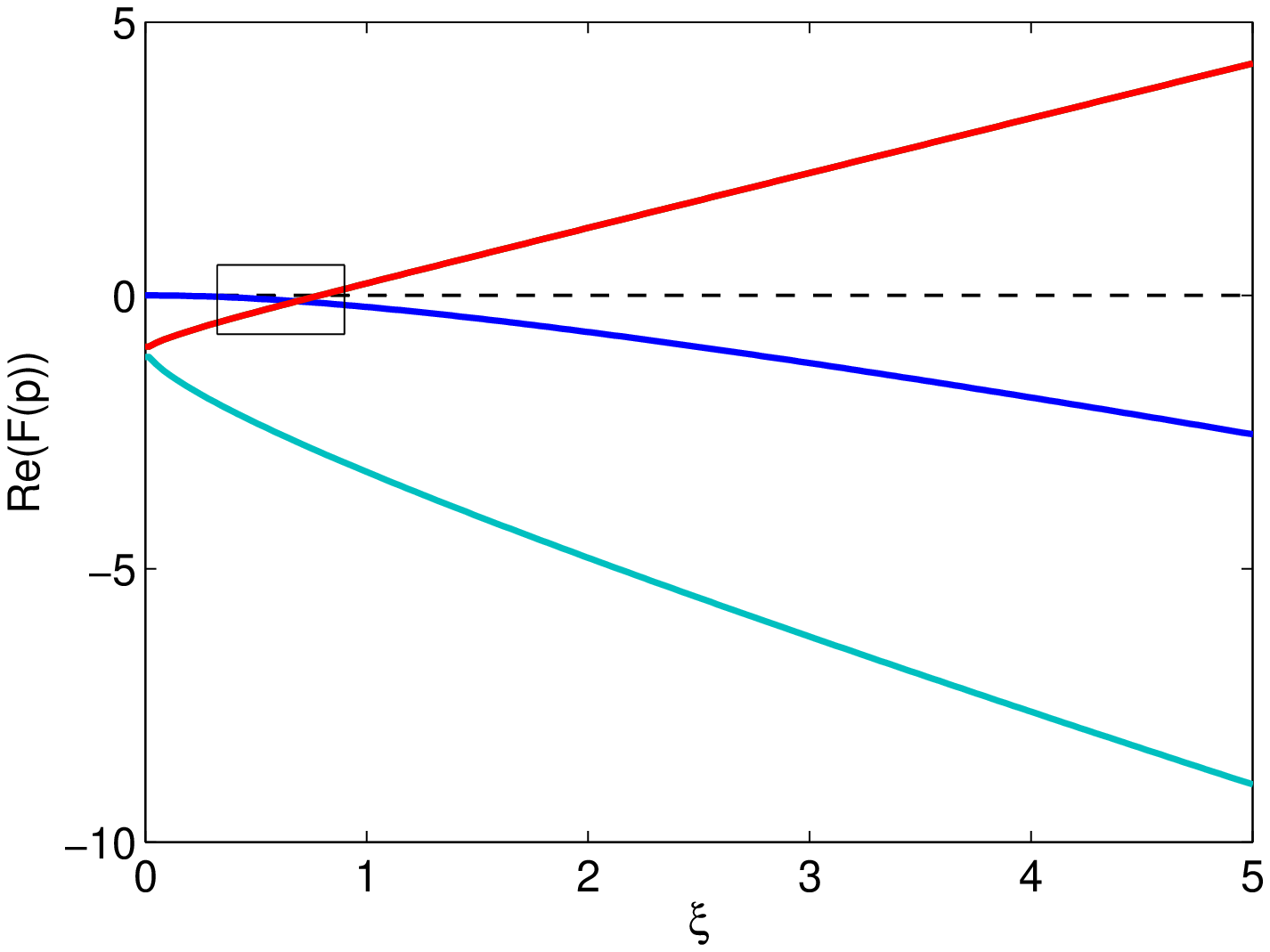}
\label{speedI}
}
}
\quad
\mbox{
\subfigure[] 
{
\includegraphics[width=0.45\textwidth,height=.30\textwidth]{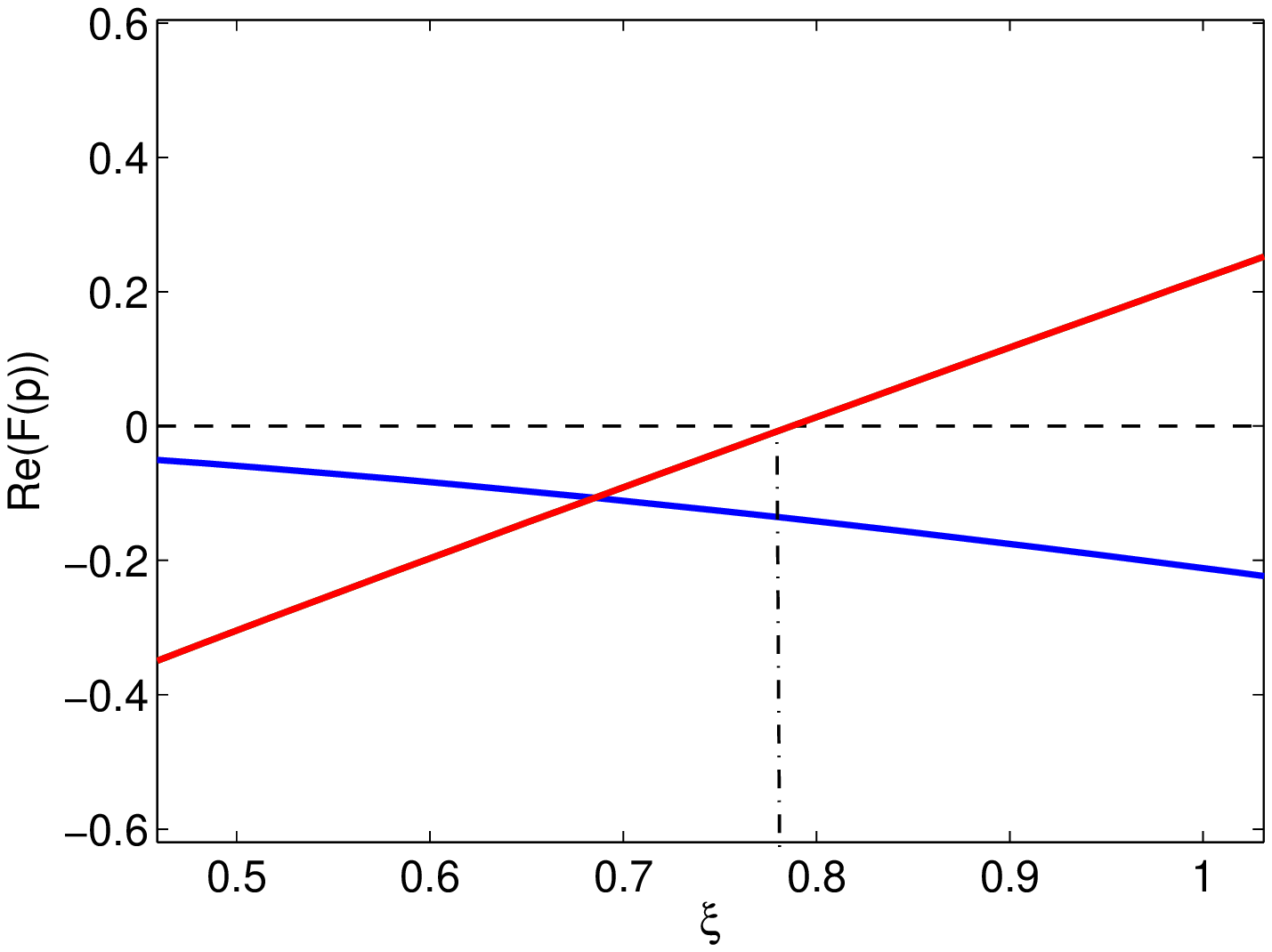}
\label{speedII}
}
}
\caption{The real parts of the branches of $F(\xi)$ for $\xi>0$ (here $\Phi_u=1$). Only three curves are seen because of the four saddle branches $p_i$ with $i=1,2,3,4$ (i.e. the solutions of (\ref{saddle:eq}) for $\xi\in \R$) two, which we take to be those with $i=2,3$, have $\mbox{Re}(F(p_2(\xi))=\mbox{Re}(F(p_3(\xi))$. \subref{speedII} shows a blow-up of the region within the small rectangle in \subref{speedI}. The plot shows that branches $p_2$ and $p_3$ have $\mbox{Re}(F(p_2(\xi^*))=\mbox{Re}(F(p_3(\xi^*))=0$ at a positive $\xi^*\approx 0.787$.}
\label{speedsearch}
\end{figure}

We observe that the condition (\ref{saddle:eq}) corresponds to $p$ being a saddle point of $F$. For real $\xi$, two of the four branches $p(\xi)$ in (\ref{saddle:eq}) have the same real 
part while their imaginary parts have opposite signs; we term these $p_2$ and $p_3$. 
They satisfy $\mbox{Re}(p_{2,3})<0$ for $\xi>0$ and $\mbox{Re}(p_{2,3})>0$ for $\xi<0$.
These branches give a pair of complex conjugate $F(p)$-branches; $F_{2,3}:=F(p_{2,3})$. 
The remaining branches, $p_1$ and $p_4$ are real for real $\xi$, and satisfy $F(p_1)$, $F(p_4)\in\R$ with $F(p_1)$, $F(p_4)\leq 0$. 
The real parts of the branches $F(p)$ for positive real $\xi$ are shown in Figure~\ref{speedsearch}. The branch points for the solutions of (\ref{saddle:eq}) are $p=-i/\sqrt{3}$, 
$\xi=i\Phi_u 9/(8\sqrt{3})$ and $p=-i/\sqrt{3}$, $\xi=-i\Phi_u 9/(8\sqrt{3})$, i.e. they occur at imaginary values of $\xi$, so the four branches 
described can be continued to $\xi<0$ for real $\xi$ through $\xi=0$. In particular, for negative $\xi$ the figures are symmetric according 
to $\mbox{Re}(F(p(\xi)))=\mbox{Re}(F(p(-\xi)))$ and 
$\mbox{Im}(F(p(\xi)))=-\mbox{Im}(F(p(-\xi)))$, thus giving complex-conjugate $F$'s. More details on the asymptotic behaviour of $p$- and $F$-branches appear in the Section~\ref{sec:expdecay}.

We assume for the moment that the initial perturbation decays faster that exponentially. For such initial data we expect that the envelope solution (\ref{clairaut:singular})-(\ref{saddle:eq}) dominates over the exponentially decaying ones (\ref{clairaut:general}). We use the neither-growth-nor-decay convention, $\mbox{Re}(F(p))=0$ and $\xi\in \R$, to locate the front and thus determine its speed and the rate of exponential decay ahead of it. 

We observe that $p_2$ and $p_3$ satisfy 
\[
F(p(\xi)) = -\Phi_u+O(\xi^{\frac{2}{3}})\quad \mbox{as}\ \xi \to 0^+ \,
\]
and 
\beqs
\mbox{Re}\left(F(p)\right)\sim \xi \quad \mbox{as}\ \xi \to +\infty \,, 
\eeqs
hence there exists a $\xi^*>0$ such that $\mbox{Re}(F(p_{2,3}(\xi^*)))=0$; see Figure~\ref{speedsearch} 
for $\Phi_u=1$ (observe that these figures represent the general case by scaling $\xi\to \Phi_u \xi$). The branches $1$ and $4$ give exponentially 
growing behaviours as $\xi\to +\infty$ and can be excluded.

At $\xi=\xi^*$ the absolute value of the exponential term in (\ref{wkbj:ansatz}) 
neither grows nor decays as $t\to +\infty$. In the case of linear selection (i.e. a pulled front) the front 
propagating to the right is located, to leading order as $t\to +\infty$, at 
$\xi=\xi^*$, i.e. in the original coordinate $x$ the front asymptotically propagates with constant speed $\xi^*$. 
The front propagating to the left is defined by the analogous argument for $\xi<0$ 
and has $x/t \sim -\xi^*$. The wave speed $\xi^*$ thus results from simultaneously solving  
(\ref{saddle:eq}) and the neither-growth-nor-decay condition 
\bequ\label{no:grow:cond}
\mbox{Re}\left(\frac{ 2p^2}{(1-p^2)^2 } - \frac{p^2}{1-p^2}
\right)=0\, . 
\eequ
Indeed, imposing $\xi\in\R$ in (\ref{saddle:eq}) and (\ref{no:grow:cond}) we get the 
following equations for $p$ corresponding to $\xi=\xi^*$
\bequ\label{speed:eqs}
\mbox{Re}(p)\mbox{Re}\left( \frac{2p}{(1-p^2)^2}\right)= \mbox{Re}\left(\frac{p^2}{1-p^2}\right)\, , 
\ \mbox{Im}\left(\frac{-2p}{(1-p^2)^2}\right)=0\,,
\eequ
which give four (two pairs of complex-conjugate) saddle points $p^*$ with
\bequ\label{the:saddle}
\mbox{Re}(p^*)\approx \pm 1.042 \quad \mbox{and} \quad \mbox{Im}(p^*)\approx \pm 0.834 \, .
\eequ
We observe that these values do not depend on $\Phi_u$, since $\xi^* \propto -\Phi_u$. 
Substituting (\ref{the:saddle}) into (\ref{saddle:eq}) gives the wave speed 
\bequ\label{speed}
\xi^*= \Phi_u \,\xi_0 \quad \mbox{with} \quad \xi_0\approx 0.787\,,
\eequ
(only the pair $p^*$ with $\mbox{Re}(p^*)>0$ gives $\xi^*>0$, the other pair giving negative $\xi^*$).

The coefficient $a$ in (\ref{wkbj:ansatz}) can be computed from (\ref{wkbj:amplitude}) 
using (\ref{clairaut}) to evaluate $f$ and its derivatives. 
Equation (\ref{wkbj:amplitude}) is a first order linear equation for $a$, whose rays are again given by (\ref{rays}). 
Since $a$ satisfies
\bequ
\frac{d a(x(t),t)}{dt} = c(x,t) \, a(x(t),t)\label{ampli:eq}
\eequ 
along rays, with 
\[
c(x,t):=  \Phi_u  \left(1+3 \left(\frac{\pa f}{ \pa x}\right)^2\right) \left( 1-\left(\frac{\pa f}{\pa x}\right)^2\right)^{-3}\frac{\pa^2 f}{\pa x^2}\,,
\]
then for $f$ given by (\ref{clairaut:general}) $a$ is constant along rays, while if $f$ is determined by (\ref{clairaut:singular})-(\ref{saddle:eq}) 
(so (\ref{rays}) with $x_0=0$ and $p$ arbitrary is an expansion fan emanating from the origin) then
$c(x,t)=-t/2$, and $a(x(t),t)= C t^{-1/2}$, where $C>0$ is constant on each ray, then
\[
a(x,t)= \frac{1}{t^{\frac{1}{2}}}\Omega\left(\frac{x}{t}\right) \quad \mbox{for some} \ \Omega\,.
\]
Since the rays are 
straight lines emerging from the origin, those that `initially' outrun the wavefront $x/t=\xi^*$ continue to do 
so for all $t$.


With (\ref{speed}) and the condition (\ref{no:grow:cond}), the solution at the leading edge of the front behaves as
\bequ\label{predict:decay}
v(x,t)\sim \frac{1}{t^\frac{1}{2}}\,\Omega\left(\frac{x}{t}\right) e^{-(x-\xi^* t)p^*}\,e^{-F(p^*)t}\, \quad t\to +\infty \quad \mbox{with} \ \frac{x}{t}=O(1)\,, 
\eequ
where
\begin{equation}\label{w:star}
 F(p^*)\approx \pm 1.1688\Phi_u\, i\,.
\end{equation}
Observe that (\ref{predict:decay}) decays exponentially like $e^{-\l^* |x|}$ as $|x| \to \infty$, with  
\bequ\label{univ:decay}
\l^*:=  \mbox{Re}(p^*) \approx 1.042 \,.
\eequ
As with Fisher's equation (cf. \cite{bramson}), the pre-exponential factor here does not remain valid 
when nonlinear effects are accounted for, but the argument for wave speed selection does.

\subsection{Exponentially decaying initial conditions}\label{sec:expdecay}
In this section we aim to clarify the front speed selection mechanism for initial conditions of the form (\ref{exp:ic}) with $\l\in \R$. 

If we consider initial perturbations of the form (\ref{exp:ic}) for (\ref{main:eq}), the linearised equation (\ref{lin:equ}) gives different possible behaviours 
far ahead of the front corresponding to (\ref{slowdecay}), each of which could give rise to a possible front 
location with speed 
\begin{equation}\label{speed:lambda}
\xi_f(\lambda)= \frac{\lambda}{\lambda^{2}-1}\,,
\end{equation}
that results by applying the neither-growth-nor-decay condition to locate the front. Under this scenario the initial data generate an exponentially dominant 
to the fast decay solution (\ref{clairaut:singular})-(\ref{saddle:eq}) behaviour in the tail. 
We find below that $\xi^*$ as in (\ref{speed}) is the front speed selected for any real $\l$: it is worth contrasting this with the linearised Fisher equation 
\[ 
\frac{\pa v}{\pa t}=\frac{\pa^2 v}{\pa x^2}+v\,,
\] 
which has separable solutions $v(x,t)=e^{(\l^2+1)t-\l x}$ and, in this case, the neither-growth-nor-decay condition gives $\xi_f(\l)=\l+1/\l$. This function grows
 unboundedly as $\l\to 0^+$ and attains its global minimum $\xi_f=2$ at $\lambda=1$. In this case the pulled-front wave speed is $2$ for 
$\l\geq 1$ (thus $\xi^*=2$ and is associated with the fastest decay rate supported by the linearised 
equation and with the minimal speed travelling wave) and $\xi_f(\l)$ ($>2$) for $\l<1$ (see e.g. \cite{frontsreview}, \cite{CK} and the references therein).

It is worth noticing that for (\ref{main:eq}) $\xi_f(\lambda)\leq 0$ for $\lambda<1$ and $\xi_f(\lambda)>0$ for $\lambda>1$, having an asymptote at $\lambda=1$, 
and hence having no global minimum. For most pattern forming systems exhibiting travelling wave fronts of the pulled type, the corresponding function $\xi_f(\lambda)$ has a global minimum for real $\lambda$, moreover, this minimum is $\xi^*$ and is attained at the associated exponential decay rate $\lambda^*$, see \cite{frontsreview} (e.g. on page 49). 
In general, the linearly selected wave speed is 
$\xi^*$ if $\lambda\geq\lambda^*$ and $\xi_f(\lambda)$ if $\lambda<\lambda^*$. In this context, (\ref{main:eq}) gives an 
exception to this {\it rule}.

In what follows we aim to discern which of the exponential behaviours is selected in the limit $t\to +\infty$ with $x/t=O(1)$, assuming that the initial exponential behaviour pertains in the tail.

\paragraph{The limit $\lambda \to 1$.} The current limit is instructive both as the borderline case and because the problem becomes more tractable analytically than for general $\lambda$. 
Starting from the linearised problem 
\bequ
\frac{\pa v}{\pa t} =-\frac{\pa^2 v}{\pa x^2} + \frac{\pa^3 v}{\pa x^2 \pa t}\,, 
\label{lin:phi:1}
\eequ 
we set $v=e^{-x}w$ to give
\bequ
2\frac{\pa^2 w}{\pa x\pa t } + w = -\frac{\pa^2 w}{\pa x^2}+ 2 \frac{\pa w}{\pa x}+ \frac{\pa^3 w}{\pa x^2 \pa t}\,,
\label{lin:w}
\eequ
where we have gathered on the left-hand side the terms that will end up being dominant. We now consider the two sets of initial data 
\bequ
\mbox{at} \quad t=0 \quad w=e^{-\e x} \quad \mbox{or}\quad \mbox{at} \quad t=0 \quad w=e^{\e x}
\label{Lin:data}
\eequ
with $0<\e\ll 1$, the former corresponding to the limit $\lambda\to 1^+$ and the latter corresponding to the limit $\lambda \to 1^-$. 
Appropriate scalings are then
\[
x=\frac{X}{\e}\,,\quad t=\e T 
\]
which furnish the leading-order problems (suppressing the subscript on $w_0$)
\beq
2\frac{\pa^2 w}{\pa X\pa T} = -w \quad \mbox{in}\quad X>0 \label{w:eq:rescal}\\
\mbox{at} \quad X=0 \quad w= 1\,, \label{w:eq:bc}\\
\mbox{at} \quad T=0 \quad w=e^{-X}\quad\mbox{or} \quad \mbox{at}\quad T=0 \quad w=e^X \label{w:eq:ic}
\eeq
where we have applied, somewhat arbitrarily, the condition (\ref{w:eq:bc}) to ensure that the solution is not simply separable (other such conditions would be 
equally instructive). 

The JWKB ansatz
\[
w\sim e^{-f(X,T)}
\]
yields the dominant balance
\bequ
2\frac{\pa f}{\pa T}\frac{\pa f}{\pa X}=-1
\label{phi:eq}
\eequ
for which 
\bequ
\frac{dX}{dT}=-\frac{1}{2p^2}
\label{X:eq}
\eequ
holds along rays, where $p\equiv\pa f/\pa X$. Particular solutions both to (\ref{w:eq:rescal}) and (\ref{phi:eq}) corresponding to the initial data 
(\ref{w:eq:ic}) have
\bequ
w=\exp(-X+\frac{T}{2})\,, \quad w=\exp(X-\frac{T}{2})
\label{w:eq:sol}
\eequ
both of which have from (\ref{X:eq}) that 
\bequ
\frac{d X}{dT}=-\frac{1}{2}
\label{X:eq:2}
\eequ
along rays.

Now seeking a solution to (\ref{phi:eq}) of the form 
\bequ
f=T F(\eta)\,, \quad \eta=\frac{X}{T}\,,
\label{phi:selfsim}
\eequ
(observe that $\eta=\e^2 \xi$ in (\ref{clairaut})) we find the general solutions 
\[
F=c\eta -\frac{1}{2c}
\]
for an arbitrary constant $c$, the solutions (\ref{w:eq:sol}) each being special cases of this, namely 
\bequ
F=\eta-\frac{1}{2}\,, \quad F=\frac{1}{2}-\eta\,,
\label{special:w:eq:sol}
\eequ
and singular (envelope) solutions
\bequ
F=\pm i\sqrt{2\eta}\,,
\label{envelope}
\eequ
it being a virtue of the current limit that these take a simple explicit form. The turning points at which 
(\ref{special:w:eq:sol}) and (\ref{envelope}) coincide are at $\eta=-1/2$; this of course lies outside 
the range in which (\ref{w:eq:rescal}) is being taken to hold, but will nevertheless play an revealing role in what follows.  

Before proceeding  further with the analysis corresponding to (\ref{w:eq:ic}), we consider the case in which 
(\ref{w:eq:rescal})-(\ref{w:eq:bc}) are subject to 
\[
\mbox{at}\quad T=0 \quad w=1\,,
\]
corresponding to $\lambda=1$, for which the solution takes the self-similar form 
\bequ
w=W(\zeta)\, , \quad \zeta=X T 
\label{w:selfsim}
\eequ
with 
\bequ
2\left(\zeta\frac{d^2 W}{d\zeta^2}+\frac{d W}{d\zeta}\right)=-W\,;
\label{W:eqII}
\eequ
a local analysis about $\zeta=0$, for which the left-hand side dominates, leads to eigenmodes proportional to 
$\zeta^0$ and $\zeta^0 \ln \zeta$ of which the latter needs to be rejected: this reveals that (\ref{W:eqII}) subject to 
\[
\mbox{at} \quad \zeta=0 \quad W=1
\]
is an initial value problem, whereby both signs in (\ref{envelope}) represent 
possible behaviours as $\zeta \to +\infty$, the magnitudes of the associated contributions being determined by the initial value problem from 
$\zeta=0$. Thus it is not possible to enforce any aspect of the behaviour of $W$ as $\zeta\to +\infty$, a result that will be instructive in 
what follows.

\begin{figure}
\centerline{
  \includegraphics[width=0.62\textwidth,height=.44\textwidth]{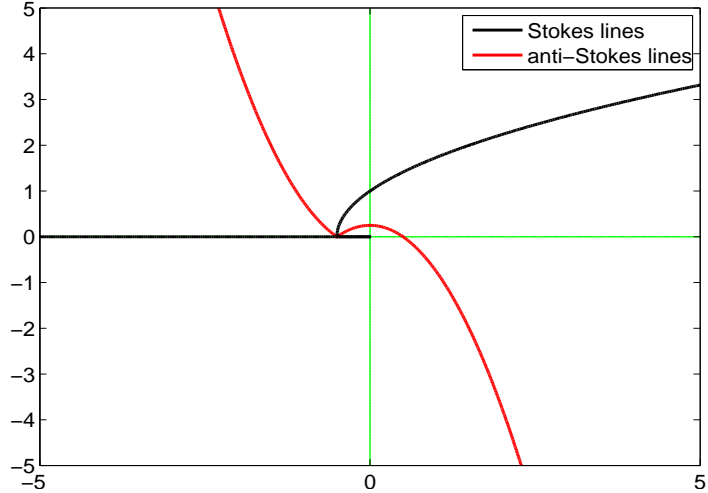}
}
\caption{Stokes and anti-Stokes lines for $\l=1$.}
\label{near1stokeslines}
\end{figure}

Returning now to (\ref{special:w:eq:sol}), the Stokes lines for these contributions relative to (\ref{envelope}) are obtained in the usual way by 
setting the imaginary parts equal (i.e. $\mbox{Im}(\eta)=\mbox{Im}(\pm i\sqrt{2\eta})$ for the contribution $F=\eta-1/2$), whereby in both cases
\bequ
\sqrt{2\rho}\sin(\theta/2)=\pm 1\,,
\label{stokesI}
\eequ
where we have set $\eta=\rho e^{i \theta}$. Similarly, the anti-Stokes (equal real parts) have
\bequ
\rho \cos \theta -\frac{1}{2}=\pm \sqrt{2\rho}\sin (\theta/2)\,.
\label{stokesII}
\eequ
The associated complex plane picture is shown in Figure~\ref{near1stokeslines}.
For $w=e^X$ at $T=0$ ($\lambda \to 1^-$) the exponential contribution is thus dominant everywhere along the Stokes line, but since its asymptotic 
series truncates, it can turn nothing on. However for more general initial data such as $w=X^\gamma e^X$, the associate series may diverge and turn on
 an envelope contribution; this contribution is associated purely with the initial data and hence has no knowledge of the boundary condition 
on $X=0$. There is accordingly a distinct envelope contribution present everywhere which can be thought of as enabling this boundary 
condition to be satisfied (a small-time analysis can be used to clarify this issue).

For $w=e^{-X}$ at $T=0$ ($\lambda \to 1^+$), the (\ref{special:w:eq:sol}) contribution is subdominant everywhere along the Stokes line and 
requiring it to be present in the bulk of the complex plane as $\mbox{Re}(\eta) \to +\infty$ implies that this contribution is switched {\it off} across the 
Stokes line on moving towards the positive real axis. This contribution is accordingly absent altogether on the real axis. 
(The width of the Stokes line presumably grows as $\rho^{1/2}$ as $\rho \to +\infty$, but since $\ln \eta$ also grows on as 
$\rho^{1/2}$ on the Stokes line this does not cause it to impinge on the real axis, unlike Stokes lines that 
run parallel to the real axis, cf. King \cite{Kingkyoto}.)

In summary, for initial conditions $w=e^X$ both exponential and envelope contributions are present everywhere, but as we show 
above the former plays no role in wave speed selection. For initial data $w=e^{-X}$, the exponential contribution is entirely absent 
from the far-field on the real line, being switched {\it off} across the Stokes lines. This somewhat surprising result is crucial to 
the selection of the wave speed.

The above analysis is in effect one for large $T$. It remains to clarify how, for small $T$, the $\exp(-X)$ contribution disappears at infinity. 
From (\ref{w:eq:rescal}) a naive small time expansion reads
\bequ
w\sim e^{-X}-\frac{1}{2}(1-e^{-X})T +\frac{1}{8}(X-1+e^{-X})T^2\,, 
\label{small:T}
\eequ
which is clearly non-uniform for large $X$ with the outer region having $X=O(1/T)$ and 
\bequ
w\sim T \Psi (\zeta)\, ,\ \zeta=XT
\label{small:T:2}
\eequ
with
\[
\zeta \frac{d^2 \Psi}{d\zeta^2}+2\frac{d\Psi}{d\zeta}=-\frac{1}{2}\Psi
\]
subject, on matching to (\ref{small:T}), to the initial data
\[
\Psi \sim -\frac{1}{2} +\frac{1}{8} \zeta \quad \mbox{as}\ \zeta\to 0^+\,.
\]
The $\exp(-X)$ terms in (\ref{small:T}) (which can of course be summed in the form $\exp(-(X-T/2))$) are exponentially 
subdominant as $T\to 0^+$ with $\zeta=O(1)$ and are turned off across the associated Stokes lines by the divergent series whose leading term is 
(\ref{small:T:2}); to this order of calculation, the Stokes lines coincide with the positive real axis: 
the situation can be clarified by noting that they are in fact described by (\ref{stokesI}) uniformly in time and that 
$\rho=O(1/T^2)$ for $\zeta=O(1)$, so that (\ref{W:eqII}) implies that $\mbox{Im}(\zeta)=O(T)$ 
(hence $\mbox{Im}(X)=O(1)$) as $T\to 0$, this illustrates how the $\exp(-X)$ term is cleared off the $X$ axis, 
being absent from a region about the axis that grows with $T$.

Two other comments are in order. Firstly, the turning point location $X=-T/2$ not surprisingly corresponds 
to the characteristic velocity of the solutions $f=\pm (X-T/2)$, the local behaviour of the turning point being 
described by the heat equation. Secondly, it is already clear from (\ref{w:eq:rescal}) being of first order in $X$ 
that the behaviour as $X\to +\infty$ cannot be imposed as a boundary condition, and the above argument implies 
that the presence of the exponential $\exp(-\lambda x)$ does not follow from the initial condition.

\begin{figure}[htb]
\centering
\mbox{
\subfigure[Real part of $p(\xi)$, $\xi\in\R$]
{
\includegraphics[width=0.46\textwidth,height=.35\textwidth]{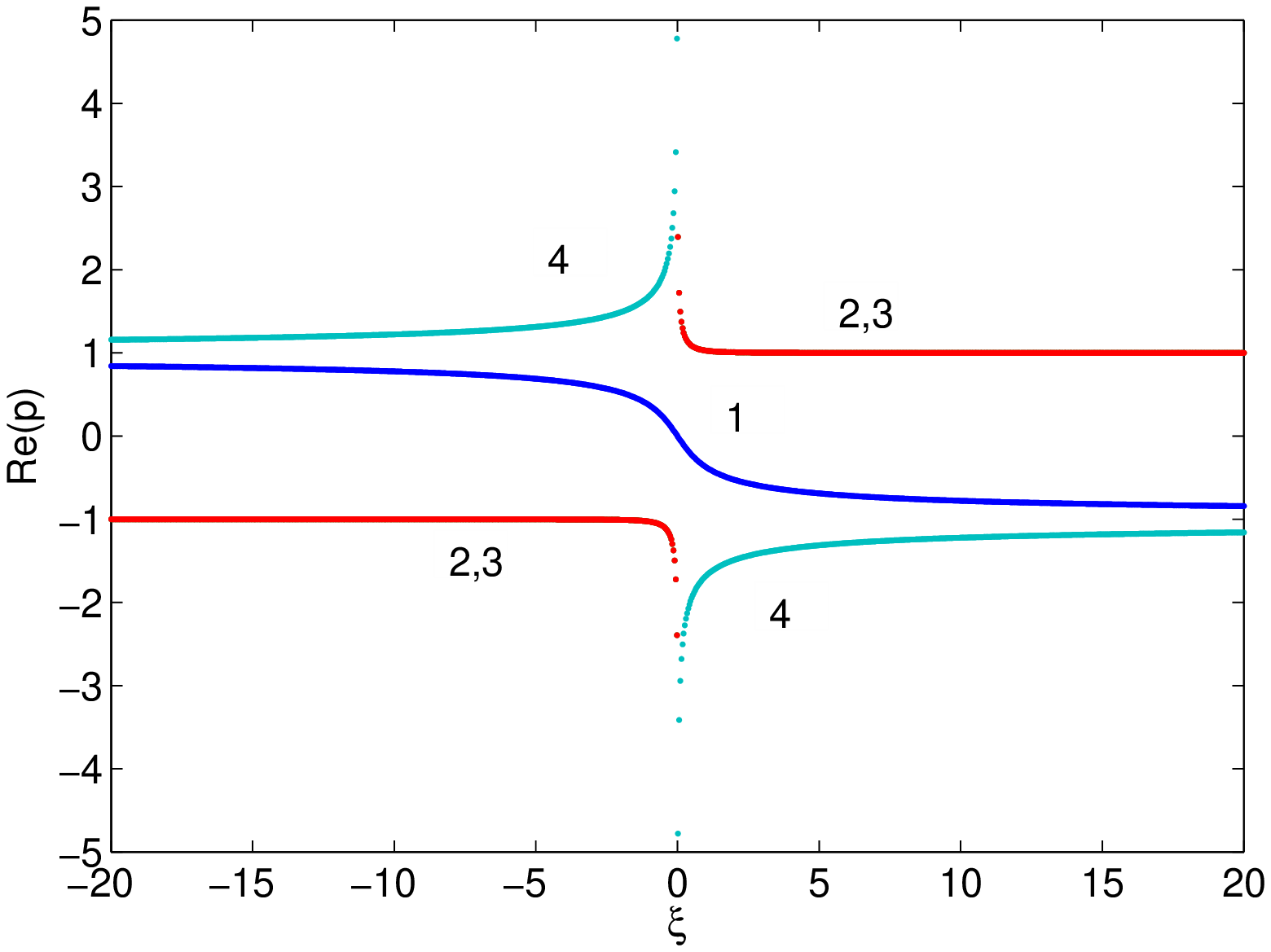}
\label{realq}}
}
\mbox{
\subfigure[Imaginary part of $p(\xi)$, $\xi\in\R$] 
{
\includegraphics[width=0.46\textwidth,height=.35\textwidth]{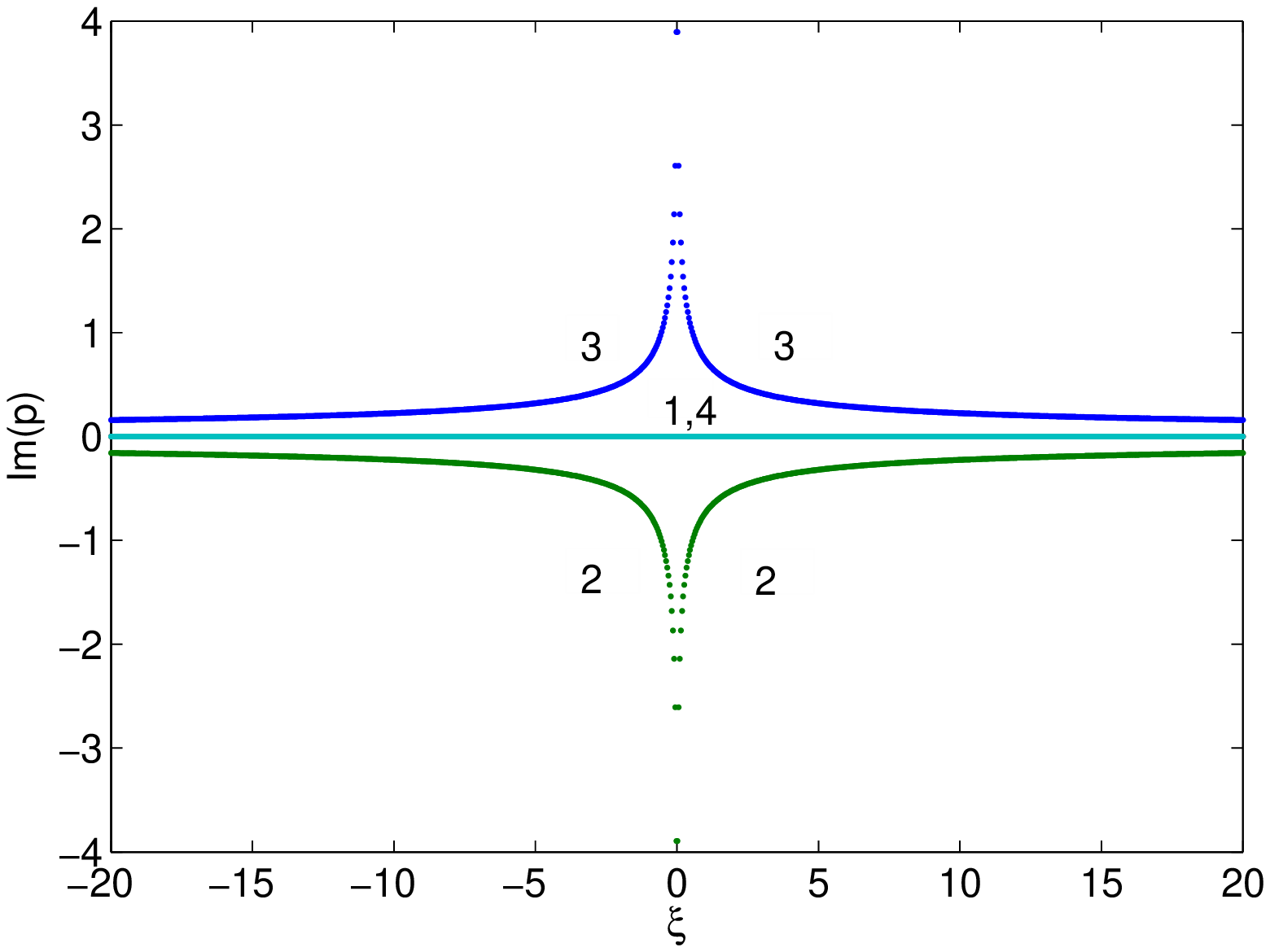}
\label{imq}}
}
\caption{Real and imaginary part of the four $p$-branches as a function of real $\xi$.
}
\label{realimagq}
\end{figure}
\begin{figure}[htb]
\centering
\mbox{
\subfigure[Real part of $F(q(\xi))$]
{
\includegraphics[width=0.46\textwidth,height=.35\textwidth]{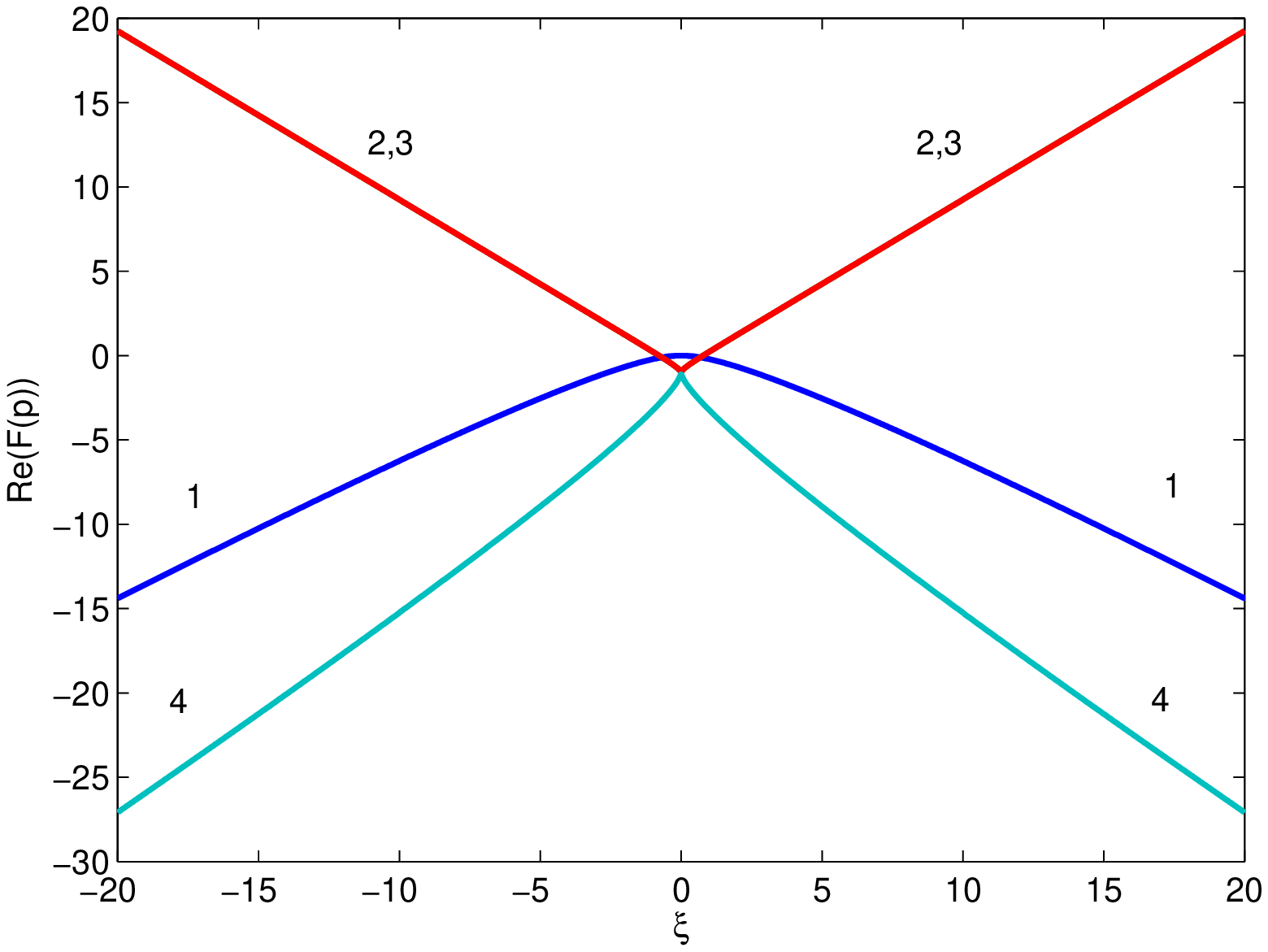}
\label{realw}}
}
\mbox{
\subfigure[Imaginary part of $F(p(\xi))$] 
{
\includegraphics[width=0.46\textwidth,height=.35\textwidth]{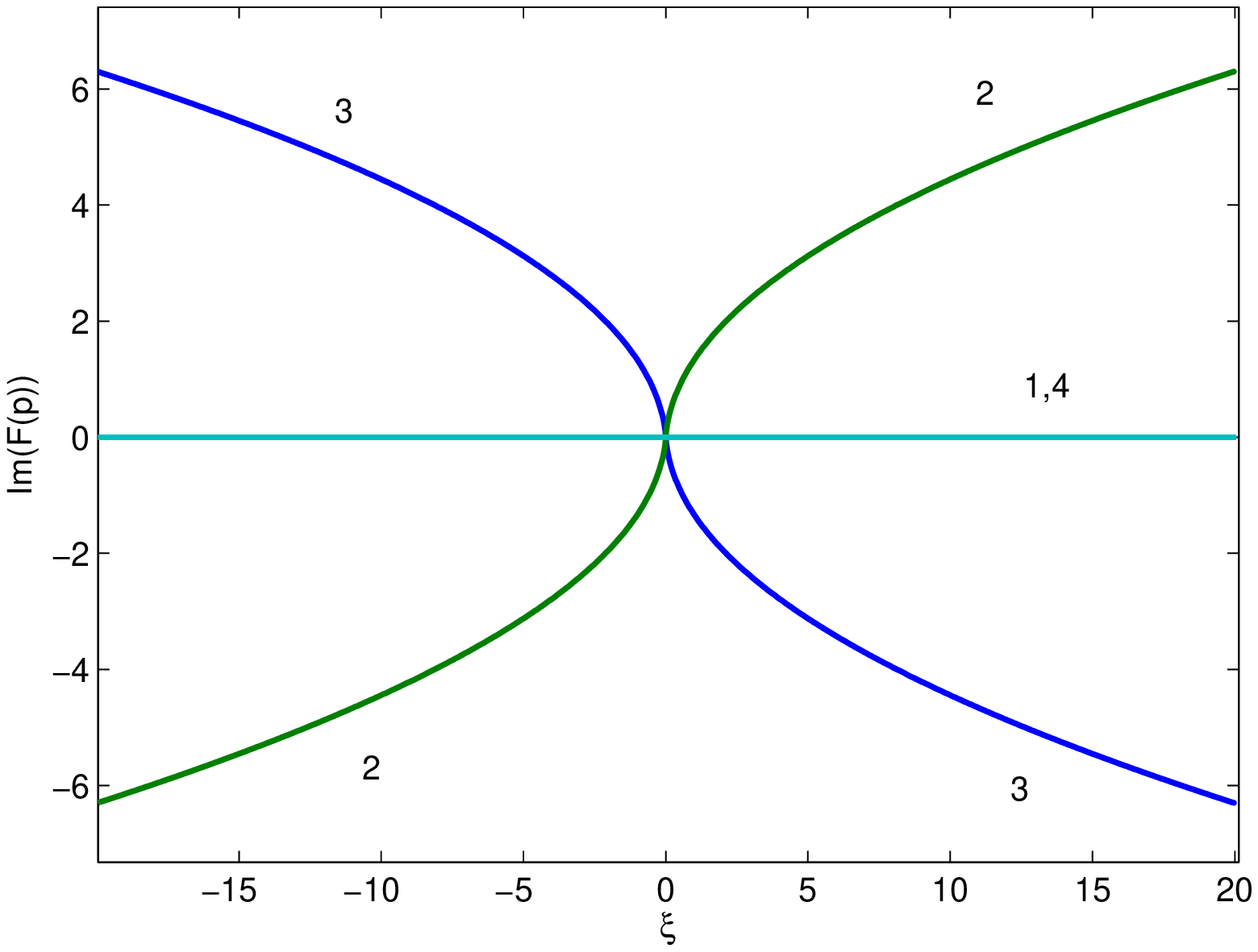}
\label{imagw}}
}
\caption{Real and imaginary part of the four $F(p)$-branches as a function of real $\xi$.
}
\label{realimagw}
\end{figure}

\paragraph{Real $\l$.} We analyse the Stokes lines associated to the exponential contributions given by the (\ref{clairaut:general}) and (\ref{clairaut:singular})-(\ref{saddle:eq}) for a fixed value of $\l$. We first need to recall and further analyse the saddle branches defined by (\ref{clairaut:singular})-(\ref{saddle:eq}). The real and imaginary parts of the branches $p$ for real $\xi$ are shown in Figure~\ref{realimagq}, and the real and imaginary parts of the branches $F(p)$ for real $\xi$ are shown in Figure~\ref{realimagw}.

\begin{figure}[htb]
\centering
\mbox{
\subfigure[Real part of $p(i\zeta)$, $\zeta\in\R$]
{
\includegraphics[width=0.46\textwidth,height=.35\textwidth]{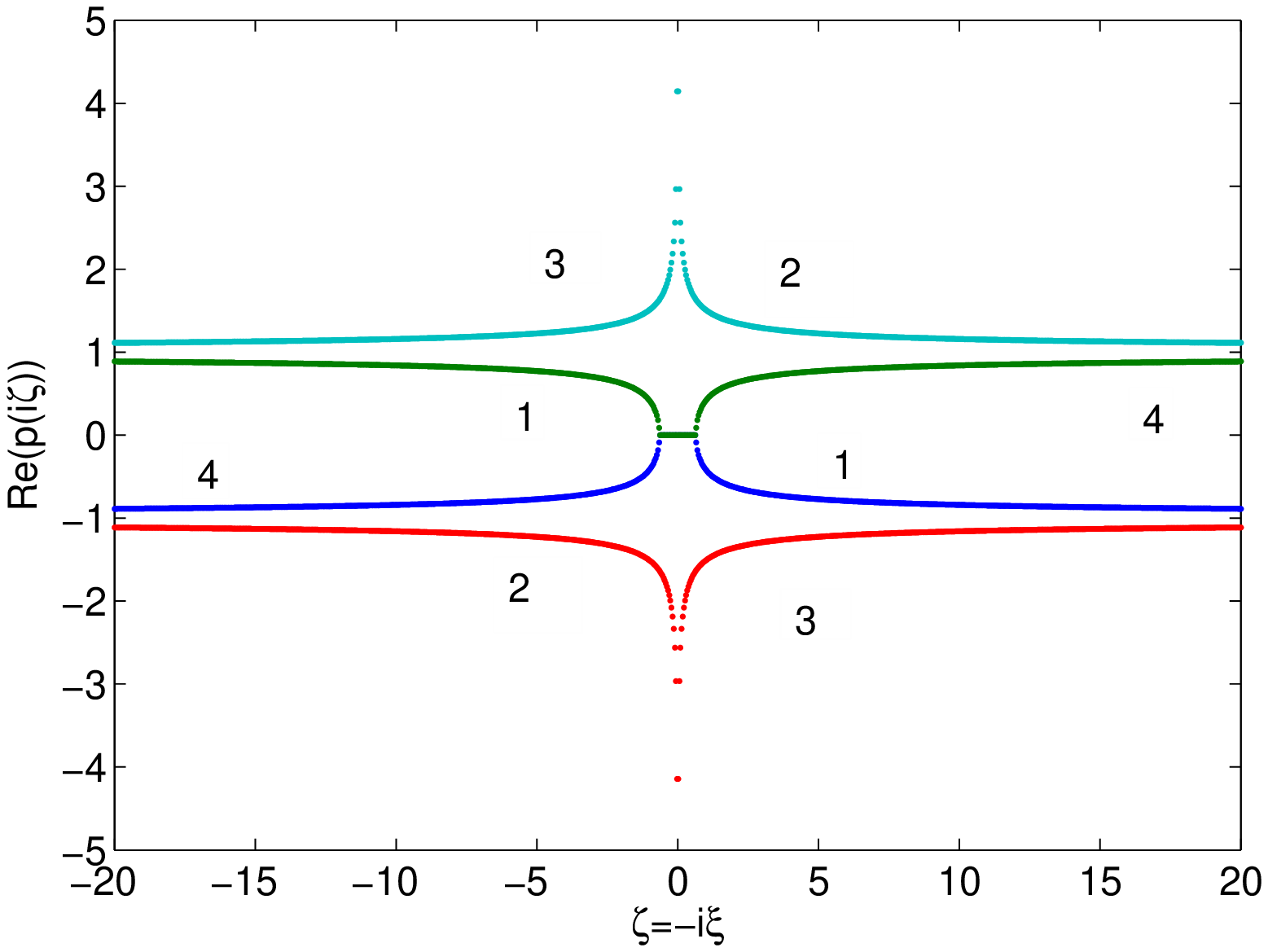}
\label{realqimagxi}}
}
\mbox{
\subfigure[Imaginary part of $p(i\zeta)$, $\zeta\in\R$] 
{
\includegraphics[width=0.46\textwidth,height=.35\textwidth]{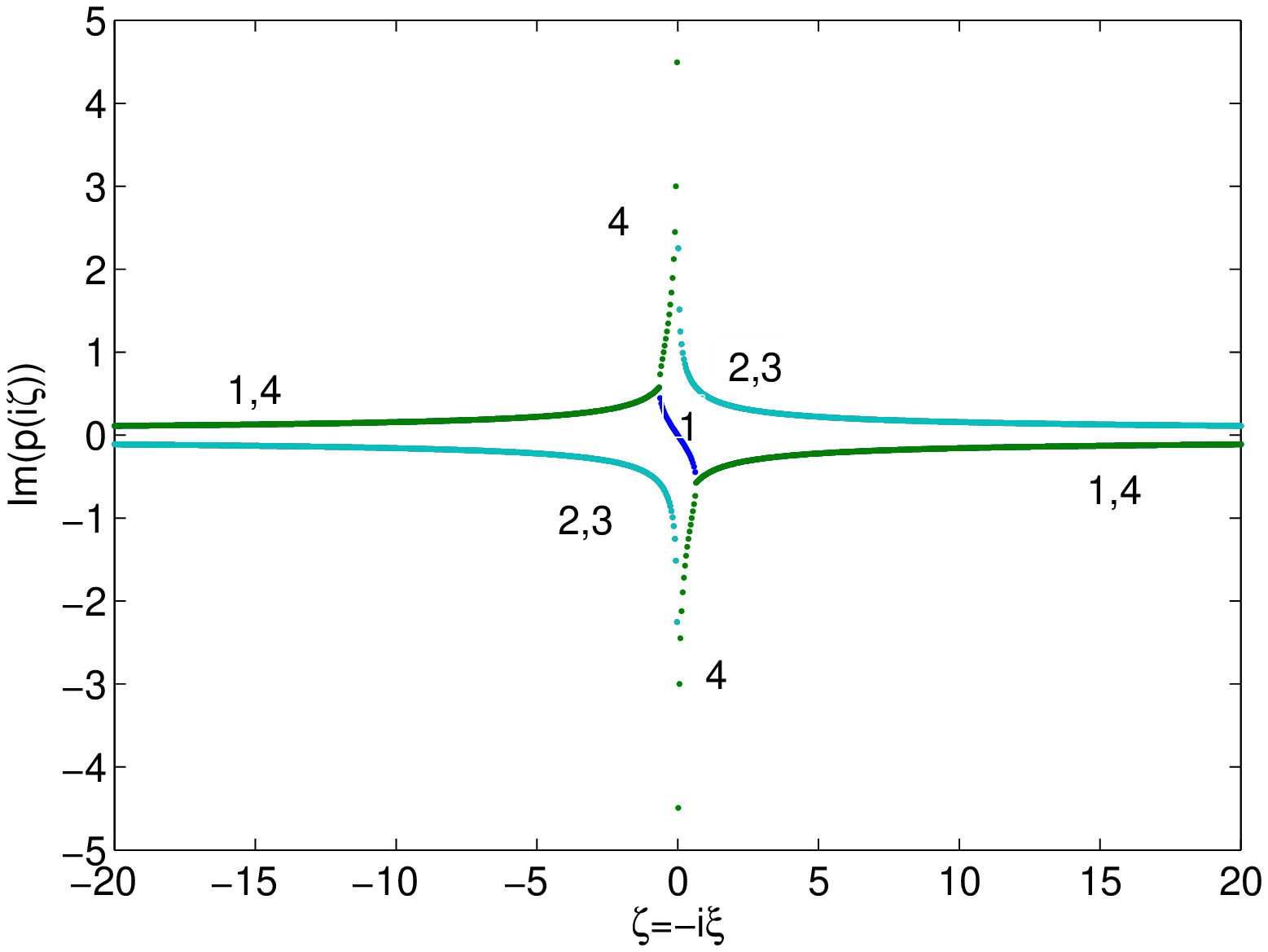}
\label{imqimagxi}}
}
\caption{Real and imaginary part of the four $p$-branches as a function of $\zeta=-i\xi \in \R$. 
}
\label{realimagqimagxi}
\end{figure}

We recall that the branch points for the solutions of (\ref{clairaut:singular})-(\ref{saddle:eq}) are 
\bequ\label{branch:points}
p=\pm \frac{i}{\sqrt{3}} \quad \mbox{at}\quad \xi= \pm i \frac{9}{8\sqrt{3}}\,.
\eequ
At these points the branches $p_1$ and $p_4$ coincide and swap identity along the branch cuts $\xi=i\zeta$ with 
$\zeta\in( -9/(8\sqrt{3}), 9/(8\sqrt{3}))$, 
see Figure~\ref{realimagqimagxi}. Observe that setting $p=i\tau $ gives $\xi=-2i\tau  /(1+\tau^2)^2$ 
thus the sign of the real parts of the $p$-branches can only change across 
this graph, which lies on the imaginary axis in the $\xi$-plane.

We now give the asymptotic expansions of these for small and large $|\xi|$ with $\xi\in \R$. We have
\[
p_1(\xi)= -\frac{1}{2}\xi + \frac{1}{4}\xi^3+ O(\xi^4)\quad \mbox{as} \ \xi \to 0^\pm\,,
\]
and
\beqs
\lefteqn{
p_2(\xi) = \frac{1}{2^{\frac{2}{3}}}(1-\sqrt{3}i)\xi^{-\frac{1}{3}}+O(\xi^{\frac{1}{3}})
\quad \mbox{as} \ \xi \to 0^\pm \,, 
}
\\
\lefteqn{
p_3(\xi) = \frac{1}{2^{\frac{2}{3}}}(1+\sqrt{3}i)\xi^{-\frac{1}{3}}+O(\xi^{\frac{1}{3}})
\quad \mbox{as} \ \xi \to 0^\pm \,, 
}
\\
\lefteqn{
p_4(\xi)= - 2^\frac{1}{3} \xi^{-\frac{1}{3}}+O(\xi^{\frac{1}{3}})  \quad \mbox{as} \ \xi \to 0^\pm\,.
}
\eeqs
Substituting these expressions into $F(p)$($=p(1+p^2)\xi/2$) we get, denoting $F_j(\xi)=F(p_j(\xi))$ for all $j=1,2,3$, that
\beq
\lefteqn{
F_1(\xi)= -\frac{1}{4}\xi^2 - \frac{1}{ 16}+O(\xi^6) \to 0 \quad \mbox{as}\ \xi \to 0^\pm\,,
}
\label{w:zero:1}\\
\lefteqn{
F_2(\xi)= -2  + \frac{1}{2^{\frac{2}{3}}}(1-\sqrt{3}i)\xi^{\frac{2}{3} }+ O(\xi)\quad \mbox{as}\ \xi \to 0^\pm \,,
}\label{w:zero:2}\\
\lefteqn{
F_3(\xi)= -2  +  \frac{1}{ 2^{\frac{2}{3}}}(1+\sqrt{3}i)\xi^{\frac{2}{3} } +O(\xi)\quad \mbox{as}\ \xi \to 0^\pm \,,
}
\label{w:zero:3}\\
\lefteqn{
F_4(\xi)= -2 -2^{\frac{1}{3}}\xi^{\frac{2}{3}}+O(\xi)\quad \mbox{as}\ \xi \to 0^\pm \,.
}
\label{w:zero:4}
\eeq
The real and imaginary parts of the $p$-branches and the corresponding $F(p)$ branches are shown 
in figures \ref{realimagq} and \ref{realimagw} respectively. They satisfy the following 
asymptotic behaviour as $\xi\to +\infty$
\beqs
p_1(\xi)=-1+ \frac{1}{\sqrt{2}}\xi^{-\frac{1}{2}}+O(\xi^{-\frac{3}{2}}) \quad\mbox{as} \ \xi \to  +\infty\,,
\\
p_2(\xi)=1-i \frac{1}{\sqrt{2}} \xi^{-\frac{1}{2}}+O(\xi^{-\frac{3}{2}}) \quad\mbox{as} \ \xi \to  +\infty\,,
\\
p_3(\xi)=1+i\frac{1}{\sqrt{2}} \xi^{-\frac{1}{2}}+O(\xi^{-\frac{3}{2}}) \quad\mbox{as} \ \xi \to  +\infty\,,
\\
p_4(\xi)=-1- \frac{1}{\sqrt{2}}\xi^{-\frac{1}{2}}+O(\xi^{-\frac{3}{2}}) \quad\mbox{as} \ \xi \to  +\infty\,,
\eeqs
giving
\beq
F(p_1(\xi))=-\xi - 2 \frac{1}{\sqrt{2}}\xi^\frac{1}{2}+O(1) \quad \mbox{as}\ \xi \to +\infty\,,
\label{w:infty:1}
\\
F(p_2(\xi))=\xi - 2i \frac{1}{\sqrt{2}} \xi^\frac{1}{2}+O(1) \quad \mbox{as}\ \xi \to +\infty\,,
\label{w:infty:2}
\\
F(p_3(\xi))=\xi + 2i  \frac{1}{\sqrt{2}} \xi^\frac{1}{2}+O(1) \quad \mbox{as}\ \xi \to +\infty\,,
\label{w:infty:3}
\\
F(p_4(\xi))=-\xi + 2 \frac{1}{\sqrt{2}} \xi^\frac{1}{2}+O(1) \quad \mbox{as}\ \xi \to +\infty\,.
\label{w:infty:4}
\eeq

We can now compute the associated Stokes lines and anti-Stokes lines.
We recall that on the Stokes lines the exponential behaviour of each $F$-branch in comparison to the exponent $\lambda$ and, respectively, 
in comparison to each other, is either maximal or minimal depending on the region of the complex $\xi$-plane. Namely, we compute numerically the contours where for $\xi\in \C$ and $\mbox{Re}(\xi)\geq 0$ (right-hand front) 
\[
\mbox{Im}(F_i(\xi))=\lambda \mbox{Im}(\xi) \quad \mbox{for} \quad i=1,2,3,4
\]
and
\[
\mbox{Im}(F_i(\xi)) = \mbox{Im}(F_j(\xi)) \quad \mbox{for} \quad i,j=1,2,3,4 \quad i\neq j\,.
\]
and those where
\[
\mbox{Re}(F_i(\xi))=\lambda \mbox{Re}(\xi) +\mbox{Re}(\l^2/(1-\l^2))\quad \mbox{for} \quad i=1,2,3,4
\]
and
\[
\mbox{Re}(F_i(\xi)) = \mbox{Re}(F_j(\xi)) \quad \mbox{for} \quad i,j=1,2,3,4 \quad i\neq j\,.
\]
The relevant branches $F_j$ are the ones that give exponential decay as $\xi\to +\infty$, i.e. $F_2$ and $F_3$.

We can make some further observations by setting $p=p_r+ip_i$ and $\xi=\xi_r+i\xi_i$ and writing
\beqs
F(\xi,p)= (p_r \xi_r -p_i\xi_i)(1+p_r^2-p_i^2) -2p_r p_i(p_r\xi_i+p_i\xi_r) \qquad \\\qquad\qquad\qquad\qquad 
+\,i\left( (p_r\xi_i+p_i\xi_r)(1+p_r^2-p_i^2) +2p_rp_i(p_r\xi_r-p_i\xi_i)\right)\,.
\eeqs
It is easy to see (with the aid of Figure~\ref{realimagqimagxi}) that the imaginary axis is an anti-Stokes line for branches 
$F_2$ and $F_3$, and for $F_1$ and $F_4$ away from the branch cut. Also the real axis is an anti-Stokes line for branches 
$F_2$ and $F_3$. Similarly, one can easily see that the real and imaginary axes are not Stokes lines for neither of the pairs 
$F_2$ and $F_3$, and $F_1$ and $F_4$. 
A computation also shows that the turning points (where an anti-Stokes line 
crosses the corresponding 
Stokes line) are the branch points (\ref{branch:points}) (setting $F_k=F_j$ for $k\neq j$ with, by the second equation in 
(\ref{clairaut:singular}), $p_k(1-(p_j)^2)^2= p_j(1-(p_k)^2)^2$ necessarily gives $p_k=p_j$). 

The turning points for the 
anti-Stokes and Stokes lines $F_k$ 
versus the $\l$ contribution are calculated by setting $F_k=\l\xi+\l^2/(1-\l^2)$ for a given $\l\in\R$ with 
\bequ\label{xi:pi}
\xi=-\frac{2p_k}{(1-(p_k)^2)^2}\,,
\eequ 
this gives the equation for $p_k$ (if $p_k\neq \l$)
\bequ\label{pol:turning}
(p_k)^2 + 2\l p_k + 1 = 0 
\eequ    
thus
\[
p_k= -\l \pm \sqrt{\l^2-1}
\]
giving the turning points upon substituting into (\ref{xi:pi}).
Obviously, for $\l>1$ 
these gives two positive real turning points. For the branches $F_1$ and $F_4$ the real axis is a Stokes line, since they are real for real $\xi$, 
see Figure~\ref{imagw}. The real positive turning points correspond to these branches when $\l>1$. 
For $\l<1$ we obtain the turning points
\[
\xi = \frac{1}{2}\left( \frac{\l}{\l^2-1}\pm i\frac{1}{\sqrt{1-\l^2}} \right)
\]
that have negative real part.

These observations and the expansions of $F$ in (\ref{w:zero:1})-(\ref{w:zero:4}) and (\ref{w:infty:1})-(\ref{w:infty:4}) 
allow one to identify the Stokes lines and anti-Stokes lines obtained 
numerically. 

\begin{figure}
\centerline{
\includegraphics[width=0.45\textwidth,height=.35\textwidth]{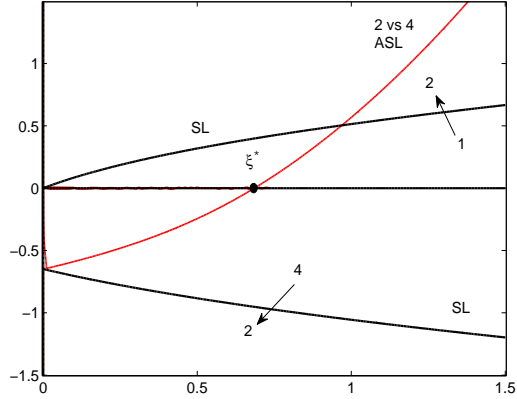}
}
\caption{The $\mbox{Im}(F_2)=\mbox{Im}(F_j)$ ($j\neq 2$) Stokes lines and $\mbox{Re}(F_2)=\mbox{Re}(F_j)$ ($j\neq 2$) anti-Stokes lines. 
The front location corresponds to where the anti-Stokes line crosses the axis. The exponentially large contributions $1$ and $4$ 
are absent so the only Stokes lines that are present are in fact inactive and the oscillatory decaying exponential $2$ and $3$ are 
present everywhere in the far-field of the right-hand $\xi$ plane.}
\label{linesF2}
\end{figure}

Figure~\ref{linesF2} shows the Stokes and anti-Stokes lines for the branch $F_2$ comparing it to $F_1$ and to $F_4$.
 The Stokes line with $F_1$ would turn $F_2$ on but $F_1$ is absent, thus this Stokes line is inactive. The same 
argument applies to the Stokes line with $F_4$. Exchange of dominance does not occur in this right-half plane. 

\begin{figure}
\centerline{
\includegraphics[width=0.45\textwidth,height=.35\textwidth]{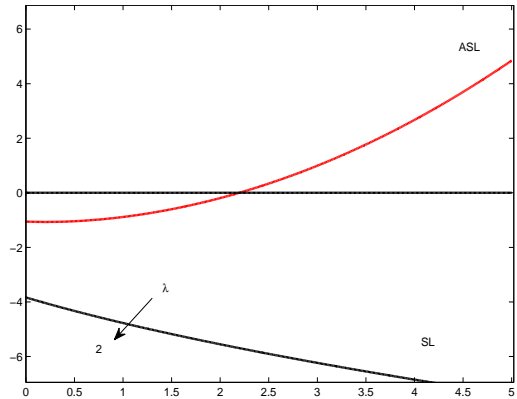}
}
\caption{$\lambda<1$: The $\mbox{Im}(F_2)=\l \mbox{Im}(\xi)$ Stokes line and $\mbox{Re}(F_2)=\l\mbox{Re}(\xi)+\l/(1-\l^2)$ anti-Stokes lines for $\l=0.5$. 
The exponential in the initial data is present on the real line (for the pure exponential the series truncates and the Stokes line is in fact inactive; it is not in any case important  for the current discussion).}
\label{linesF2lambda05}
\end{figure}

\begin{figure}
\centerline{
\includegraphics[width=0.45\textwidth,height=.35\textwidth]{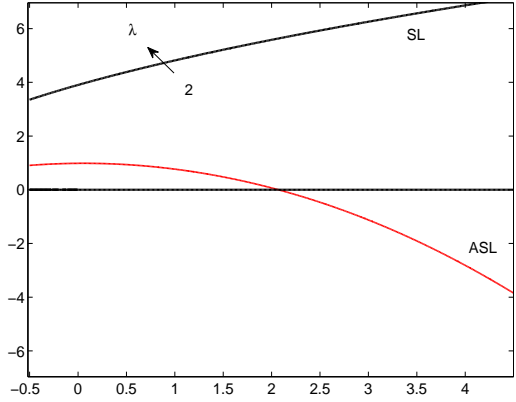}
}
\caption{$\lambda>1$: The $\mbox{Im}(F_2)=\l \mbox{Im}(\xi)$ Stokes line and $\mbox{Re}(F_2)=\l\mbox{Re}(\xi+\l/(1-\l^2))$ anti-Stokes lines for $\l=1.5$. The contribution $2$ switches off the exponential from the initial data across the Stokes lines shown, as does the contribution $3$ across the Stokes line that is the mirror image in the real $\xi$ axis of the one shown.}
\label{linesF2lambda15}
\end{figure}

The corresponding lines 
for $F_3$ are the mirror image in the real axis. Figures \ref{linesF2lambda05}($\l<1$) and \ref{linesF2lambda15} ($\l>1$) show the 
Stokes and anti-Stokes lines of $F_2$ compared to the exponential contribution of the initial data (again the lines for $F_3$ are the 
mirror image in the real axis). For $\l>1$ the branch $F_2$ that gives a slower exponential decay is turned on across the Stokes line. 
Similarly, $F_3$ is turned on across the Stokes line lying in the fourth quadrant. For $\l<1$, the contribution from the initial data 
would turn on $F_2$ across the Stokes line, but it series does not diverge, the Stokes line is inactive. We note that 
it would give a negative wave speed, so in any case it plays no role on front selection.

\section{Analysis of the asymptotic regions}\label{sec:3}
\subsection{Analysis of the leading edge of the front}\label{section:2.3new}
While the analysis that follows is in some respects standard, one of its implications is not. Setting 
\begin{equation}\label{2.20}
x = z + ct\,,\ v = e^{\mu t} e^{-\lambda z} V\,,
\end{equation}
where $c$, $\mu$, $\lambda$ are for the moment arbitrary, in the linearised equation (\ref{lin:equ}) yields 
\begin{eqnarray}
(1-\lambda^2) \frac{\pa V}{\pa t} &=& (-\mu -c\lambda -\Phi_u\l^2 +\mu\l^2 + c\l^3)V -(-c-2\Phi_u\l+2\mu\l+3c\l^2)\frac{\pa V}{\pa z} \nonumber \\
&+&(-\Phi_u +\mu +3c\l)\frac{\pa^2 V}{\pa z^2}-c\frac{\pa^3 V}{\pa z^3}
-2\l \frac{\pa^2 V}{\pa z\pa t} +  \frac{\pa^3 V}{\pa z^2\pa t} \,.
\label{2.21}
\end{eqnarray}
By choosing $\mu$ and $\l$ to satisfy 
\begin{eqnarray}
\mu + c \l & = & -\Phi_u \l^2 +\mu \l^2 +c\l^3 \,,\label{2.22}
\\
c  & = & - 2\Phi_u \l + 2\mu \l +3c\l^2\,,
\label{2.23}
\end{eqnarray} 
((\ref{2.23}) corresponding to a repeated root condition of (\ref{2.22}))
so that $c(1-\l^2)^2 +2\Phi_u\l=0$. We eliminate the $V$ and $\pa V/\pa z$ 
terms from (\ref{2.21}) to get 
\begin{equation}\label{2.24}
(1-\l^2) \frac{\pa V}{\pa t} \sim (-\Phi_u+\mu+3c\l)\frac{\pa^2 V}{\pa z^2}
\end{equation}
as a putative large-time balance for $z=O(t^{1/2})$. Note that (\ref{2.22}) and (\ref{2.23}) are equivalent to (\ref{clairaut:singular}) on identifying $\mu$ with $-F$, $\l$ with $p$ and $c$ with $\xi$. 

Now if in (\ref{2.22}) we require $\mu$ to be imaginary, as well as $c$ 
necessarily being real, we have two complex equations for $i\mu$, $c\in\R$ and 
$\l\in\C$ with solutions 
\begin{eqnarray}
(c,\mu,\l)&=&( 0\,,\ 0,\, 0) \,,\label{2.25}  \\
(c,\mu,\l) &\approx& (+0.7872\Phi_u\,,\  -1.1688 \Phi_u i\,, \ +1.0419+0.8337 i)\label{2.26} \,, \\
(c,\mu,\l)&\approx& (+0.7872\Phi_u\,,\   +1.1688\Phi_ui\,,\ +1.0419-0.8337 i)
\label{2.27}\,, \\
(c,\mu,\l) &\approx&( -0.7872\Phi_u\,,\  +1.1688\Phi_ui\,,\ -1.0419+0.8337 i)
\label{2.28}\,, \\
(c,\mu,\l) &\approx &(-0.7872\Phi_u\,,\  -1.1688\Phi_ui\,,\ -1.0419-0.8337 i)
\label{2.29}\,,
\end{eqnarray}
corresponding to the results on Section\ref{section:2.2} on identifying $\lambda$ with $\mbox{Re}(p^*)$, $c$ with $\xi^*$ and $\mu$ with $F(p^*)$ (see (\ref{speed}), (\ref{w:star}) and (\ref{univ:decay}). 
In the two relevant cases (\ref{2.26}) and (\ref{2.27}) (with $\xi^*>0$), equation (\ref{2.24}) reads 
\begin{equation}\label{2.30}
\frac{\pa V}{\pa t} = D \frac{\pa^2 V}{\pa z^2} \quad \mbox{with} \ D = \frac{-\Phi_u +\mu +3c\l}{1-\l^2}\approx(-0.1474\pm 0.8923 i)\Phi_u 
\end{equation} 
an unusual aspect of which is that the real part of the diffusivity is slightly negative (observe that $\Phi_u=1$ is the maximum possible value of $\Phi_u$; $\max_{\{u: \ \phi'(u)<0\}}\{-\phi'(u)\}=1=-\phi'(0)$), so that (\ref{2.30}) is of backward heat equation type. The appropriate boundary condition on (\ref{2.24}) describing the large-time behaviour, $t\to+\infty$ with $z=O(t^{1/2})$, is (cf. \cite{CK})
\begin{equation}\label{2.31}
\mbox{on} \quad z=0\,, \quad V=0\,, 
\end{equation}
corresponding to matching into a modulated travelling wave (pertaining for $z=O(1)$) with repeated-root far-field behaviour of the form 
\[
ze^{ i w(q^*) t}\e^{- \lambda^* z} \quad\mbox{as}\ z\to +\infty\,.
\]
Intriguingly, the Stokes lines analysis outlined in Section~\ref{sec:expdecay} implies that the solution to (\ref{2.30}) can have no steady far-field behaviour so its large-time behaviour is generically of the form 
\begin{equation}\label{2.32} 
V\sim \frac{Iz}{t^\frac{2}{3}} e^{\frac{z^2}{4Dt}} \quad \mbox{as} \ z\to +\infty
\end{equation}
for some constant $I$: this can be viewed as a self-similar solution of the second kind 
, the upshot being that in the current context it does not matter that the equation is of backward type, it seems likely that there are numerical implications, however (see below). Equivalently, it is important to keep in mind the properties of the solution in $\C$, not simply in $\R$: (\ref{2.30}) can then be viewed as a forward equation in suitable directions in $\C$. Matching into (\ref{2.32}) implies that the modulated travelling wave is of the form
\[
u\sim U\left(x-\xi^* t + \frac{3}{2 \l^*}\ln t; t+\frac{3i\ln \l^*}{2F(p^*)}\ln t\right)\,,
\]
being periodic in its second argument with period $T=2\pi i /F(p^*)$ ($\in \R$).

\begin{figure}
\centering
\mbox{
\subfigure[]
{\includegraphics[width=0.45\textwidth,height=.35\textwidth]{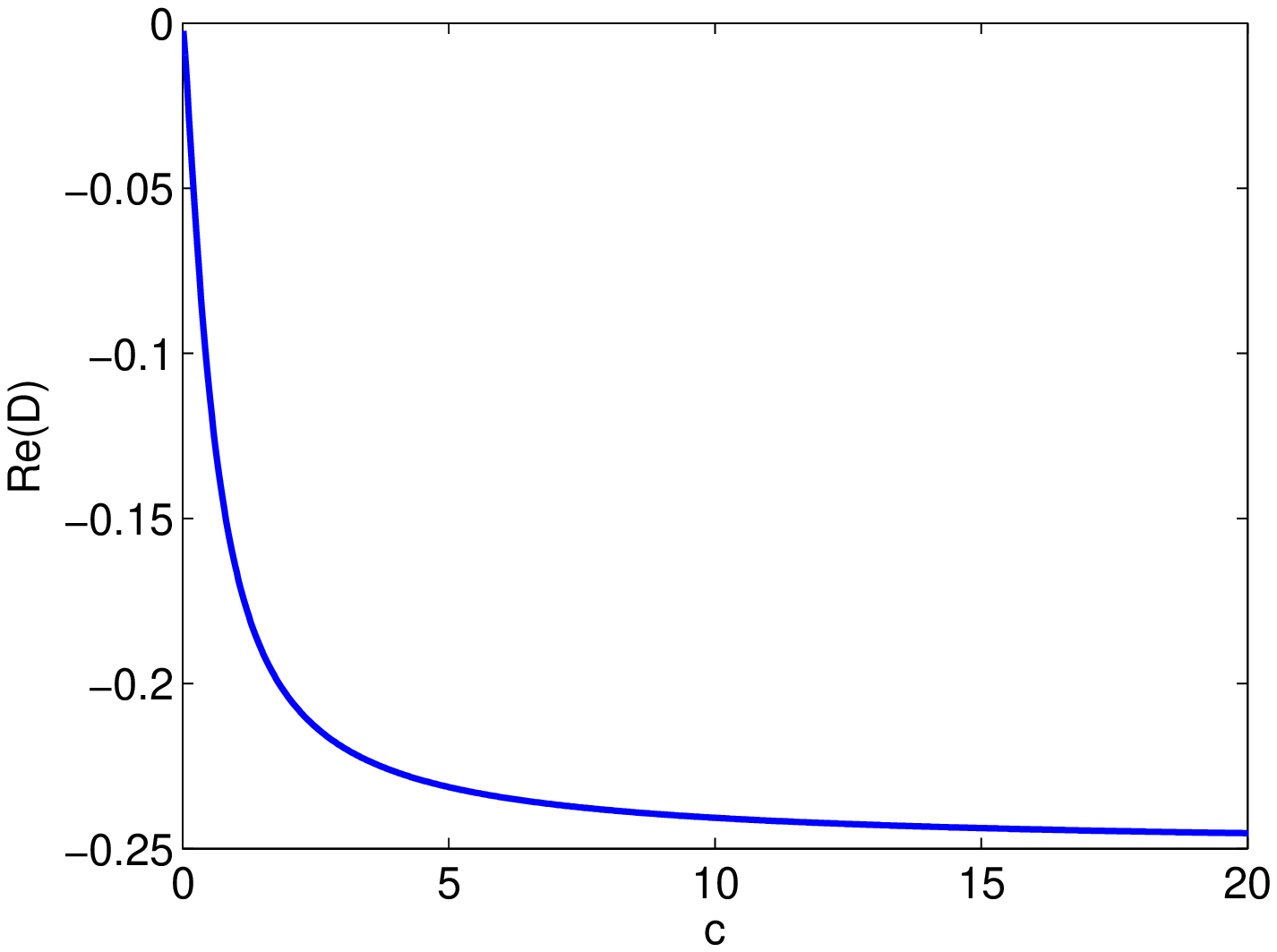}
\label{D:real}}
} 
\quad
\mbox{
\subfigure[]
{\includegraphics[width=0.45\textwidth,height=.35\textwidth]{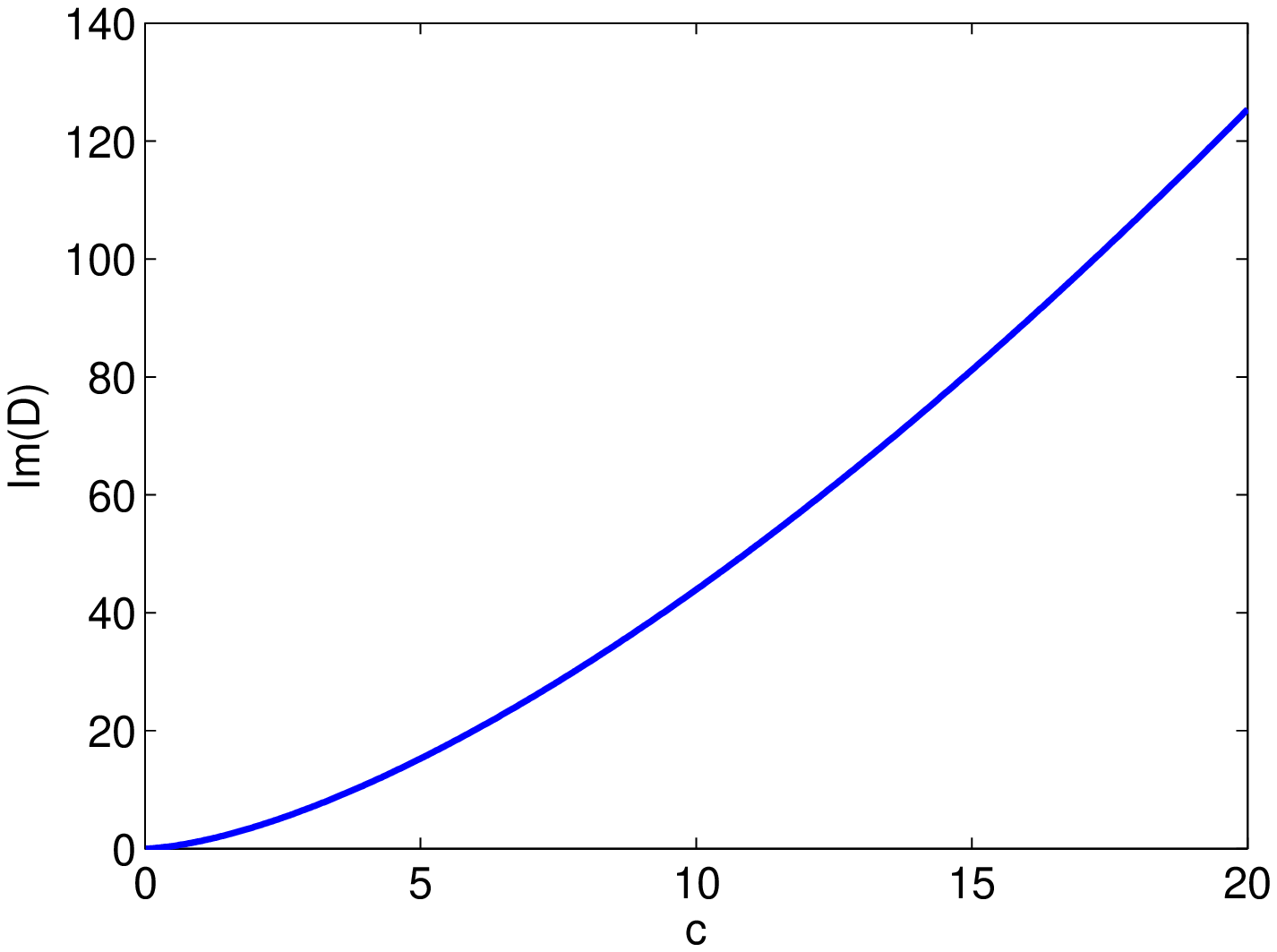}
\label{D:imag}}
}
\caption{Real and imaginary parts of $D$ for the branch solutions of (\ref{2.22})-(\ref{2.23}) that select the 
critical wave speed $\xi^*$. Observe that the complex conjugate of $D$ is also 
relevant and that $\mbox{Re}(D)<0$ and $\mbox{Re}(D)/\mbox{Im}(D)\to 0$ as $c\to +\infty$.
}
\label{D:complex}
\end{figure}

A second application of (\ref{2.22})-(\ref{2.23}) is to identify the location of the 
transition between backward and forward diffusion if it exists. If now we fix $c\in\R$ we can determine 
$\mu$ and $\l\in \C$ and identify for what $c$ the diffusivity $(-\Phi_u+\mu +3c\l)/(1-\l^2)$ is purely imaginary. It turns out, however, that for the relevant branch solutions of (\ref{2.22})-(\ref{2.23}) $\mbox{Re}(D)$ remains negative for all $c$ and asymptotes to a small negative value, whereas $\mbox{Im}(D)$ becomes unbounded as $c\to +\infty$. Figure~\ref{D:complex} shows a picture of $\mbox{Re}(D)$ and $\mbox{Im}(D)$ for the relevant solutions of (\ref{2.22})-(\ref{2.23}) as a function of $c$.

\subsection{Wavelength selection behind the front}\label{section:2:3}
In this section we give some ingredients for analysing the pattern behind the front. 
We assume that the approximation (\ref{predict:decay}) in the linear regime is valid for all $x$, at least initially. 
In the discussion that follows we appeal to symmetry in concentrating on the right-hand side of the growing perturbation, 
i.e. on the front that propagates to the right with speed $\xi^*$.

Observe that in the moving frame $z=x-s(t)$, where $s(t)=\xi^*t+o(t)$ as $t\to +\infty$, 
the approximation (\ref{predict:decay}) is periodic in $t$ with period
\bequ\label{T:period}
T:= \frac{2 \pi i}{F(p^*)}\approx 5.375\,.
\eequ
This suggests that in the transition regime the solution can be described by the modulated travelling wave
\bequ
u(x,t) \sim U(z,t)\quad \mbox{as} \  t\to +\infty \quad \mbox{with} \ z=O(1)\, ,
\label{F:long:time}
\eequ
whereby
\bequ
U(z,t+T)=U(z,t)\quad \mbox{and} \quad  U(+\infty,t)=u_u\, . \label{F:conditions} 
\eequ
Further, (\ref{F:long:time})-(\ref{F:conditions}) are consistent with the solution approaching a periodic-in-$x$ pattern with wavelength
\[
X:=  \xi^* T 
\]
behind the front, whereby
\bequ\label{F:left:boundary}
U(z,t)\sim u_s(z+\xi^* t) \quad \mbox{as} \ z \to -\infty\,;  
\eequ
where $u_s(\zeta)$ is a $X$-periodic function of $\zeta$ and will be piecewise constant for (\ref{phi:cubic}). 
Assuming that 
\bequ\label{log:correct}
s(t) = \xi^* t -\nu \ln t + x_0 + o(t) \quad \mbox{as} \ t\to +\infty 
\eequ
for constants $x_0$ and $\nu$ (with $\nu=3/(2\lambda^*)$, cf. \cite{frontsreview} and see Section~\ref{section:2.3new} above), this implies that the steady-state solution that is left behind takes the form
\bequ\label{logcorrect}
u \sim u_s(x+\nu \ln(x/\xi^*)-x_0) \quad \mbox{for} \quad x=O(t)
\eequ
which in turn has spatial period approaching $X$ as $x \to +\infty$. Observe that $X$ does not depend on $u_u$ and $\phi$, 
and by (\ref{speed}) and (\ref{T:period}) we have
\bequ\label{X:period2}
X= \xi_0 \frac{2\pi i}{ p^* \xi_0 +\frac{(p^*)^2}{1-(p^*)^2} }\approx 4.232 \,.
\eequ

Let us now address how the value $\phi(u_s)$ in (\ref{F:left:boundary}), 
which is necessarily a constant, is related to $\phi(u_u)$. 
Suppose we are in the transition regime where (\ref{F:long:time}) 
is valid, then to leading order $U$ satisfies
\bequ\label{mtw}
\frac{\pa U}{\pa t}-\xi^* \frac{\pa U}{\pa z} =\frac{\pa^2}{\pa z^2}\left(\phi(U) + \frac{\pa U}{\pa t} -\xi^* \frac{\pa U}{\pa z} \right) \,,
\eequ
and integration with respect to $t$ over a temporal period gives
\[
-\, \xi^* \frac{d}{dz}\langle U\rangle =\frac{d^2}{dz^2}\left( \langle\phi(U)\rangle -\xi^* \frac{d}{dz} \langle U\rangle \right) \,;
\]
here $\langle\cdot\rangle=1/T \,\int_0^T \cdot \,dt$ (the `cell' average). Now integrating with respect to $z$ subject to 
$\langle U \rangle \to u_u$ as $z\to +\infty$ we obtain 
\[
-\xi^* (\langle U\rangle -u_u) =\frac{d}{dz}\left( \langle \phi(U)\rangle -\xi^* \frac{d}{dz} \langle U\rangle \right)\,,
\]
so that
\[
\langle  U \rangle\to u_u \quad \mbox{as} \quad z\to -\infty \,,
\]
which in turn implies that (by (\ref{F:left:boundary})) 
\bequ\label{cons:symmetric}
\frac{1}{X}\int_{x_0}^{x_0+X} u_s\,dz= u_u 
\eequ
(corresponding to conservation of mass) and that
\bequ\label{us:relates:ends}
\xi^* \int_{-\infty}^{\infty} (\langle U\rangle - u_u) \,dz = A - \phi(u_u)\quad \mbox{as} \ t \to +\infty\,,
\eequ
with 
\[
A:=\phi(u_s(x))\,. 
\]
Let us now look at the $X$-periodic pattern. First we observe that, by (\ref{mass}), periodicity (in $x$) implies mass 
conservation in a periodicity cell, as expected by (\ref{cons:symmetric}). For (\ref{phi:cubic}), we expect an oscillatory 
pattern alternating between two constant values $u_-$ and $u_+$ with $u_-<u_+$ such that 
(\ref{phi:const}) holds and $\phi'(u_-)$, $\phi'(u_+)>0$. This, however, cannot happen if 
$\phi$ is of the form (\ref{phi:model2}) or is given by (\ref{phi:exposym}), because two values in the stable region 
satisfying (\ref{phi:const}) do not exist, and we expect $u_-$ and $u_+$ to depend on $t$. 
In the case $\phi(u)=-ue^{-u^2}$ we anticipate that $u_-(t) \to -\infty$ and $u_+(t) \to +\infty$. 

The relation (\ref{us:relates:ends}) simplifies in the symmetric cases with $u_u=0$: equation (\ref{mtw}) is invariant under the transformation $U\to-U$, 
since these $\phi$'s are odd functions, 
implying (under a suitable uniqueness assumption) that $-U$ is a translation of $U$ (by half the period $T$) in $t$, so that $\langle U\rangle =-\langle U\rangle$ and hence $\langle U\rangle =0$. Together with (\ref{us:relates:ends}), this implies that
\bequ\label{symmetric:A}
A=\phi(0)=0\, . 
\eequ
More generally, it is unclear that $A$ can be determined explicitly, though it is 
presumably a property of the modulated travelling wave; this is the one respect 
in which we leave the non-symmetric case open.

\subsection{The approach to a steady pattern}\label{section:3.2}
We now analyse how the left and right values of a `jump' from $u_-(t)$ to $u_+(t)$ are approached as $t \to +\infty$ (assuming there 
is a sharp transition as suggested by the construction of steady states for $\phi$ 
of the form (\ref{phi:model1})). 
In this narrow region the dominant balance is given as $t\to+\infty$ by
\bequ\label{jump:balance}
0 = \frac{\pa^2}{\pa x^2}\left(\phi(u) +\frac{\pa u}{\pa t}\right)\, , 
\quad \mbox{with} \quad u(-\infty,t)=u_-(t), \ u(+\infty,t)=u_+(t)\,,
\eequ
so integrating twice with respect to $x$ gives
\bequ\label{t:ode}
\phi(u) + \frac{\pa u}{\pa t} =\phi(u_+(t)) +\frac{d}{dt}u_+(t)\,,
\eequ
and
\bequ\label{t:ode:int}
\phi(u_+(t))+\frac{d}{dt}u_+(t)=\phi(u_-(t))+\frac{d}{dt} u_-(t)\,.
\eequ

Observe that (\ref{t:ode}) is an ODE in $t$ and as it stands contains no information about the 
$x$-dependence, though this does enter through the initial data $u(x,0)$, which we shall take smooth (given that 
the solution inherits the regularity of the initial data, a smoothness condition can be imposed for large time if it holds initially).

We now outline in more detail how a piecewise-constant steady state $u_s(x)$ is attained for a 
nonlinearity $\phi$ of the from (\ref{phi:model1}). We take the limit profile as $t\to +\infty$ to be 
\beqs
\lefteqn{u_s(x)=u_+ \quad \mbox{for}\ x \in (x_{2n},x_{2n+1})\,,}\\
\lefteqn{u_s(x)=u_u \quad \mbox{for} \ x=x_m\,,}\\
\lefteqn{u_s(x)=u_- \quad \mbox{for}\ x \in (x_{2n+1},x_{2n+2})\,,}
\eeqs
for integers $n$ and $m$, with $\phi(u_+)=\phi(u_u)=\phi(u_-)=0$ and with $x_{2n+1}-x_{2n}\to X$ as $n\to+\infty$.
We concentrate on the range $x \in [x_{2n},x_{2n+1})$, with the remainder 
following by obvious symmetry arguments. We also introduce
\[
\Phi_s= \phi'(u_+)=\phi'(u_-)\, ,
\]
wherein the final equality holds because we are restricting attention to symmetric cases.
This simplifies a number of considerations that follow; in particular, the outer solution has
\[
u=u_+ + W \quad x\in (x_{2n},x_{2n+1}), \ u=u_-+ W \quad x\in(x_{2n+1},x_{2n+2})\,,
\]
giving to leading order
\bequ\label{v:eq} 
\frac{\pa W}{\pa t}=\Phi_s\, \frac{\pa^2 W}{\pa x^2} + \frac{\pa^3 W}{\pa x^2 \pa t}
\eequ
for all $x$, with the associated continuity conditions
\[
\left[\Phi_s\, W + \frac{\pa W}{\pa t}\right]_-^+=\left[\Phi_s \,\frac{\pa W}{\pa x} +\frac{\pa^2 W}{\pa x\pa t}\right]_-^+=0 
\quad \mbox{at} \ x=x_m
\]
following on matching into the inner regions below (or on intuitive grounds). In non-symmetric cases, 
the diffusivity in the second term in (\ref{v:eq}) depends on whether $u\sim u_+$ or $u\sim u_-$, further complicating the analysis.

\paragraph{The inner region $x-x_{2n}=O(e^{-\omega t})$}

Introducing the (exponentially-narrow) large-time inner scaling
\[
x=x_{2n}+ Z e^{-\omega t}\,,
\]
for some constant $\omega$, gives
\bequ\label{hat:kappa:eq}
\frac{\pa u}{\pa t} + \omega Z \frac{\pa u}{\pa Z} = e^{2\omega t}\frac{\pa^2}{\pa Z^2}\left( \phi(u) 
+ \frac{\pa u}{\pa t} + \omega Z \frac{\pa u}{\pa Z}\right)\,.
\eequ
Setting $u \sim u_0(Z)$ as $t \to +\infty$ implies (in order to match outwards) that 
\bequ\label{hat:kappa:match}
\omega Z \frac{d u_0}{dZ}= A -\phi(u_0)\,. 
\eequ
Linear behaviour, $u_0-u_u\sim K(t)\, (x-x_m)$ as $x \to x_m$ then requires, generically, 
that $\omega=\Phi_u$; the first non-generic case, in which $K(t)=o(e^{\Phi_u t})$ as $t \to +\infty$, will 
instead have $u_0-u_u = O(Z^3)$ as $Z \to 0$, so that $\omega=\Phi_u/3$. Henceforth we set $\omega =\Phi_u$. 
Equation (\ref{hat:kappa:match}) then fixes $u_0$ up to a rescaling of $Z$, 
corresponding to a translation of $t$, in the form
\bequ\label{kappa:int}
\int_{u_u}^{u_0} \left( \frac{\Phi_u}{A-\phi(u')}-\frac{1}{u'-u_u}\right)\, du' + \ln(u_0-u_u) =\ln(\alpha_n Z) \,
\eequ
wherein the rescaling $\alpha_n$ depends on the initial data. Moreover,
\bequ\label{moreover}
u\sim u_+ - \beta_n Z^{-\kappa} \quad \mbox{as} \ Z \to +\infty \quad \mbox{with} 
\quad \kappa:=\frac{\Phi_s}{\Phi_u}
\eequ
for some positive constant $\beta_n  \propto \alpha_n^{-\kappa}$.

We remark that such inner regions can be viewed as being initiated by the modulated travelling wave at successive 
time intervals of $T$, so that in (\ref{kappa:int}) we have 
\[
\alpha_n \sim \hat{\alpha} \, e^{-n \Phi_u  T}, \ \beta_n \sim \hat{\beta} e^{n\Phi_s  T} \quad  \mbox{as} \ n \to +\infty
\]
for some constants $\hat{\alpha}$, $\hat{\beta}$. For similar reasons we also have 
\[
x_{2n} \sim n X + \nu \ln (n T) + \tilde{x} \quad \mbox{as} \ n \to +\infty
\] 
for some $\tilde{x}$, so the quantity $\alpha_n Z$ in (\ref{kappa:int}) becomes 
\bequ\label{zeta:n}
\alpha_n Z \sim  \hat{\alpha} (x- n\, \xi^* T -\nu\ln (nT) - \tilde{x}) e^{\Phi_u t} e^{- n \Phi_u  T} \quad \mbox{for large} \ n \,. 
\eequ
We can confirm the associated periodicity constraint by taking $t \to t+T$, $x \to x +\xi^* T$ in (\ref{zeta:n}) to yield $\alpha_{n-1} Z$ to leading order, as required for matching with the tail of the modulated travelling wave.

\paragraph{The outer region, $x-x_{2n}=O(1)$.} There are two distinct contributions to the solution. 
Equation (\ref{v:eq}) has separable solutions 
\bequ\label{separable}
e^{-   \Phi_s \frac{ k^2}{1+k^2} t} \{ \cos (kx), \sin (kx)\}\, ,
\eequ
with $k^2 \in\R$, whose amplitude will depend on the initial data. 
It is noteworthy that the decay rate in (\ref{separable}) saturates as $k^2\to +\infty$ 
(to $e^{-\Phi_s t}$, which will be significant later), as does the corresponding growth rate in the 
backward-diffusion range of $\phi$.

The separable solutions relevant to matching forward into 
the modulated travelling wave take the form
\[
e^\frac{2\pi n i t}{T} \, e^{\sigma (x-\xi^* t)}, \ \mbox{Re} \,\sigma > 0
\]
for an integer $n$, yielding the dispersion relation
\[
(\sigma^2 -1) \left(\xi^* \sigma -\frac{ 2\pi n i }{T}\right) = \Phi_s \sigma^2
\]
for $\sigma$ with $\mbox{Re} \sigma >0$ (for $n\neq 0$ there are two roots of this polynomial with positive real part, as can be seen 
by computing the Cauchy index, which is $-1$, see e.g. \cite{coppel}; for $n=0$, one root is $\sigma=0$, the other two being real and having opposite sign). The amplitudes of these modes will in effect be fixed by 
the tail of the modulated travelling wave.

\paragraph{The intermediate region, $x-x_{2n}=O(t^{-1/2})$.} In order to match with (\ref{moreover}) 
in this intermediate layer, we set
\bequ\label{w:def}
W=e^{-\Phi_s t} G
\eequ
in (\ref{v:eq}) to give 
\bequ\label{w:eq}
\frac{\pa G}{\pa t} - \Phi_s G= \frac{\pa^3 G}{\pa x^2 \pa t}\,. 
\eequ
Seeking an asymptotically-self-similar solution
\bequ\label{asymp:self:similar}
G \sim t^\frac{\kappa}{2} {\cal G}(\eta), \ \eta =(x-x_{2n})t^\frac{1}{2}
\eequ
as $t \to +\infty$, the first term in (\ref{w:eq}) is negligible and
\bequ\label{W:eq}
\eta \frac{d^3 {\cal G}}{d \eta^3} + ( \kappa + 2) \frac{d^2 {\cal G}}{d \eta^2} = -2 \Phi_s {\cal G}\,. 
\eequ
We require, to match with (\ref{moreover}), together with higher-order matching into the inner region, that
\bequ\label{W:match}
{\cal G} (\eta)= -\beta_n \eta^{-\kappa} + O(\eta) \quad \mbox{as} \ \eta \to 0
\eequ
(i.e. that no $\eta^0$ term be present). As $\eta \to +\infty$, 
the JWKB method yields the three possible asymptotic forms to be
\bequ\label{W:ln:form}
\ln {\cal G}(\eta) \sim -\frac{3}{2} (2 \Phi_s )^\frac{1}{3}\eta^\frac{2}{3} \{ 1, e^{\pm \frac{2\pi i}{3}}\}\, .
\eequ
Two of these are exponentially growing but, perhaps unexpectedly, 
these cannot be suppressed; a boundary condition count implies 
that a single constraint (on the ratio of the two exponentially growing terms) is required, 
and a full analysis of this would necessitate discussion of the associated Stokes phenomenon, 
whereby these growing terms are ultimately switched off, (\ref{w:def}), 
(\ref{asymp:self:similar}) 
being exponentially subdominant to (\ref{separable}) as $t\to +\infty$ for all finite $k$. We shall not pursue such 
an analysis, but emphasise that this mechanism of selecting boundary conditions for the intermediate-asymptotic 
similarity solution (\ref{asymp:self:similar}) is an unusual one. While the similarity solution 
decays faster than any of the modes (\ref{separable}), 
it is of interest in view of its singular matching condition (\ref{W:match}) and should be visible over appropriate 
scalings lying between inner and intermediate. The oscillatory exponential growth associated 
with (\ref{W:ln:form}), however, occurs in a regime over which (\ref{asymp:self:similar}) is negligible.

\section{Comparison with numerical results}\label{section:3}
In this section we complete the asymptotic analysis of the patterns. 
We further compare the analytical predictions of this and the previous sections with numerical results. In the numerical examples we choose fast decaying initial conditions. Further examples with exponentially decaying initial conditions are considered in Section~\ref{section:new}.

For both (\ref{phi:cubic}) and (\ref{phi:exposym}) 
we solve equation (\ref{main:eq}) 
numerically on an interval $(-L,L)$ truncated for numerical purposes, with $L=50$ or $L=100$, subject to the initial data
\bequ\label{init:symm:nume}
u_0(x) = 0.1 \,e^{-x^2} \,,
\eequ
(a perturbation of $u_u=0$), and the symmetric boundary conditions
\bequ\label{bc}
\frac{\pa}{\pa x}\left(\phi(u) + \frac{\pa u}{\pa t}\right) =0 \quad \mbox{at} \ x=\pm L\,. 
\eequ
We use an explicit-in-time method. For the examples computed in this section we introduce the unknown
\[
g := \phi(u) + \frac{\pa u}{\pa t}\,,
\]
and solve the elliptic problem
\beqs
\lefteqn{ - \frac{\pa^2 g}{\pa x^2} + g = \phi(u),  \quad x \in (-L,L) \, ,} \label{ellip} \\
\lefteqn{ \frac{\pa g}{\pa x}(\pm L) = 0 \, ,}\label{ellip:bc}
\eeqs
by Gauss elimination at every time step, where $u$ is the value of the solution at the previous time step. 
The value at the next step is then obtained by solving 
\bequ
\frac{\pa u}{\pa t} = \frac{\pa^2 g}{\pa x^2} \quad x \in (-L,L), \ t>0\,, \label{evol}
\eequ
by forward Euler time discretisation. Explicit in time methods for one-dimensional pseudo-parabolic 
equations have been proved to be stable for a small enough temporal step (even allowing backward diffusion) and regardless of the 
spatial step size, 
assuming a priori that the solution remains bounded, see e.g. \cite{thomee}, where a single-step method is discussed, and also 
the earlier works \cite{Ewing1} and \cite{Ewing2}. The convergence results for the semidiscrete (discrete in space) problem can be found in \cite{pierre} 
(for the nonlinearity (\ref{cubic})).

The cases we consider share some features. 
Firstly, for both $\phi$'s (i.e. (\ref{phi:cubic}) and (\ref{phi:exposym})) we have
\bequ\label{phi:values@0}
\phi(0)=0,\quad \phi'(0)=1\,. 
\eequ
Also the numerical results confirm that the wave speed is here of the pulled-front type, i.e. is given by $\xi^*=\xi_0$, see (\ref{speed}).
Moreover, (\ref{symmetric:A}) and (\ref{X:period2}) are confirmed numerically in these cases (see below).

\begin{figure}[htb]
\centering
\mbox{
\subfigure[The initial data $u_0$ and the solution at $t=30$.]
{
\includegraphics[width=0.46\textwidth,height=.35\textwidth]{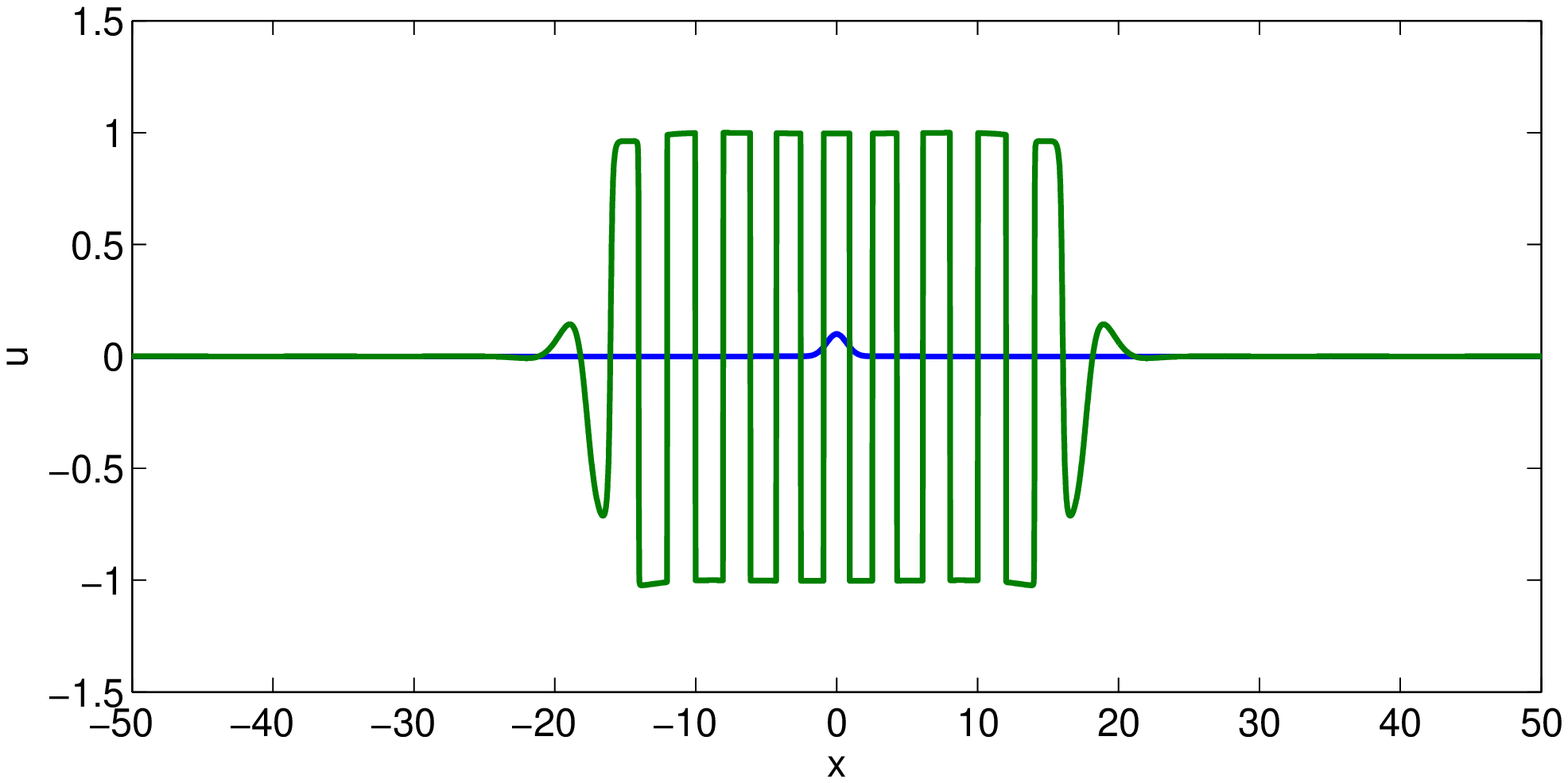}
\label{uns0:1}}
}
\mbox{
\subfigure[The solution at $t=60$.] 
{
\includegraphics[width=0.46\textwidth,height=.35\textwidth]{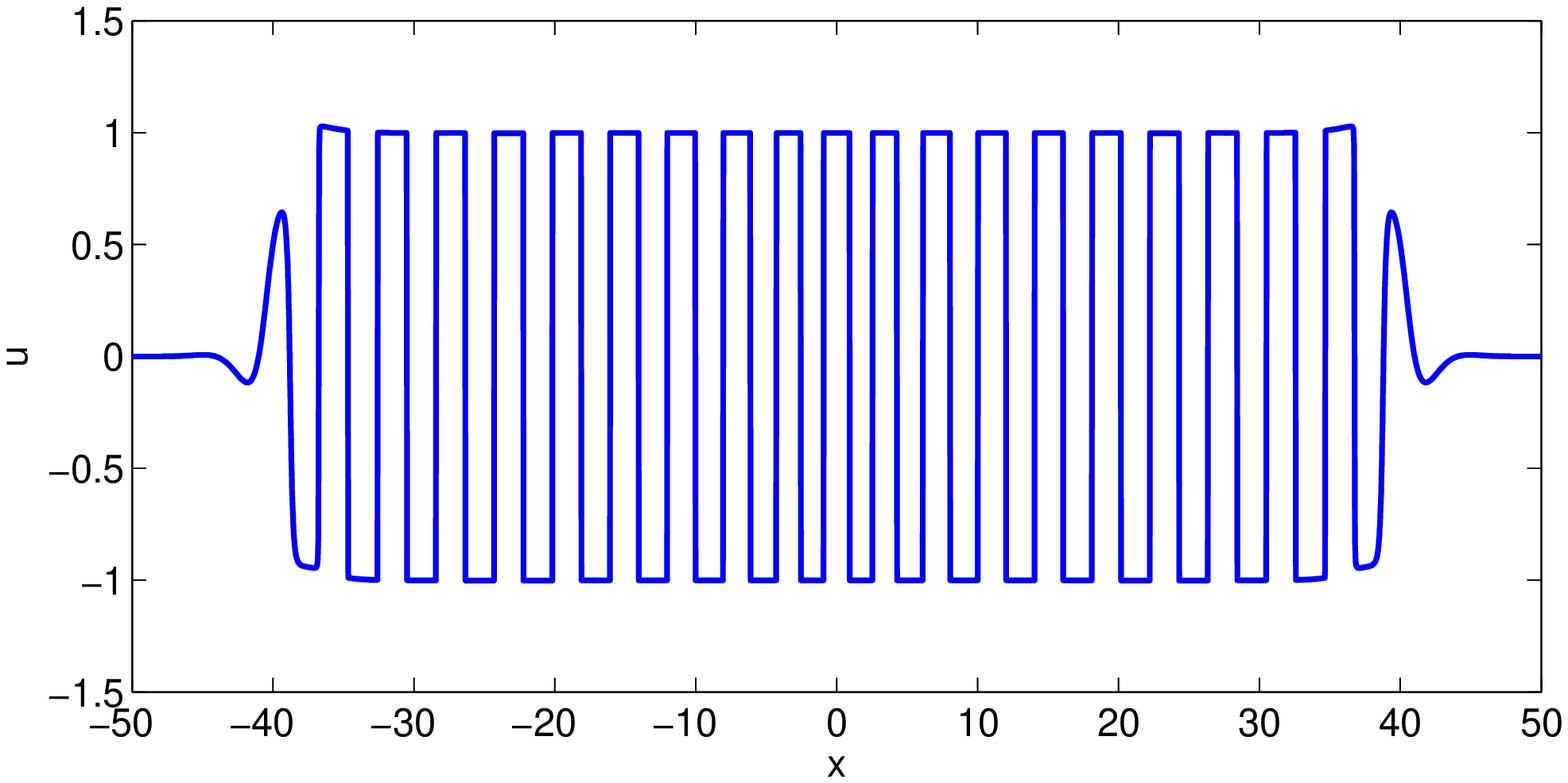}
\label{uns0:2}}
}
\caption{Numerical solutions of (\ref{main:eq}) with $\phi(u)=u^3-u$ and $u_0=0.1 e^{-x^2}$ at times $t=26$ and $t=54$, illustrating 
how the domain is 
invaded by a periodic piecewise constant solution, with constant values alternating from $-1$ to $1$. The time step is $\Delta t=0.01$ 
and the spatial one is $\Delta x=0.1$.}\label{unstable:0}
\end{figure}

\begin{figure}[htb]
\centering
\mbox{
\subfigure[The initial data $u_0$ and the solution at $t=30$.]
{
\includegraphics[width=0.46\textwidth,height=.35\textwidth]{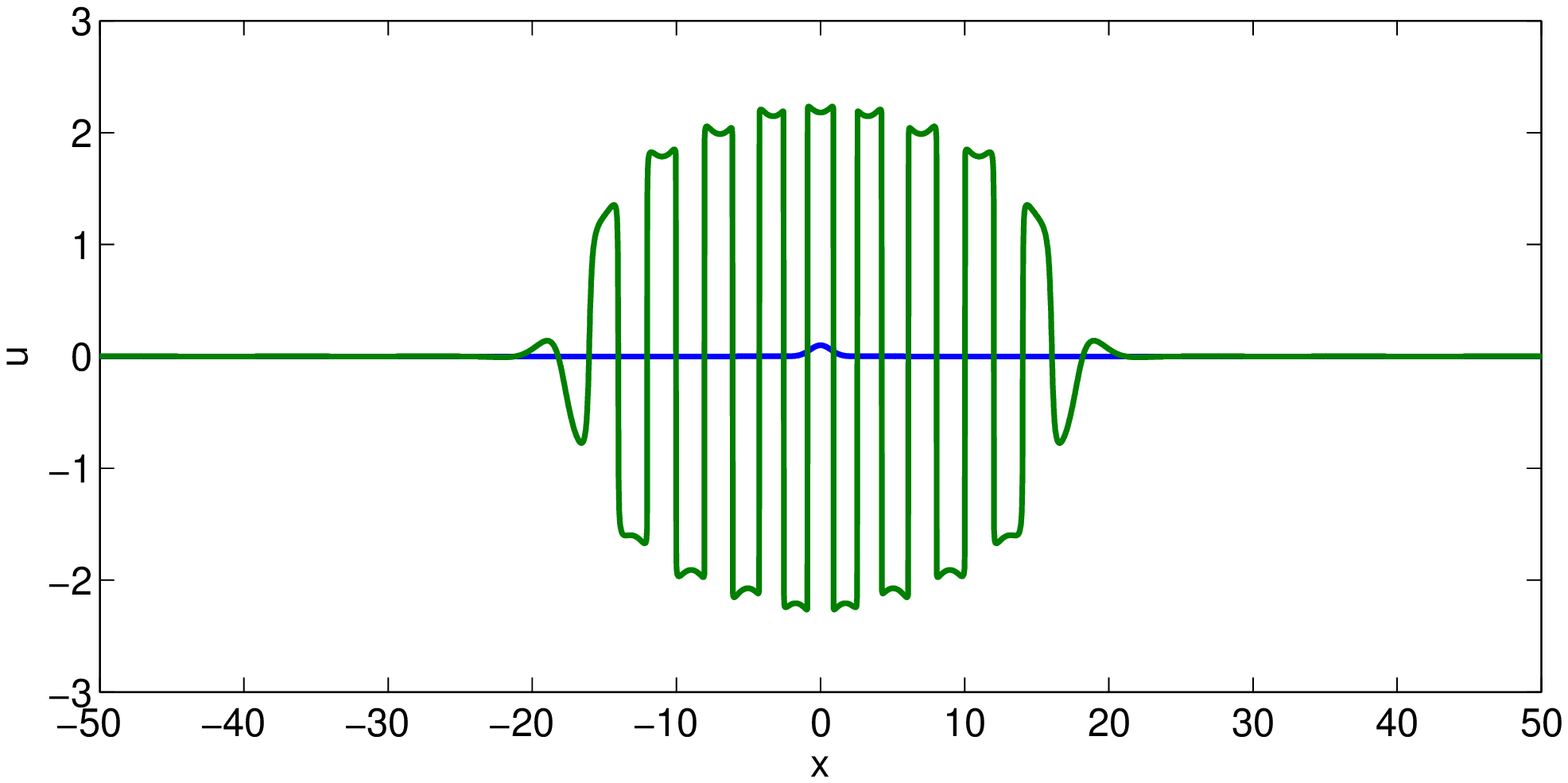}
\label{uns0:3}}
}
\mbox{
\subfigure[The solution at $t=60$.]
{
\includegraphics[width=0.46\textwidth,height=.35\textwidth]{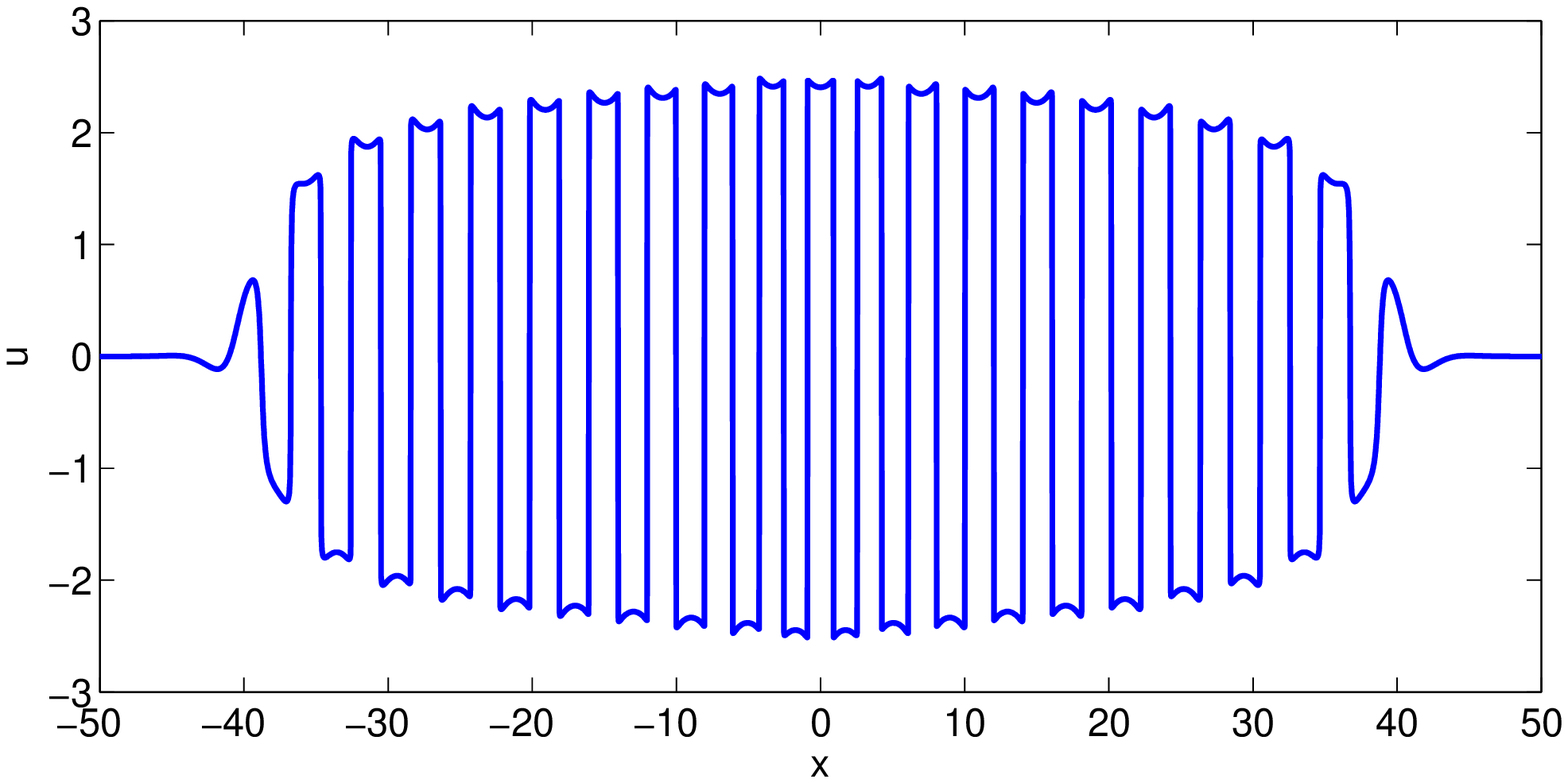}
\label{uns0:4}}
}
\caption{Numerical solutions of (\ref{main:eq}) with $\phi(u)=-u\,e^{-u^2}$ and $u_0=0.1 e^{-x^2}$ at times $t=30$ and $t=60$, 
illustrating how the domain is invaded by a periodic solutions, where the absolute value of the maximum and minimum values of 
the solution in a periodic cell increase with $t$. The time step is $\Delta t=0.01$ and the spatial one $\Delta x=0.1$.
}
\label{2unstable:0}
\end{figure}

As an illustration we start by showing numerical results for (\ref{phi:cubic}) 
in Figure~\ref{unstable:0}, plotting the profile $u(x,t)$ against $x$ for various values of $t$. 
There are two fronts, one moving to the left and one to the right. 
The simulation shows how the domain is ultimately filled (away from the origin) by a near-periodic solution. 
Solutions get close to a piecewise constant solution with spatial oscillations 
between the values $u=-1$ and $u=1$. A similar simulation is shown in Figure~\ref{2unstable:0} for (\ref{phi:exposym}). 
In this case $u_+(t) \to +\infty$ and $u_-(t)\to -\infty$ as $t\to +\infty$.

\begin{figure}
\centering
\mbox{
\subfigure[Numerical solution for $\phi(u)=u^3-u$ at $t=58$.]
{\includegraphics[width=0.45\textwidth,height=.35\textwidth]{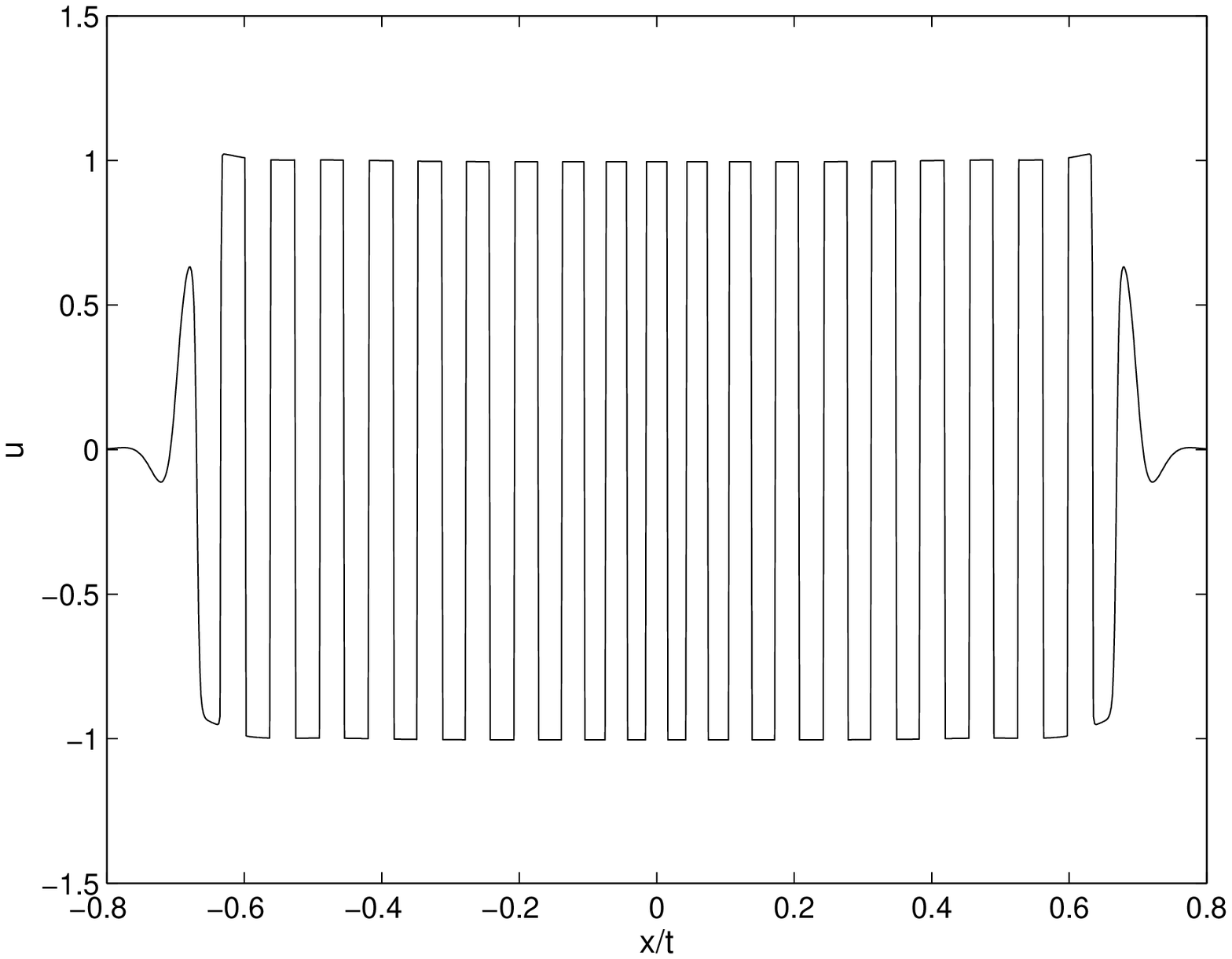}
\label{speed:cubic}}
}
\quad
\mbox{
\subfigure[Numerical solution for $\phi(u)=-u \, e^{-u^2}$ at $t=58$.]
{\includegraphics[width=0.45\textwidth,height=.35\textwidth]{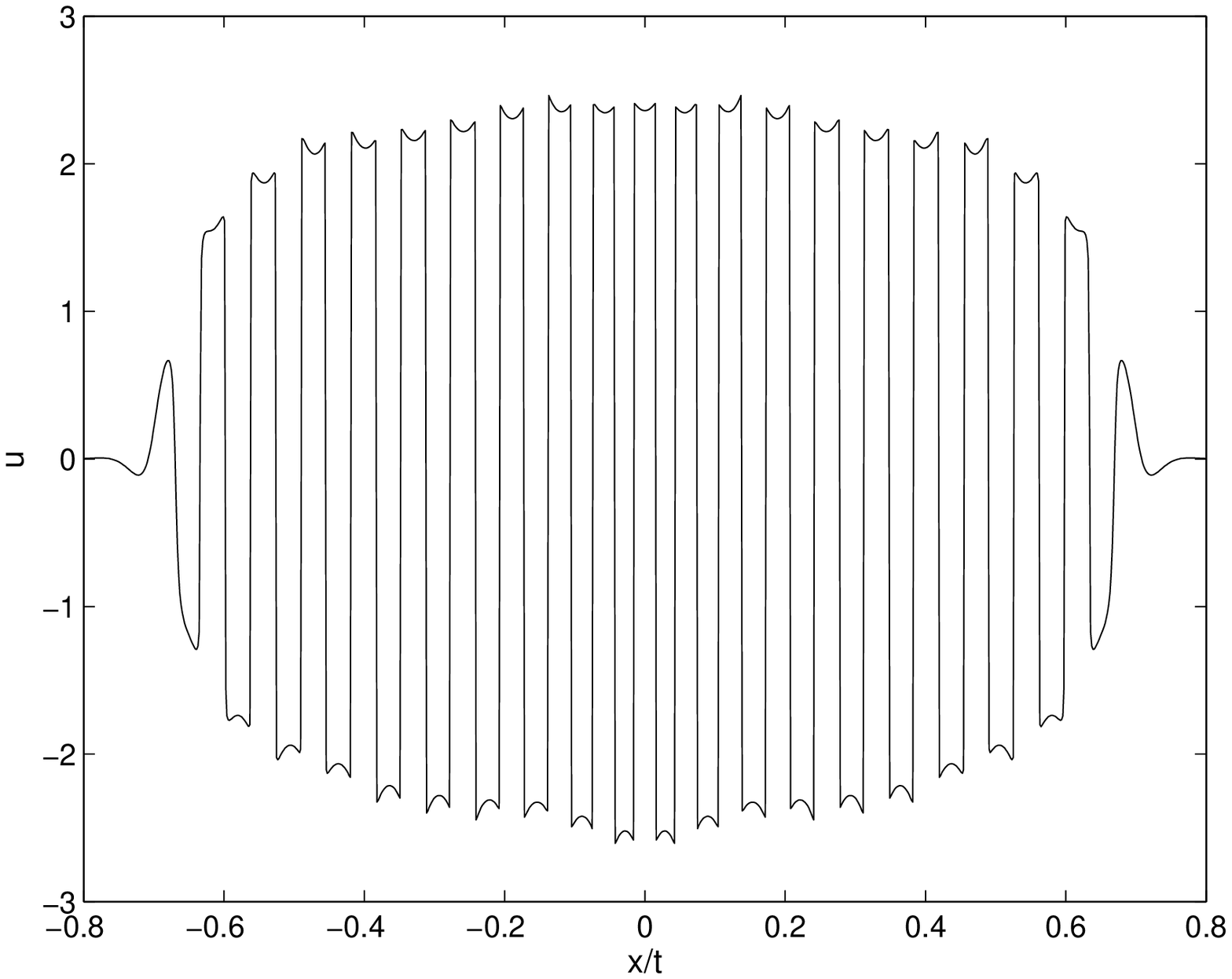}
\label{speed:exposym}}
}
\\
\mbox{
\subfigure[Numerical solution for $\phi(u)=u^3-u$ at $t=68$.]
{\includegraphics[width=0.45\textwidth,height=.35\textwidth]{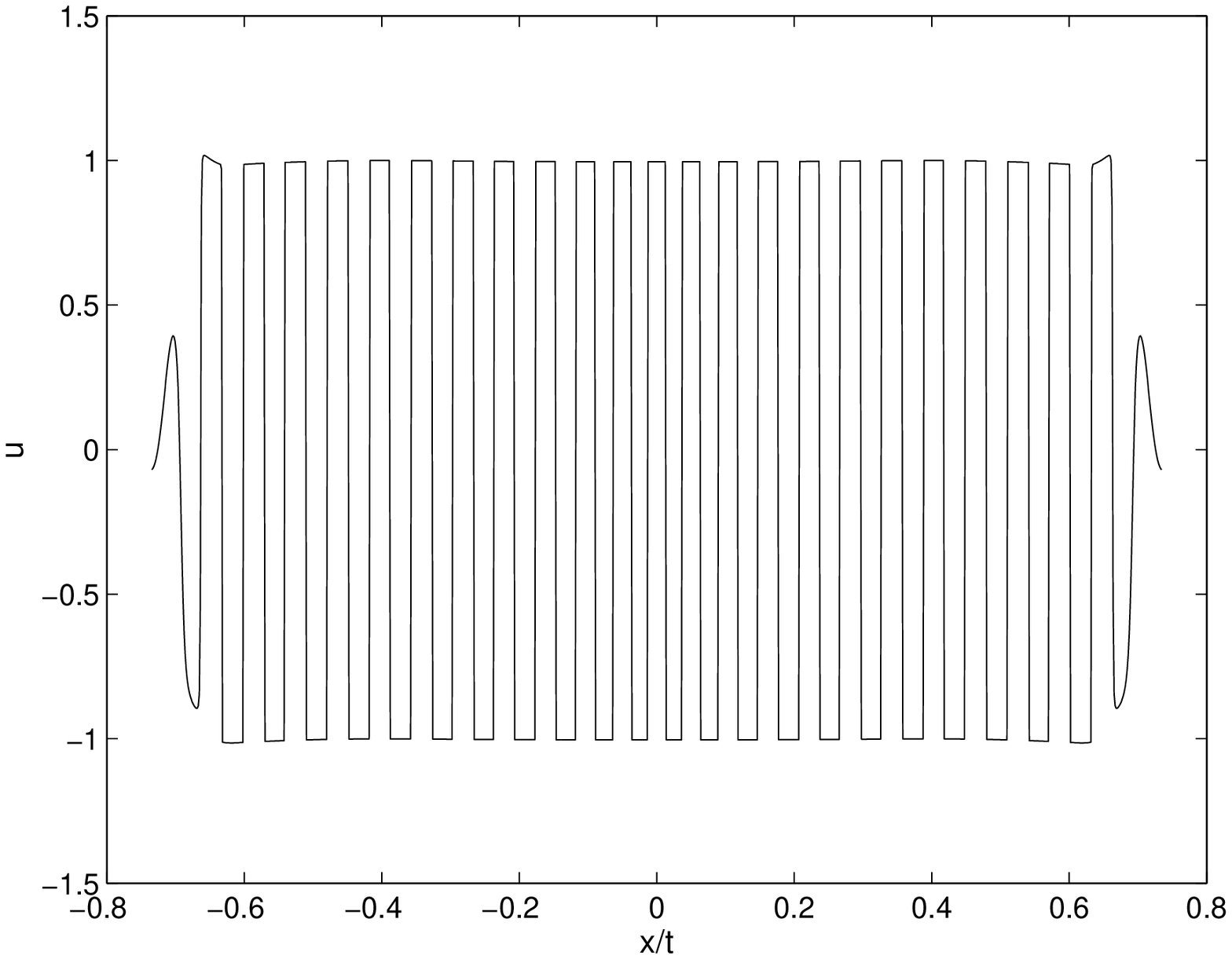}
\label{speed:cubic2} }
} 
\quad
\mbox{
\subfigure[Numerical solution for $\phi(u)=-u \, e^{-u^2}$ at $t=68$.]
{\includegraphics[width=0.45\textwidth,height=.35\textwidth]{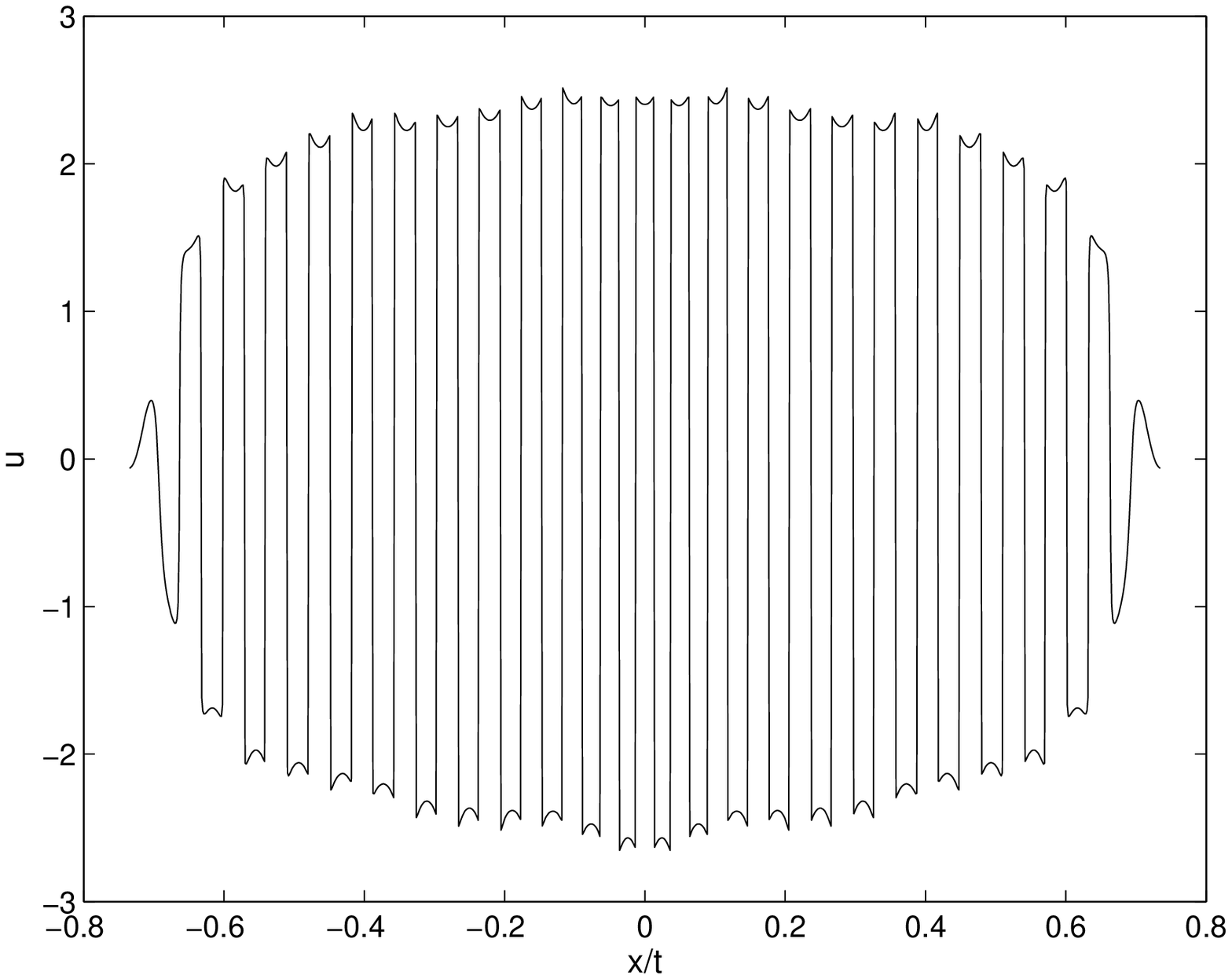}
\label{speed:exposym2}}
}
\caption{Numerical solutions at times $t=58$ and $t=68$ against the velocity variable $\xi=x/t$. The solution is in effect confined in this 
variable to the interval $(-0.8,0.8)$, consistent with the approximated value of $\xi^* (=\xi_0)$ given in (\ref{speed}). 
\subref{speed:cubic} and \subref{speed:cubic2} show the solutions at these times for $\phi(u)= u^3-u$, and \subref{speed:exposym} 
and \subref{speed:exposym2} show the solutions at these times for $\phi(u)= -ue^{-u^2}$.}
\label{unstable:0:speed}
\end{figure}

We first check the prediction for the front speed $\xi^*=\xi_0\approx 0.8$ (see (\ref{speed})). 
In Figure~\ref{unstable:0:speed}, we plot the solution at 
$t=66$ in the {\it velocity} coordinate $\xi=x/t$. This shows that  
the solution profiles are indeed confined in the spatial interval $(-0.8,0.8)$. 

\begin{figure}
\centering
\mbox{
\subfigure[Numerical period for $\phi(u) = u^3 - u$.]
{
\includegraphics[width=0.45\textwidth,height=.35\textwidth]{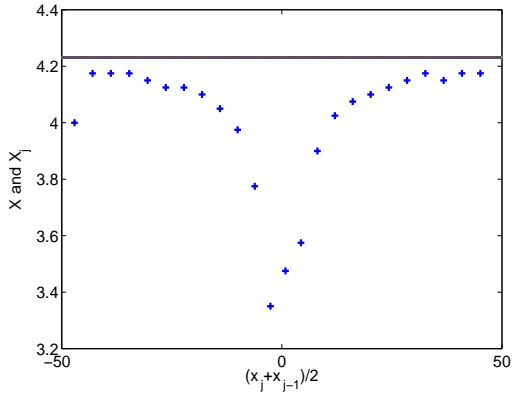}
\label{unstable:0:period1}}
}
\quad
\mbox
{
\subfigure[Numerical period for $\phi(u) = -u \, e^{-u^2}$.]
{
\includegraphics[width=0.45\textwidth,height=.35\textwidth]{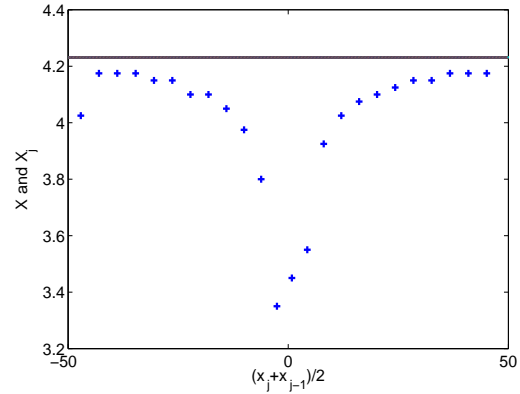}
\label{unstable:0:period2}}
}
\caption{Comparison of the numerical periods $X_j$ and the estimated period $X\approx 4.232$ (solid line). 
The set of values $x_j$ are the grid points where the left upper corner of each plateau is located for the 
numerical solution at $t_f=70$, the numerical periods are $X_j = x_j - x_{j-1}$, here shown against 
$(x_j +x_{j-1})/2$.
}\label{Xperiod}
\end{figure}

In order to check the spatial period (\ref{X:period2}), we let the program run until the front just exits the domain, so that the 
period laid down 
everywhere except in the middle and edges of the domain is that laid down by the advancing modulated travelling waves.
 We then apply the Heaviside function $H$ to the solution at the last time step $t_f$, 
and locate the values $x_i$ of the spatial grid such that $H (u(t_f,x_i))=1$ but $H(u(t_f,x_{i-1}))=0$.  
Let $x_{j}$ denote these values. The distance between two such consecutive values, 
i.e. $x_j-x_{j-1}=:X_j$, should approximate the period (\ref{X:period2}). 
In Figure~\ref{unstable:0:period1} we plot $X_j$ (crosses) and the period $X$ (dotted line) 
against $(x_j +x_{j+1})/2$ (the mean value of two consecutive $x_j$'s). This strongly suggests that the values away from the origin 
(and just before reaching the boundaries) are in agreement with (\ref{X:period2}). It is clear that the grid size imposes a minimum 
error of accuracy. For the numerical simulation in Figure~\ref{unstable:0:period1} the spatial grid is $\Delta x=0.025$, and the best approximation obtained is $4.175$; comparing this with the predicted value (\ref{X:period2}) indicates the significance of the logarithmic correction in (\ref{logcorrect}) (which implies $X_j\sim X -\nu \ln(x_j/x_{j-1})<X$ as $t\to+\infty$), the numerical characterisation of which would require a much finer spatial step and a larger domain (we can confirm the that $\nu>0$, however). 

The nature of the filter applied to determine the $x_j$ implies that the results in Figure~\ref{Xperiod} are not quite symmetric. However, if we compute $X_j$ by locating the $x_j$'s such that 
$H (u(t_f,x_i))=1$ but $H(u(t_f,x_{i+1}))=0$ instead, 
we obtain Figure~\ref{Xperiod} reflected about the $y$-axis.

\begin{figure}
\centering
\mbox{
\subfigure[Results for $\phi(u) = u^3 - u$.]
{
\includegraphics[width=0.45\textwidth,height=.35\textwidth]{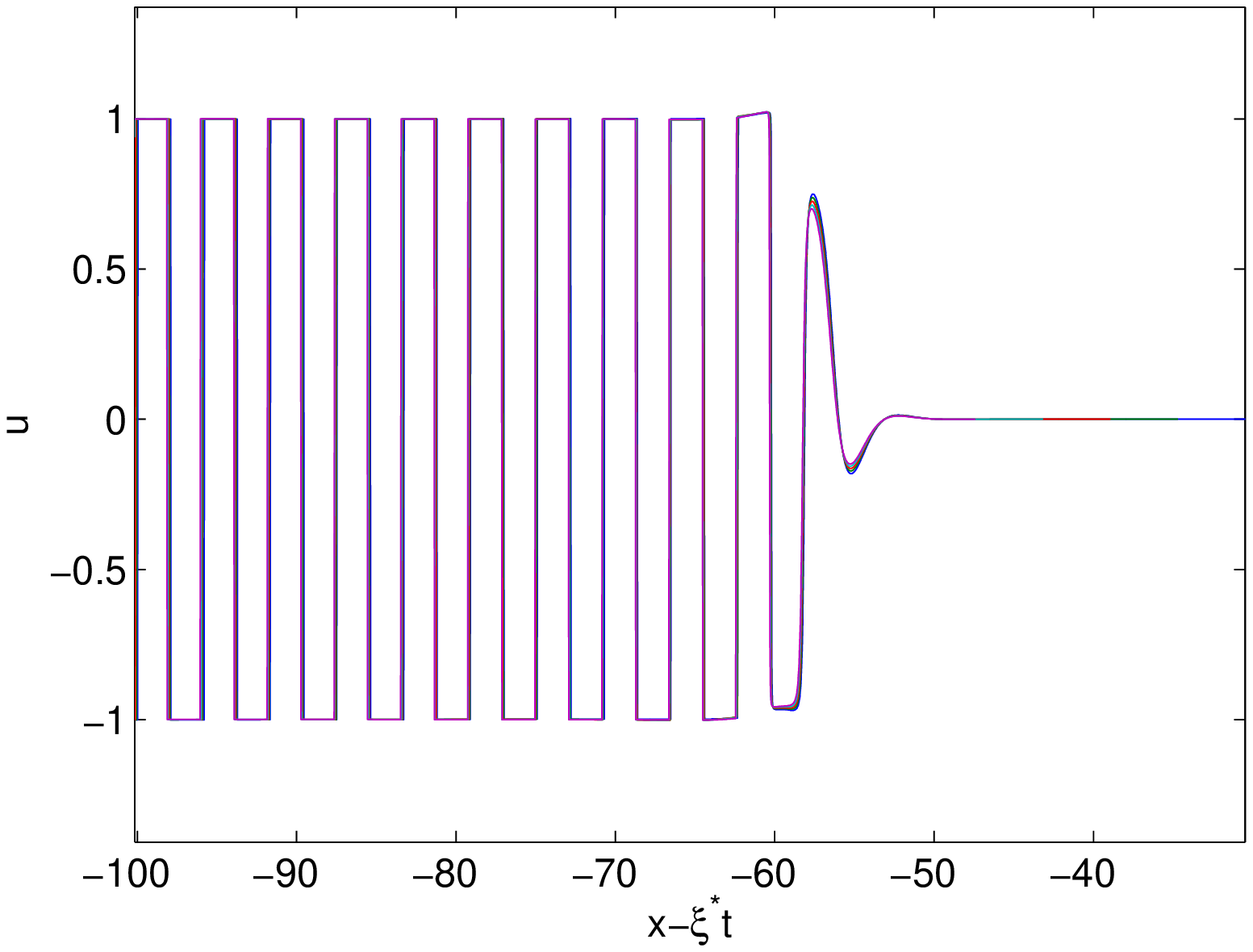}
\label{Tperiod1}}
}
\quad
\mbox
{
\subfigure[Results for $\phi(u) = -u \, e^{-u^2}$.]
{
\includegraphics[width=0.45\textwidth,height=.35\textwidth]{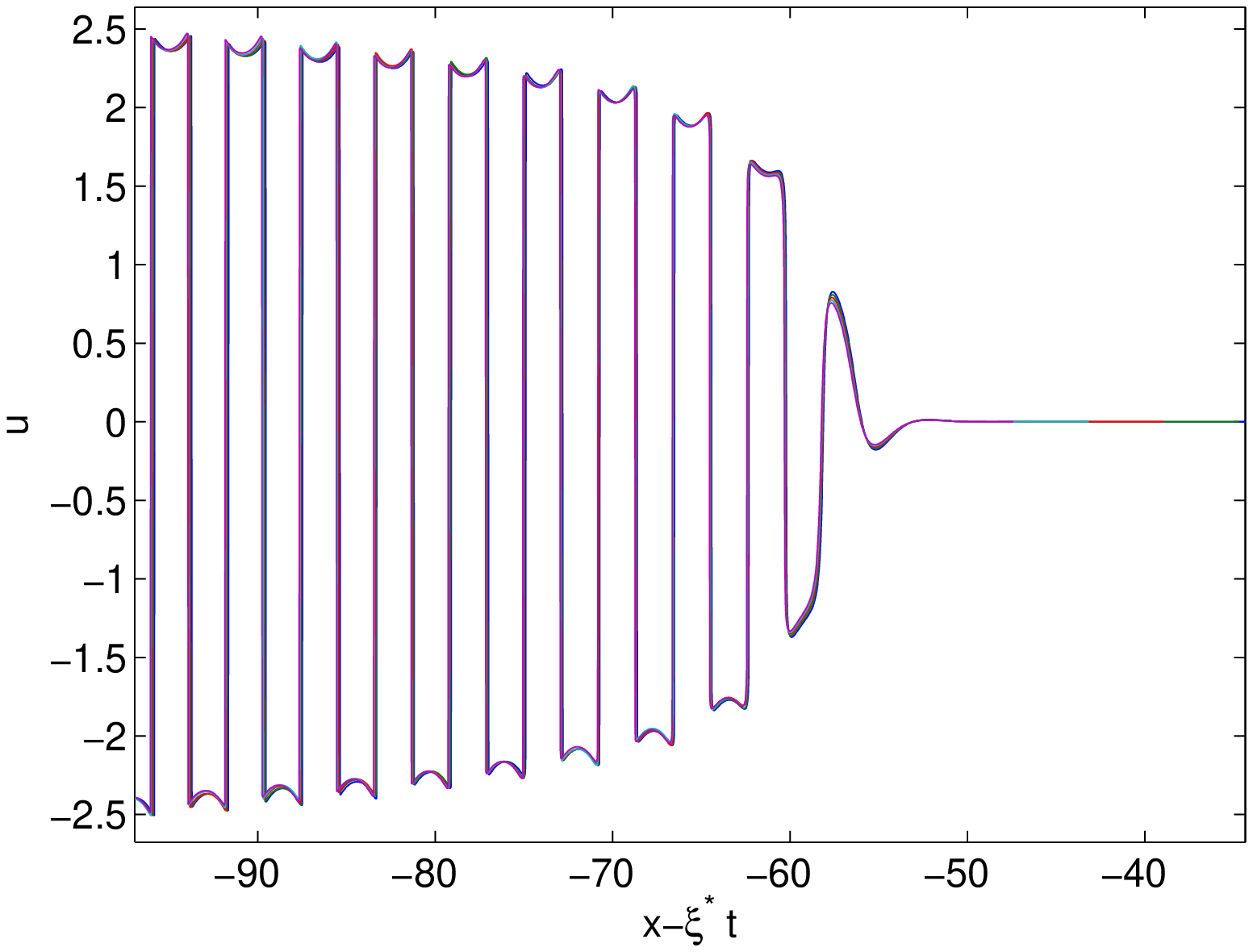}
\label{Tperiod2}}
}
\caption{Time period validation. With a time step size $\Delta t=0.005$ we can verify the approximation of the time period of the modulated travelling wave (\ref{T:period}). 
The pictures show solutions at times $t=5.375 \, k$, with $k=20$, $21$, $22$, $23$ and $24$, against the moving coordinate $z=x-\xi^*t$, in a domain around the front 
(the spatial domain here has $L=100$, however, and $\Delta x=0.05$). \subref{Tperiod1} shows results for (\ref{phi:cubic}) and \subref{Tperiod2} for (\ref{phi:exposym}).
 In both cases the solutions show near overlap around the front.
}\label{Tperiod}
\end{figure}

In order to (approximately) verify the time period of the modulated travelling wave (\ref{T:period}), we take a finer time step, namely $\Delta t= 0.005$, and approximate $T$ numerically by $T\approx 5.375$.  
We then compute solutions at times $t=5.375 \, k$ with $k\in \N$. Figure~\ref{Tperiod} shows results for $k=20$, $21$, $22$, $23$ and $24$, 
where we have plotted the solutions against the moving coordinate $z=x-\xi^*t$ 
(with $\xi^*=\xi_0$ approximated as in (\ref{speed})); \subref{Tperiod1} shows a computation for (\ref{phi:cubic}) and \subref{Tperiod2} a computation for 
(\ref{phi:exposym}). Only part of the domain near the front is shown. Here we have used the spatial step size $\Delta x=0.05$ and the spatial domain has $L=100$.

\paragraph{The case $\phi(u)=u^3-u$ with $u_u=0$.}
As explained above we expect the pattern to approach a steady state, i.e. $u_s$ 
in (\ref{F:left:boundary}) is a piecewise constant function, 
therefore $u_+$ and $u_-$ are constant values and the condition (\ref{t:ode:int}) becomes
\bequ\label{r:l:values}
\phi(u_+) = \phi(u_-) = A \,(=0)\,,
\eequ
and (\ref{symmetric:A}) implies that (\ref{r:l:values}) is fulfilled with
\[
u_- = -1, \ \mbox{and} \  u_+ = 1\,.
\]
\begin{figure}
\centerline{
\includegraphics[width=0.62\textwidth,height=.44\textwidth]{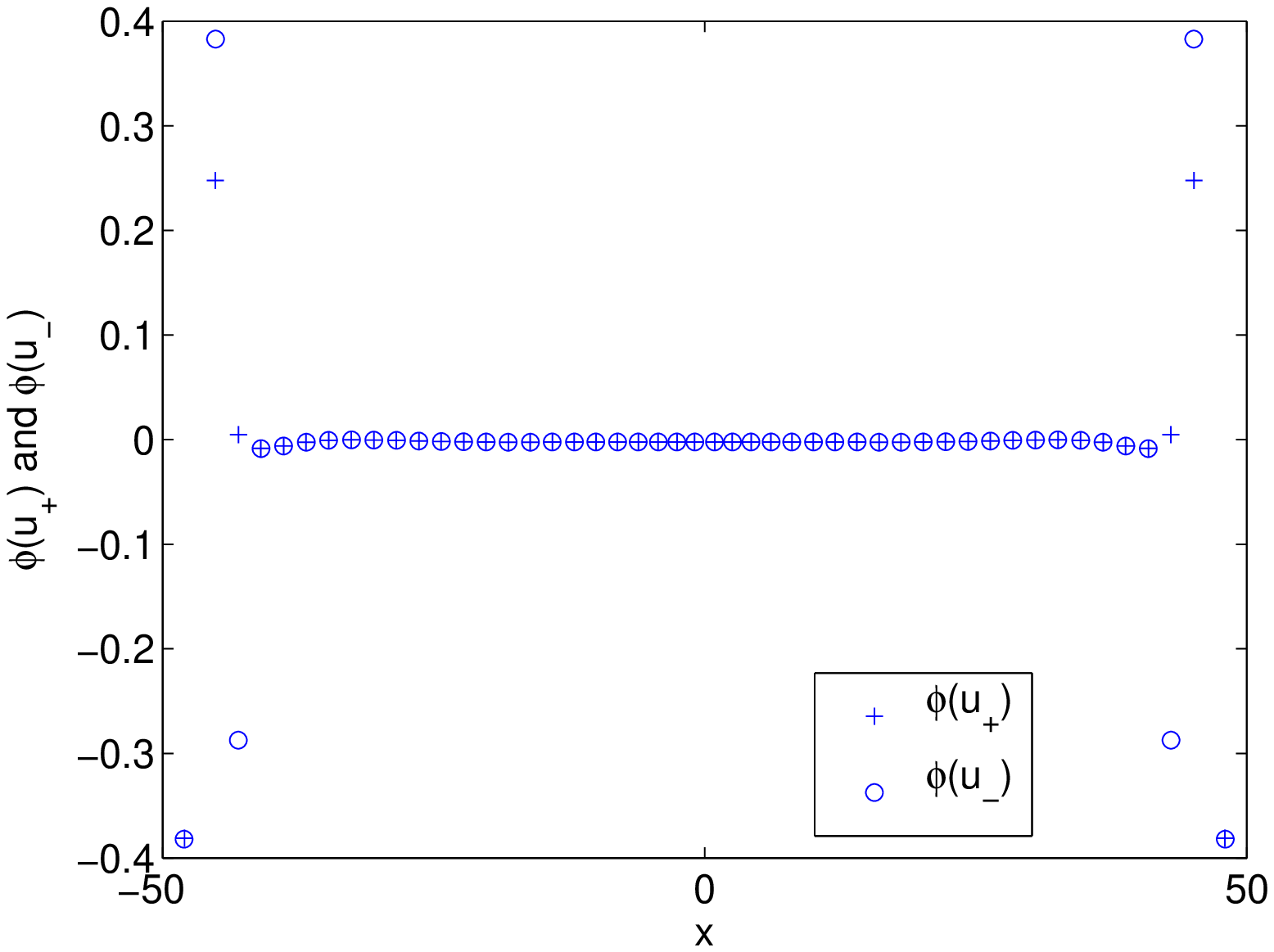}
}
\caption{Comparison of $\phi(u_-)$ and $\phi(u_+)$ at $t_f=70$ for $\phi(u)=u^3-u$.}
\label{phiulur:plot}
\end{figure}
Figure~\ref{phiulur:plot} shows a good agreement of the values $\phi(u_+)$ (circles) and $\phi(u_-)$ 
(crosses) along the domain. These are computed at a late time step, when the fronts have reached the 
boundaries. The values $u_+$ are the values of $u$ at the $x_j$'s obtained above. The values $u_-$ are 
obtained analogously (by applying $-H$ in place of $H$).

\begin{figure}
\centerline{
\includegraphics[width=0.62\textwidth,height=.49\textwidth]{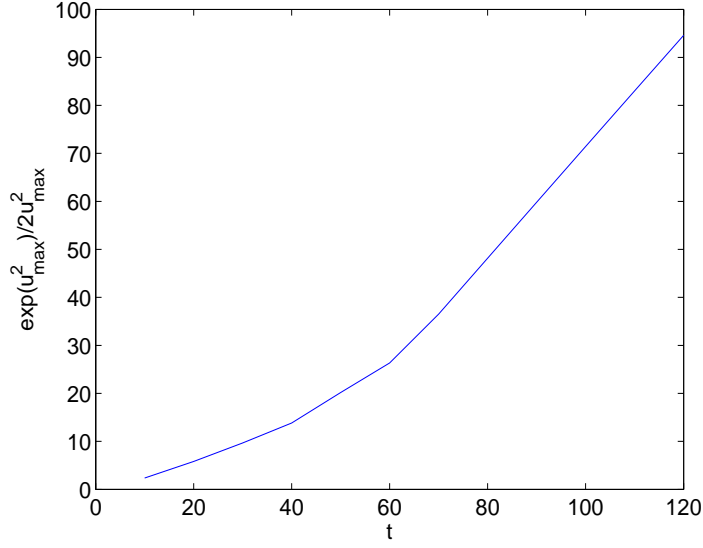}
}
\caption{The value $e^{u_{max}^2}/2u_{max}^2$ against $t$ (see \ref{symmetric:ur1}), 
where $u_{max}$ indicates the maximum of the solution at each time step.}
\label{log:1_2}
\end{figure}

\paragraph{The case $\phi(u)=-ue^{-u^2}$ with $u_u=0$.}
Here (\ref{symmetric:A}) implies that, 
since $\lim_{u\to\pm \infty}\phi(u)=0^\pm$ and there are no values $u$ in the stable region with $\phi(u)=0$, 
$\lim_{t\to +\infty} \phi(u_+(t))=\lim_{t\to +\infty} \phi(u_-(t))=0$. The symmetry of $\phi$ (odd) then suggests
that $u_+(t)=-u_-(t)$ as $t\to +\infty$ and the condition (\ref{t:ode:int}) becomes
\bequ\label{symmetric:ur}
\frac{d}{dt} u_+(t) \sim u_+(t) e^{-u_+(t)^2}\,,
\eequ
hence
\bequ\label{symmetric:ur1}
\frac{e^{u_+^2(t)}}{u_+^2(t)}\sim 2 t\,, \quad u_+(t)\sim \sqrt{\ln t} \quad \mbox{as} \ t\to +\infty\,. 
\eequ
To see this numerically, we have checked that (\ref{symmetric:ur1}) is satisfied by the global maximum of the solution; 
it is reasonable to expect that one of the $u_+$'s attains 
the maximum, see Figure~\ref{2unstable:0}, and that the others grow at a comparable rate. Figure~\ref{log:1_2} shows this computation. 
Although the 
graph is not a straight line of slope $1$, a linear fitting shows that its tail ($t>60$) has approximate slope $1.1$ (the correction terms to the first of (\ref{symmetric:ur1}) 
are only logarithmically smaller, so slow approach to the asymptotic behaviour is expected). 
Here $L=120$ and at $t=120$ the fronts have not yet reached the boundaries.

The interior layers are of a similar structure to those discussed above for $\phi(u)=u-u^3$. There is a noteworthy difference in the `plateau' regions between the interior layers, 
however, whereby their approach to spatially uniform is logarithmic rather than exponential in $t$ as $t\to+\infty$: thus if we set
\[
u\sim u_+(t)+ W
\]
we find for large time that 
\[
\frac{d u_+}{dt}\sim \frac{1}{t}\frac{\partial^2 W}{\partial x^2}
\]
(neglecting terms being only logarithmically smaller than these)
and hence 
\bequ\label{inner:uplus}
W\sim -\frac{1}{4\sqrt{\ln t}}(x-x_j)(x_{j+1}-x) \quad \mbox{as} \quad t\to+\infty \quad \mbox{for} \ x_j<x<x_{j+1}\,.
\eequ
We note that different values in the constant of integration in (\ref{symmetric:ur}) lead to $O(1/t\sqrt{\ln t})$ differences in $u_+$ as $t\to+\infty$ and that these deviations are negligible compared to (\ref{inner:uplus}) and that the numerical observations are consistent with the very slow decay in (\ref{inner:uplus}).

\begin{figure}
\centerline{
\includegraphics[width=0.62\textwidth,height=.44\textwidth]{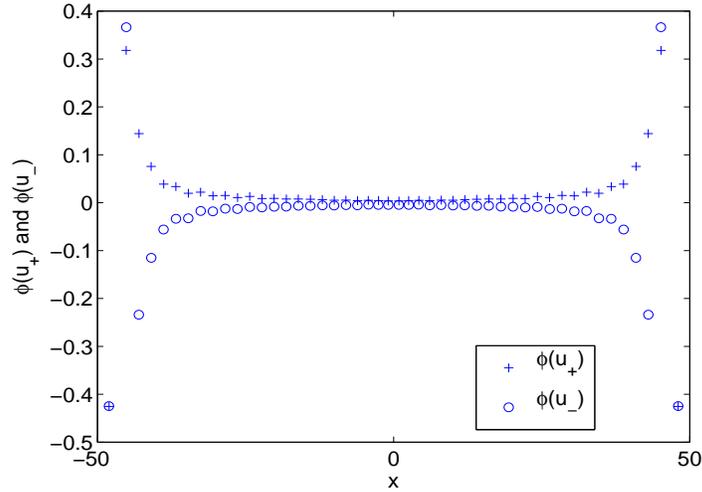}
}
\caption{Comparison of $\phi(u_-)$ and $\phi(u_+)$ at $t_f=70$ for $\phi(u)=-u\, e^{-u^2}$.}
\label{phiulur:plot2}
\end{figure}



\section{Remarks on oscillatory initial conditions}\label{section:new}
In this section we outline the implications of considering initial conditions of the form (\ref{exp:ic}) with $\lambda\in\C$, for brevity setting $\Phi_u=1$. 
Before we continue and to make the exposition clearer, we give some unifying notation that applies to systems exhibiting pulled fronts. 
Let $\bar{\xi}(\lambda)$ denote the function that assigns the linearly selected wave front speed to the exponent $\lambda\in\C$ in (\ref{exp:ic}). 
As before, $\xi^*$ is the critical speed selected by fast decaying initial perturbations (\ref{speed}) and $\lambda^*$ is the associated (real) 
exponential decay rate (\ref{univ:decay}). Let also $\xi_f(\lambda)$ be the function that results from applying the neither-growth-nor-decay condition applied 
to the separation-of-variables solution (\ref{slowdecay}), i.e.
\bequ\label{poss:fast:speed}
\xi_f(\l)=  \mbox{Re}\left(\frac{\l^2}{\l^2-1}\right)/\mbox{Re}(\l)\,.
\eequ 

For comparison purposes we describe the front speed selected $\bar{\xi}(\lambda)$ associated to Fisher's equation (see \cite{PhDWood} and \cite{WKC}). In this case $\xi_f(\lambda)$ with $\l\in\C$ is given by
\bequ\label{xif:fisher}
\xi_f(\l)=\frac{\mbox{Re}(\l^2+1)}{\mbox{Re}(\l)} \quad \l\in\C\,,
\eequ
and $\bar{\xi}(\l)=\xi_f(\l)$ if $\mbox{Re}(\l)<1$ and if $\xi_f(\l)>2$, 
$\bar{\xi}(\l)=2$ otherwise. The analysis for $\l\in\mathbb{R}$ applying in the 
connected components of the complex plane where $\xi_f(\l)>2$.

\begin{figure}
\centerline{
\includegraphics[width=0.62\textwidth,height=.44\textwidth]{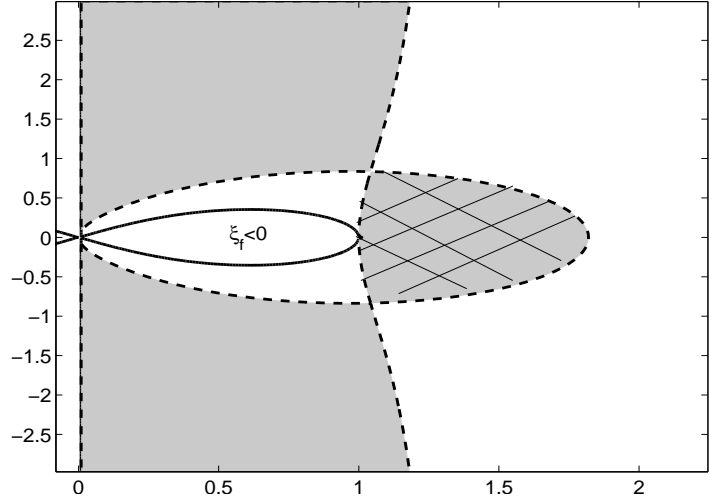}
}
\caption{Contour plot for $\xi_f(\l)$, the possible wave speed for fronts for exponentially decaying initial data. The shaded regions correspond to $\l\in\C$ with $\xi_f(\l)>\xi^*$, 
the dashed contour line shows $\xi_f(\l)=\xi^*$, and the solid one $\xi_f(\l)=0$ (outside this contour $\xi_f>0$). For initial data 
(\ref{osci:slow:data}) with $\l$ in the meshed region to the right the wave speed selected is nevertheless $\xi^*$ (as it is in the unshaded regions), the analysis for $\l\in\R$ applying here. 
In the other shaded regions, however, front velocities greater than $\xi^*$ are indeed realised.}
\label{pseudo-l-complex}
\end{figure}

Let us analyse the current case. The front selection analysis already reveals the oscillatory nature of the pattern that is laid behind, which is also suggested by (\ref{mass}). The matching into the transition region is thus into a 
modulated travelling wave, see Section~\ref{section:2:3}, while the solutions (\ref{slowdecay}) are not oscillatory for $\l\in\R$. 
A local analysis for $\l\in \C$ reveals that there are regions in the complex plane for which $\xi_f(\l)>\xi^*$. 
It is easily verified that $\xi_f(p^*)=\xi^*$ 
In particular this computation gives two extrema for $\xi_f(\l)$ that are attained at
complex conjugates values of $\l$, their approximate values being given in (\ref{the:saddle}).
Inspection of the equation $\xi_f(\l)=\xi^*$, that after setting $\l=\alpha+ i\beta$ can be written as
\begin{equation}\label{contour}
\frac{(\alpha^2+\beta^2)^2-(\alpha^2-\beta^2)}{\alpha\left(  (\alpha^2+\beta^2)^2 - 2(\alpha^2-\beta^2)+1   \right)}=\xi^*\,,
\end{equation}
by using (\ref{poss:fast:speed}) and that $\xi^*$ is obtained at a saddle point of $F$, shows that there are three connected 
components of $\l\in\C$ in which $\xi_f(\l)>\xi^*$. One of them intersects the real axis, let it be denoted by $\Omega_r$. 
We let $\Omega_l^1$ (for $\l$ with $\beta>0$) and $\Omega_l^2$ (for $\l$ with $\beta<0$) denote the other two components. Notice
that these are mirror images in the real axis, by the symmetry of $\xi_f$. A contour plot of (\ref{poss:fast:speed}) is shown in 
Figure~\ref{pseudo-l-complex}: the components $\Omega_r$, $\Omega_l^1$, and $\Omega_l^2$ are those shaded, $\Omega_r$ being 
highlighted with a mesh. An analysis of the type pursued in Section~\ref{sec:expdecay} for $\l\in \R$ applies (e.g. by continuity) to the connected 
components intersecting the real axis. We infer that for initial data of the form
\bequ\label{osci:slow:data}
u_0(x) \sim \e \, e^{-\alpha|x|}\cos(\beta(x+x_0))\,,\quad 0<\e\ll 1
\eequ
for constant $x_0$, then with $\l=\alpha+i \beta$ in $\Omega_r$ the wave speed is $\bar{\xi}(\lambda)=\xi^*$, 
whereas for $\l$ in $\Omega_l^1$ or in $\Omega_l^2$ a modulated travelling with speed 
$\bar{\xi}(\lambda)=\xi_f(\lambda)>\xi^*$ ensues.

Another interesting observation that emerges from Figure~\ref{pseudo-l-complex} is that there 
are values of $\l$ in $\Omega_l^1$ and $\Omega_l^2$ 
(i.e. with $\bar{\xi}(\l)=\xi_f(\lambda)>\xi^*$) that have $\alpha>\l^*$. We verify this with an example below; the maximum of $\alpha$ in the sets $\Omega_l^1$ and $\Omega_l^2$ asymptotes to $1/\xi^*$ ($\approx 1.269$) as $|\beta|\to\infty$ and is larger than $\l^*$, see (\ref{contour}) and Figure~\ref{pseudo-l-complex}. This gives a different scenario than that of the paradigm Fisher's case. 

\begin{figure}
\centering
\mbox{
\subfigure[Results for $\phi(u) = u^3 - u$.]
{
\includegraphics[width=0.45\textwidth,height=.35\textwidth]{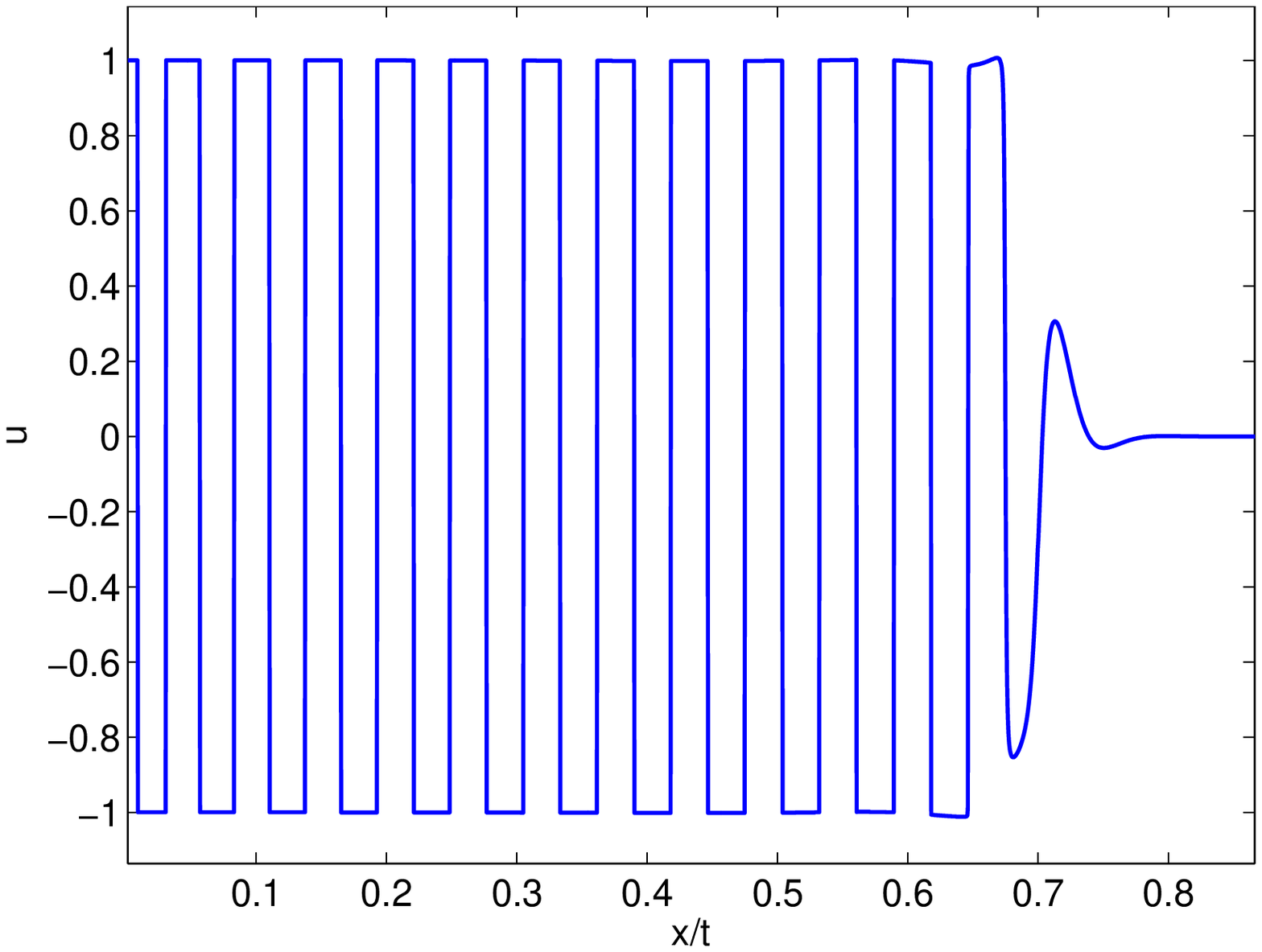}
\label{osci:crit:cubic}
}
}
\quad
\mbox
{
\subfigure[Results for $\phi(u) = -u \, e^{-u^2}$.]
{
\includegraphics[width=0.45\textwidth,height=.35\textwidth]{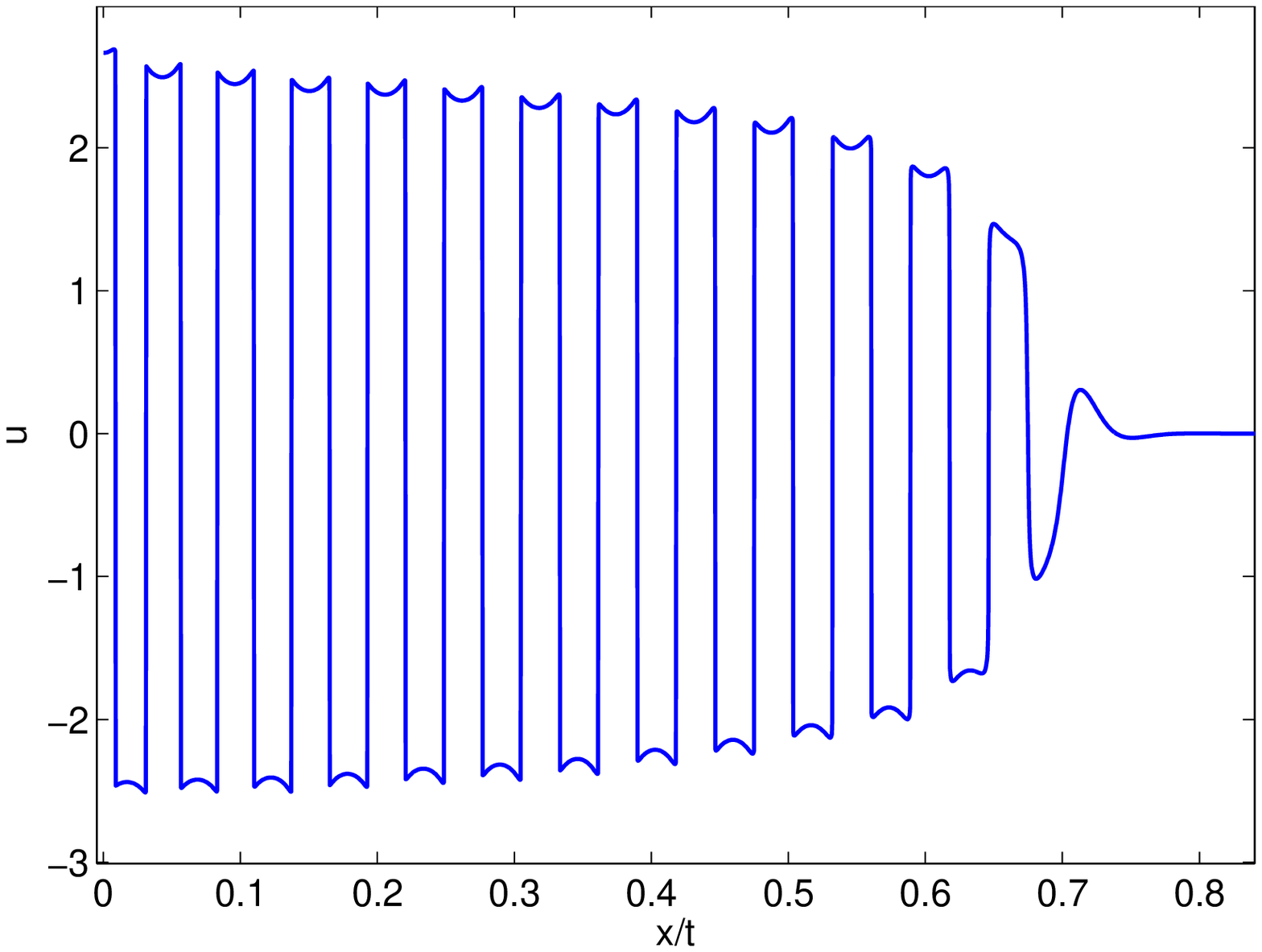}
\label{osci:crit:exp}}
}
\caption{
Simulations for initial data of the form (\ref{exp:ic}) with $\l=1.5+0.5\,i\in\Omega_r$. The figures show the solution at time 
$t=74$ against the velocity variable $x/t$, the solutions up to the front being confined to the interval $[0,0.8]$. 
\subref{osci:crit:cubic} Shows the result for $\phi$ given by (\ref{phi:cubic}) and \subref{osci:crit:exp} the one with $\phi$ given by (\ref{phi:exposym}).
 The computations are done with the domain $[0,100]$ taking $\Delta t= 0.001$, $\Delta x=0.05$. 
}\label{osci:crit}
\end{figure}

\begin{figure}
\centering
\mbox{
\subfigure[Results for $\phi(u) = u^3 - u$.]
{
\includegraphics[width=0.45\textwidth,height=.35\textwidth]{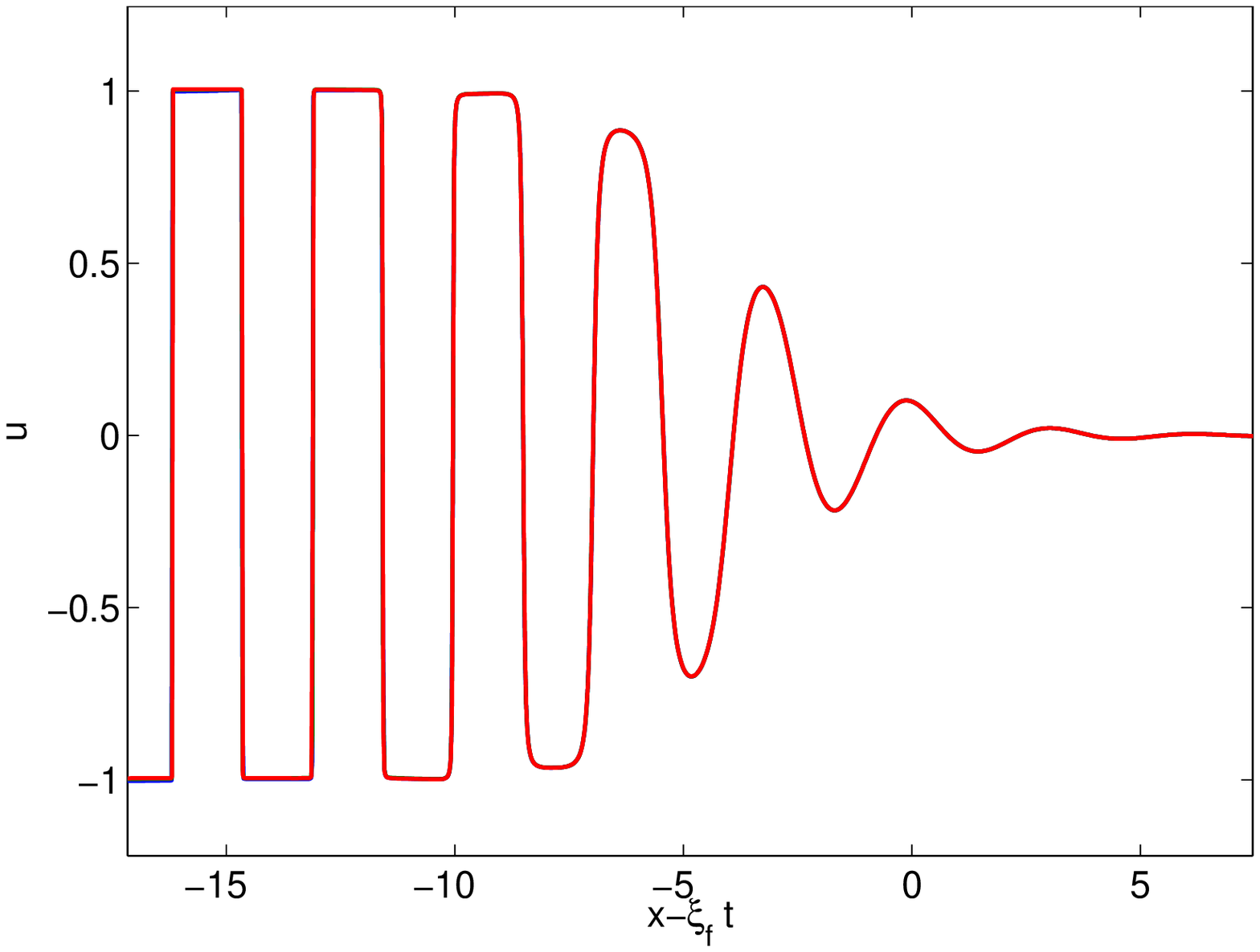}
\label{cubicfastMTWperiod}}
}
\quad
\mbox
{
\subfigure[Results for $\phi(u) = -u \, e^{-u^2}$.]
{
\includegraphics[width=0.45\textwidth,height=.35\textwidth]{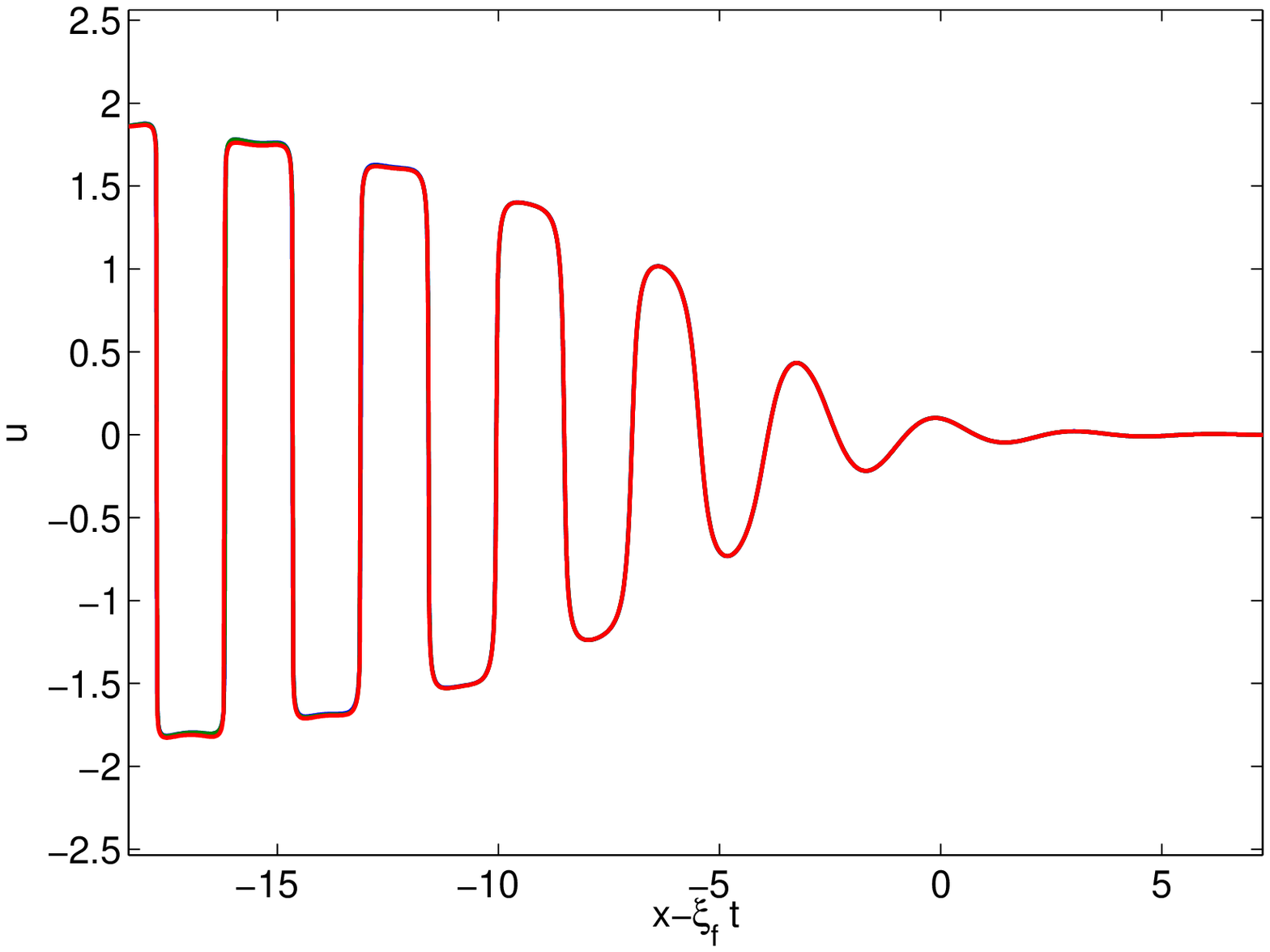}
\label{expofastMTWperiod}}
}
\caption{Time period validation for solutions with initial data (\ref{osci:slow:data}) and $\l=0.5+2i\in\Omega_l^1$. With a time step size 
$\Delta t=0.001$ we can verify the approximation of the time period of the modulated travelling wave (\ref{T_f:explicit}). 
The pictures show solutions at times $t=1.870 \, k$, with $k=12$, $13$, $14$ against the moving coordinate $z=x-\xi_f \,t$, 
in a domain around the front (the spatial domain is $[0,100]$ and $\Delta x=0.05$). Here $\xi_f$ is approximated as in (\ref{T_f:explicit}). 
\subref{cubicfastMTWperiod} shows results for (\ref{phi:cubic}) and \subref{expofastMTWperiod} for (\ref{phi:exposym}).
 In both cases the solutions show near overlap around the front.
}\label{fastMTWperiod}
\end{figure}

\begin{figure}
\centering
\mbox{
\subfigure[Numerical period for $\phi(u) = u^3 - u$.]
{
\includegraphics[width=0.45\textwidth,height=.35\textwidth]{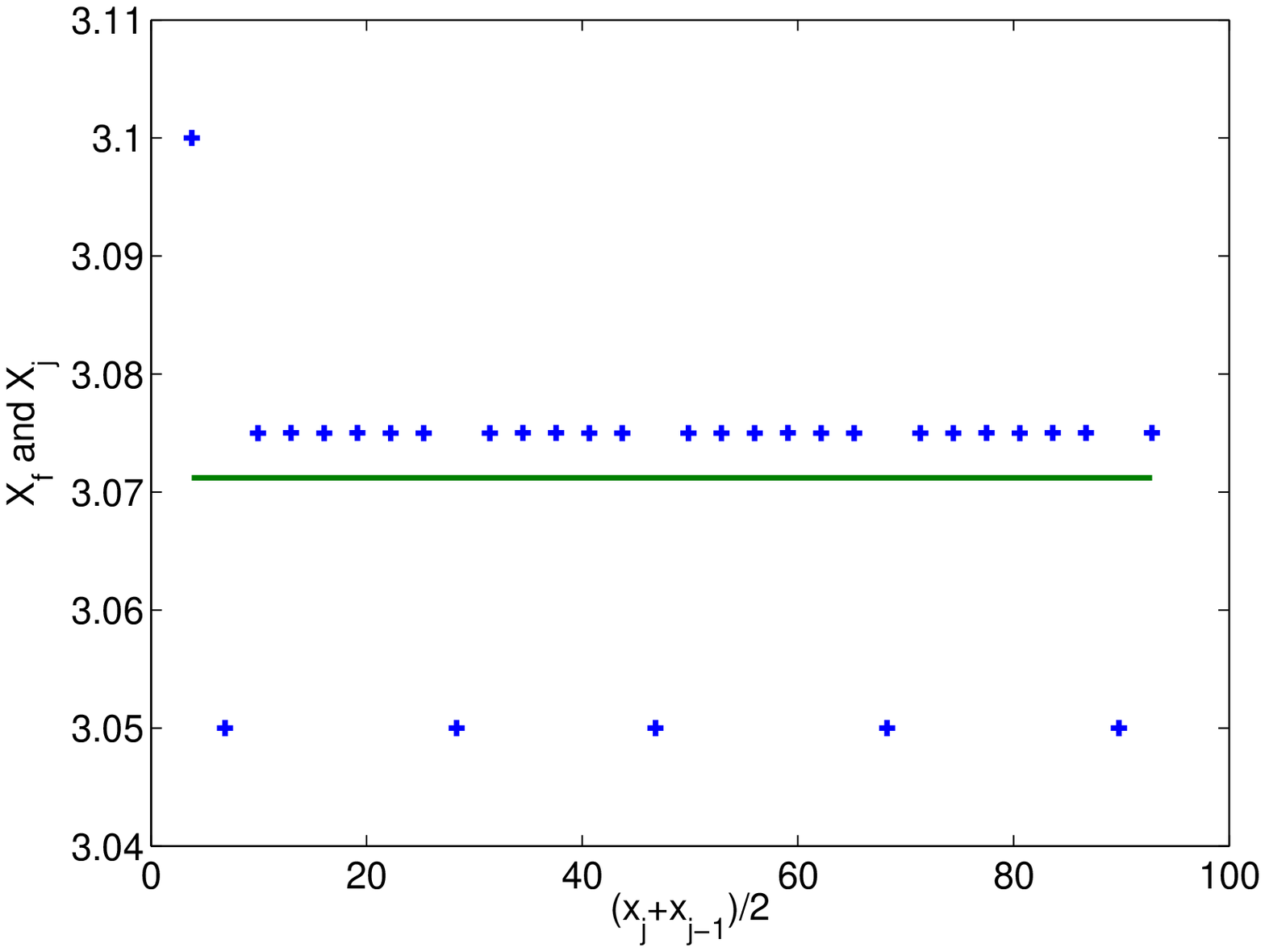}
\label{unstable:f:period1}}
}
\quad
\mbox
{
\subfigure[Numerical period for $\phi(u) = -u \, e^{-u^2}$.]
{
\includegraphics[width=0.45\textwidth,height=.35\textwidth]{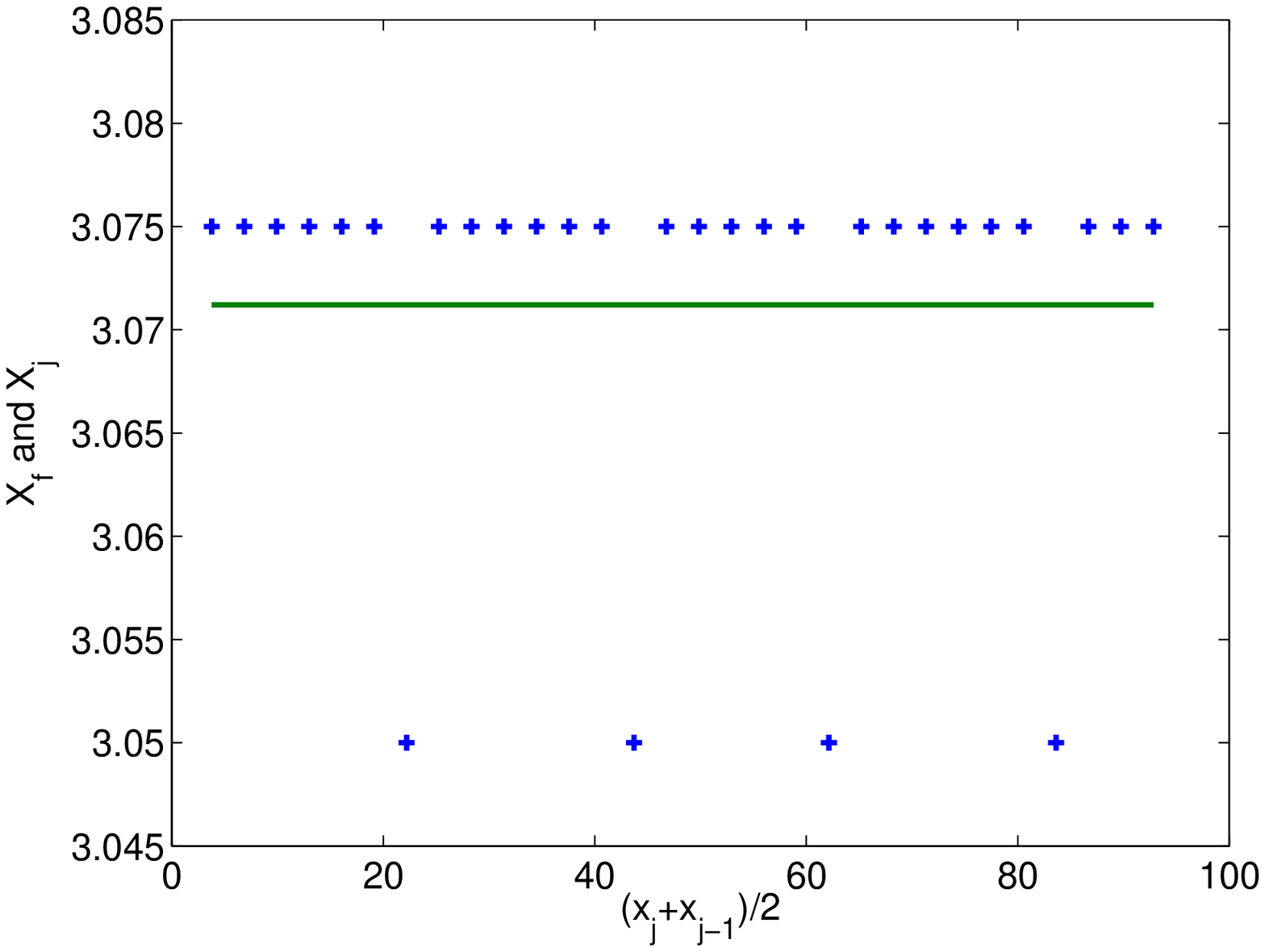}
\label{unstable:f:period2}}
}
\caption{Comparison of the numerical periods $X_j$ and the estimated period $X_f\approx 3.0712$ (solid line) for solutions with 
initial condition (\ref{osci:slow:data}) and $\l=0.5+2i\in\Omega_l^1$ . 
The set of values $x_j$ are the grid points where the left upper corner of each plateau is located for the 
numerical solution at $t=52$, the numerical periods are $X_j = x_j - x_{j-1}$, here shown against 
$(x_j +x_{j-1})/2$. Note the small scale of the vertical axis.
}\label{Xfperiod}
\end{figure}

The leading-order behaviour near a front that propagates with speed $\xi_f(\l)$ corresponds to a solution (\ref{slowdecay}), thus
\[
v(x,t)\sim e^{-(x -\xi_f t)\l } e^{-\left(\xi_f\, \l+ \,\frac{\l^2}{1-\l^2} \right)t }  \quad \mbox{as}\quad t\to+\infty\,,
 \ \frac{x}{t}=O(1)\,,
\]
which has decay rate $\mbox{Re}(\l)$ in $z=x -\xi_f t$ and is periodic in $t$ with period
\[
T_f = \frac{2\pi }{\mbox{Im}\left(\xi_f\, \l +  \frac{\l^2}{1-\l^2}\right)  }
\] 
and so a modulated travelling wave in the transition region would have temporal period $T_f$ and the pattern laid behind would have spatial period
\[
X_f = \xi_f T_f\,.
\]
We now verify these conjectures numerically. First take $\l\in \Omega_r$ not real, namely $\l=1.5+0.5i$. In this case the front propagates with speed $\xi^*$ and not $\xi_f\approx 0.8718$. Simulations for 
(\ref{phi:cubic}) and (\ref{phi:exposym}) are shown in Figure~\ref{osci:crit}.
The solution at time $t=74$ is depicted against the velocity 
variable $x/t$, the solution up to the front being in effect confined in the interval $(0,0.8)$ ($\xi^*\approx 0.787$).
We now take $\l=0.5+ 2i\in \Omega_l^1$. In this case we obtain 
\bequ\label{T_f:explicit}
\xi_f\approx 1.6424 \,,\ T_f\approx 1.8700\,,\ X_f\approx 3.0712\,.
\eequ
We have verified this for both nonlinearities (\ref{phi:cubic}) and (\ref{phi:exposym}), the numerical results being depicted in Figure~\ref{fastMTWperiod}. 
We compute solutions at consecutive times $t=1.87 \, k$ with $k\in \N $ (approximating $T_f$ by $1.87$) . 
Figure~\ref{fastMTWperiod} shows results for $k=12$, $13$ and $14$, where we have plotted the solutions against the moving coordinate $z=x-\xi_f t$ 
for $\xi_f$ approximated as in (\ref{T_f:explicit}). Only a domain near the front is shown. The value of $X_f$ can be verified as in 
Section~\ref{section:3}, cf. Figure~\ref{Xperiod}, and the results for both $\phi$'s are shown in Figure~\ref{Xfperiod}.

\begin{figure}
\centering
\mbox{
\subfigure[Moving frame and time period.]
{
\includegraphics[width=0.45\textwidth,height=.35\textwidth]{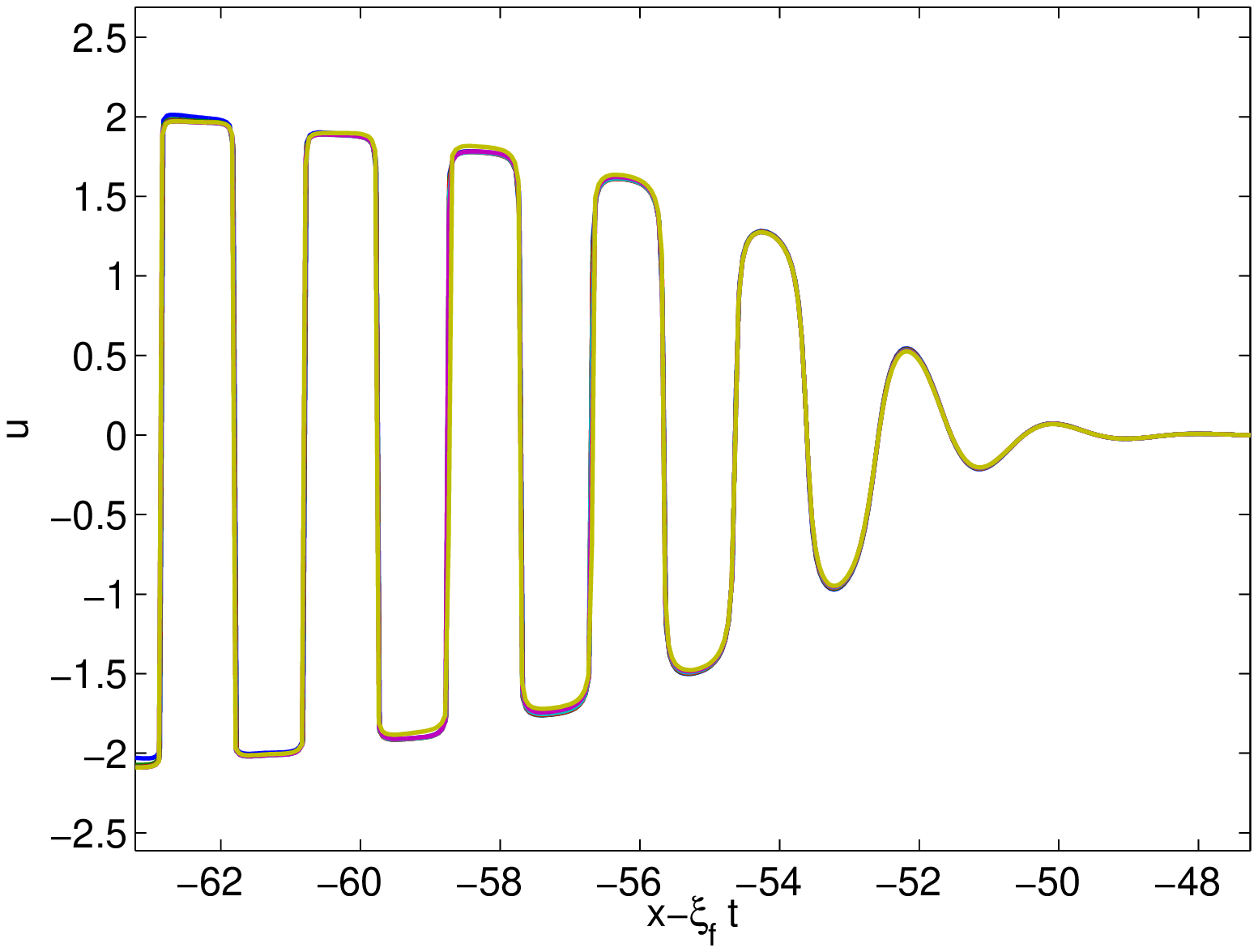}
\label{mtw:fasterlambdastar}
}
}
\quad
\mbox
{
\subfigure[Pattern wavelength.]
{
\includegraphics[width=0.45\textwidth,height=.35\textwidth]{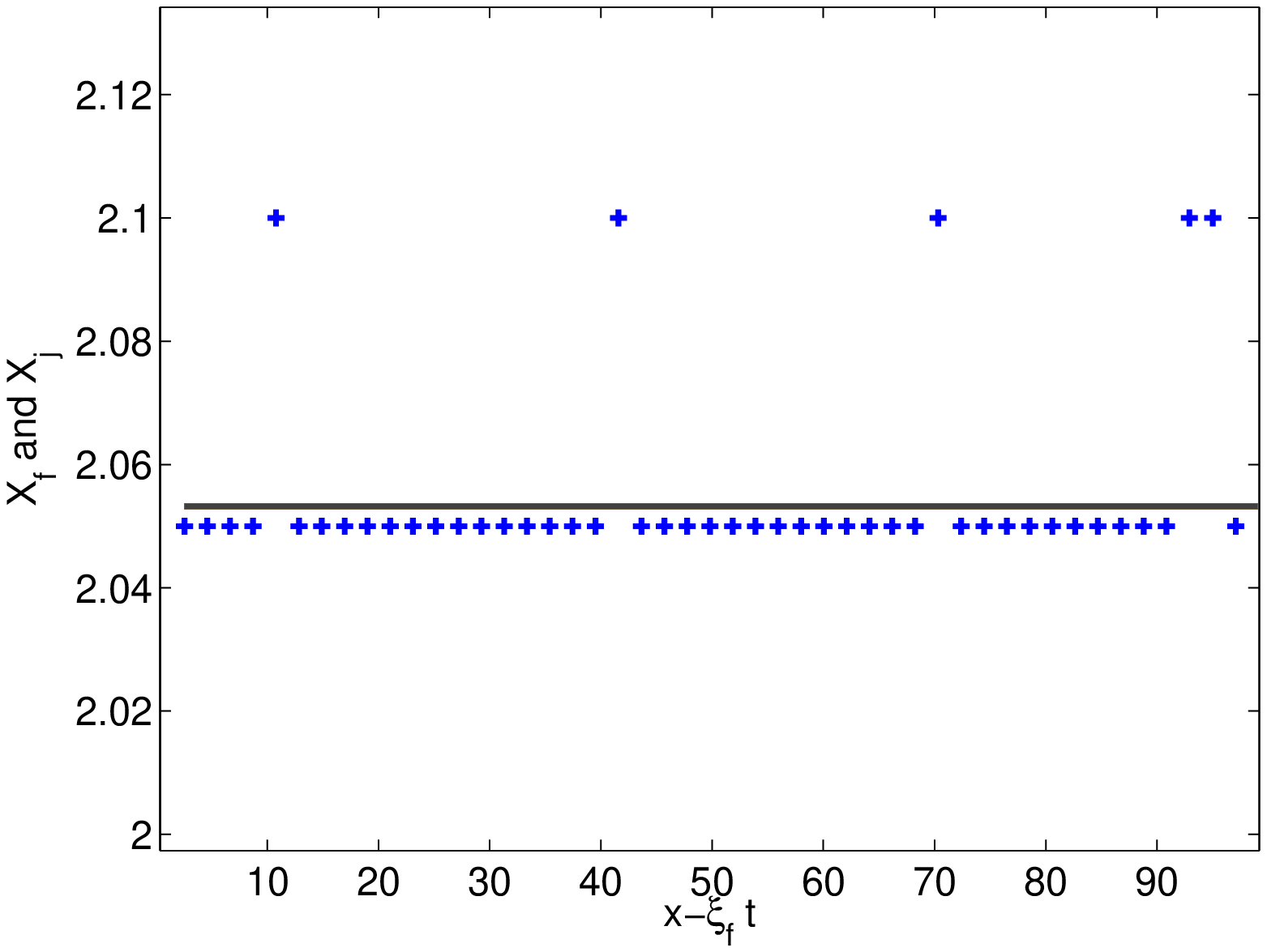}
\label{period:fasterlambdastar}}
}
\caption{
Example with $\mbox{Re}(\l)>\l^*$ and $\xi(\l)>\xi^*$. Simulations for initial data of the form (\ref{exp:ic}) with $\l=1.05+3\,i\in\Omega_l^1$
 and $\phi$ given by (\ref{phi:cubic}). 
\subref{mtw:fasterlambdastar} shows results at time steps $t=2.33k$ with $k=35$ to $40$, with increment $1$, in the variable $x-\xi_f t$ where 
$\xi_f$ is approximated as in (\ref{osci:fast}). As before, the solutions are depicted around the front, showing near overlap, 
thus confirming that $\xi(\l)=\xi_f$ in this case, and that the time period of the modulated travelling wave is that given in (\ref{osci:fast}). 
\subref{period:fasterlambdastar} Comparison of the numerical periods $X_j$ and the estimated period $X_f$ (solid line) in (\ref{osci:fast})
for solutions with initial condition (\ref{osci:slow:data}) and $\l=0.5+2i\in\Omega_l^1$. The numerical period is obtained as in 
Section~\ref{section:3}, see also Figure~\ref{Xfperiod}. Here we have computed in the domain $[0,100]$, with spatial step $\Delta x=0.05$ 
and temporal one $\Delta t=0.01$.
}\label{fasterlambdastar}
\end{figure}

Finally, we take $\l=1.05+ 3i\in \Omega_l^1$, so that $\mbox{Re}(\l)>\l^*$, and compute 
\bequ\label{osci:fast}
\xi_f\approx  0.8811\,, \ T_f\approx 2.3303\,,\ X_f\approx\,2.0532.
\eequ 
Clearly $\xi_f>\xi^*$, $T_f<T$ and $X_f<X$. Numerical results confirming that $\bar{\xi}(\l)=\xi_f(\lambda)$ for this value 
of $\lambda$ are shown in Figure~\ref{fasterlambdastar}. 
Figure~\subref{mtw:fasterlambdastar} 
shows solutions in the coordinate $x-\xi_f t$ at times $t=2.33 k$ with $k=35$ to $40$ (approximating $T_f$ by $2.33$), 
and Figure~\ref{period:fasterlambdastar} compares the approximation of 
$X_f$ above to the period found numerically at $t=107.18$ (computed as in Section~\ref{section:3}). We recall that this gives an exception to the behaviour of more widely studied systems for which solutions with initial exponential decay larger than $\l^*$ propagate at a speed $\xi^*$ (cf. \cite{frontsreview}, page 49).

\section{Discussion}\label{sec:discu}
We have presented an analysis that aims to clarify pattern formation in the 
pseudo-parabolic equation (\ref{main:eq}) when an initial disturbance of an unstable state is introduced. 
This and related models, such as the Cahn-Hilliard equation (\ref{cahn-hilliard}), have not been analysed in great detail in this context, 
in contrast to other PDEs, most notably semilinear parabolic ones, cf. \cite{frontsreview} and the references therein. 

We have focused on the odd nonlinearities (\ref{phi:cubic}) and (\ref{phi:exposym}) and taken the unstable state to be $u_u=0$, in what we called the 
{\it symmetric} cases. The analysis implies that there are three asymptotic regimes as $t\to +\infty$ that have to be matched. 
The linear regime ahead of the leading edge of the propagating disturbance, contains information that carries over to the nonlinear regimes behind of it, indicating in particular how pattern formation is initiated by a modulated travelling wave in the subsequent transition regime, 
although the pattern laid behind in the last regime is ultimately shaped by the specific form of the nonlinearity and the conserved moments 
(\ref{mass}). The success of the approach is supported by the numerical results presented in sections \ref{section:3} and \ref{section:new}. 

To deal with the linear regime we have outlined the JWKB approach, that aims to identify
 the exponential contribution that selects the front. For fast decaying initial conditions 
(such as a Gaussian, as outlined in Section~\ref{section:2}) the front location is determined, up to a logarithmic correction, by identifying the 
transition from exponential decay to exponential growth in the solution associated with fast decaying initial conditions of the linearised equation.
The front advances at a linear rate as 
$t\to+\infty$ and has speed $\xi^*$. This is not a universal rule-of-thumb; there are well-known examples of second order semilinear parabolic 
equations for which the wave speed is faster than that given by the linear regime: 
these are nonlinearly selected fronts or of the `pushed' type. In such problems an a priori knowledge of the, usually, a one-parameter family of travelling waves is an advantage; the linear-regime behaviours ahead of the front are matched as appropriate 
back into a front location that travels at the speed of the relevant travelling wave, cf. \cite{frontsreview}, \cite{CK}. 
When the speed of the travelling wave is that given by linear 
arguments, the fronts are often said to be of the `pulled' type. That the fronts are of this type for (\ref{main:eq}) is in fact confirmed numerically 
(sections \ref{section:3} and \ref{section:new}).

The matching into the second transition region is here into a modulated travelling wave 
(travelling at speed $\xi^*$ and periodic in time with period $T$ given in (\ref{T:period}), 
see Section~\ref{section:2:3}). We have mentioned 
that there is an intermediate region giving rise to the logarithmic correction in (\ref{log:correct}) with $x-\xi^*t=O(t^{1/2})$ in which the dominant balance is a complex heat equation (with diffusion constant $D\approx -0.147+0.892i$ in this case). The matching condition in this region to the modulated travelling wave is the asymptotic behaviour  
\[
v(x,t)\sim e^{- p^* Z-F(p^*) t } (B +C\, Z) \quad\mbox{as}\quad Z\to +\infty\,,
\]
with  $Z=x-\xi^*t$ 
for some constants $B$ and $C$ (here $B$ and $C$ result from the repeated root condition, this is equivalent to a 'degenerate node' in the travelling wave case). 
If $C>0$, then  
\bequ\label{dipole:bc}
s(t) \sim x-\xi^* t + \frac{3}{2\l^*} \ln t + O(1) \quad \mbox{as} \quad t\to+\infty\,,
\eequ
corresponding to a dipole solution of the heat equation, 
whereas if $C=0$, 
\bequ\label{error:bc}
s(t) \sim x-\xi^* t + \frac{1}{2\l^*} \ln t + O(1) \quad \mbox{as} \quad t\to+\infty
\eequ
(see, e.g. \cite{frontsreview} and \cite{CK} for details). 
Numerical results have been so far inconclusive as to whether (\ref{dipole:bc}) or (\ref{error:bc}) 
applies in the current case.

In Section~\ref{sec:expdecay} and in Section~\ref{section:new} we have analysed the front speed selection mechanism for slowly decaying initial perturbations.  
The analysis of Stokes lines in Section~\ref{sec:expdecay} is done for real $\lambda$ in (\ref{exp:ic}) and applies to non-symmetric cases as well and can be generalised to $\l\in\C$. 
It shows that exponentially decaying initial perturbations lead to fronts that propagate with the critical speed $\xi^*$ (\ref{speed}) and 
decay exponentially with rate $\lambda^*$ (\ref{univ:decay}). This is supported by numerical results shown in Section~\ref{section:new}. 
For complex $\lambda$ we have discerned the front speed selected for each $\lambda$, by analysing the level sets associated to (\ref{poss:fast:speed}) and assuming that one can extend the results for real $\lambda$ by continuity into the pertinent connected components. In particular, we have found that there are 
regimes of the decay rate ($\mbox{Re}(\lambda)$) and the wavelength ($\mbox{Im}(\lambda)$) for which the front 
propagates at a wave speed faster than $\xi^*$. This is investigated numerically (see Figure~\ref{fasterlambdastar}) and is worth emphasising since it gives a different scenario for front selection 
mechanism than that exhibited by well-studied semilinear reaction-diffusion equations: there are initial conditions with exponential decay faster 
than the critical one for which the front propagates with a speed faster than the critical one.

\begin{figure}
\centering
\mbox{
\subfigure[]
{
\includegraphics[width=0.45\textwidth,height=.35\textwidth]{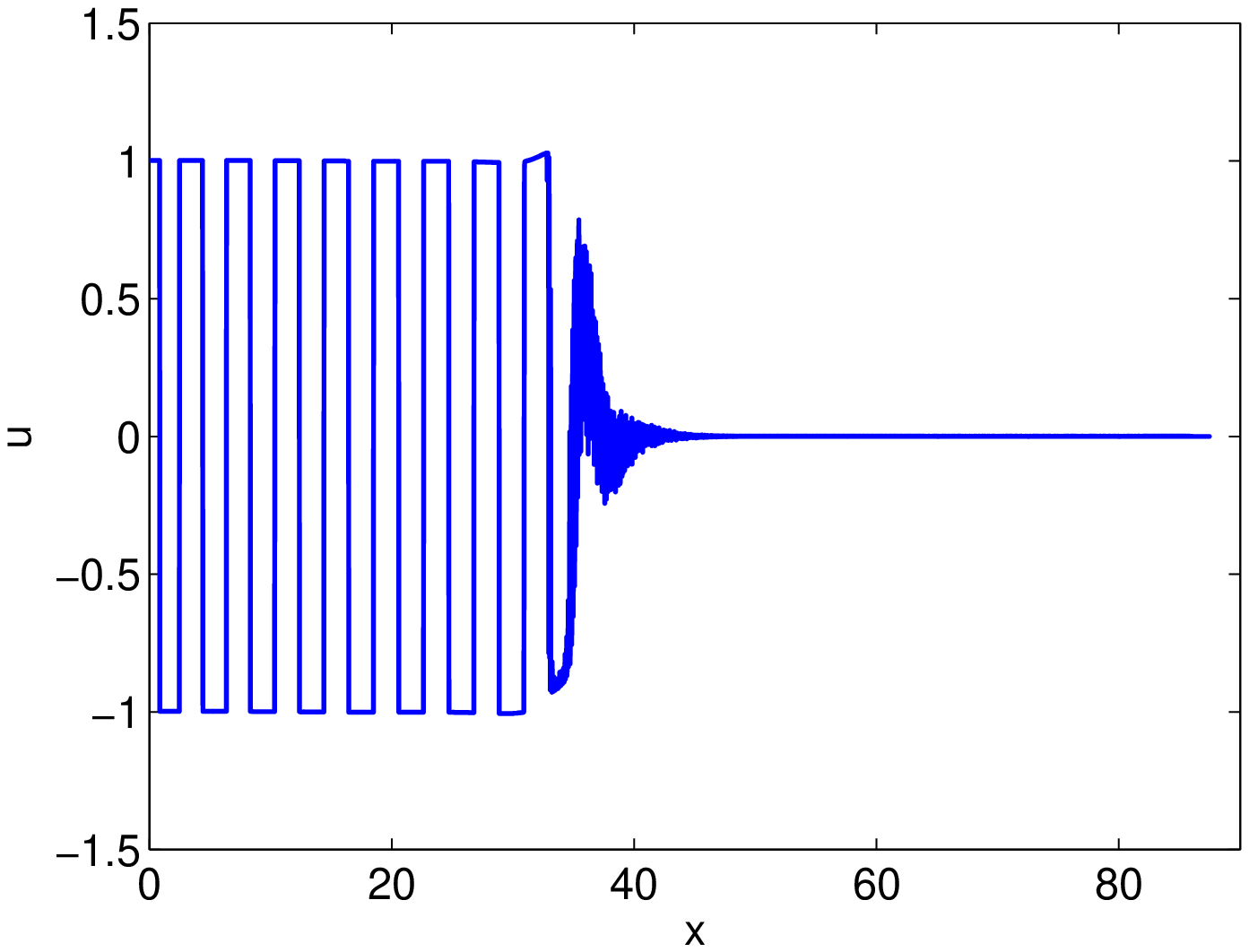}
\label{rosci:plain}}
}
\quad
\mbox{
\subfigure[]
{
\includegraphics[width=0.45\textwidth,height=.35\textwidth]{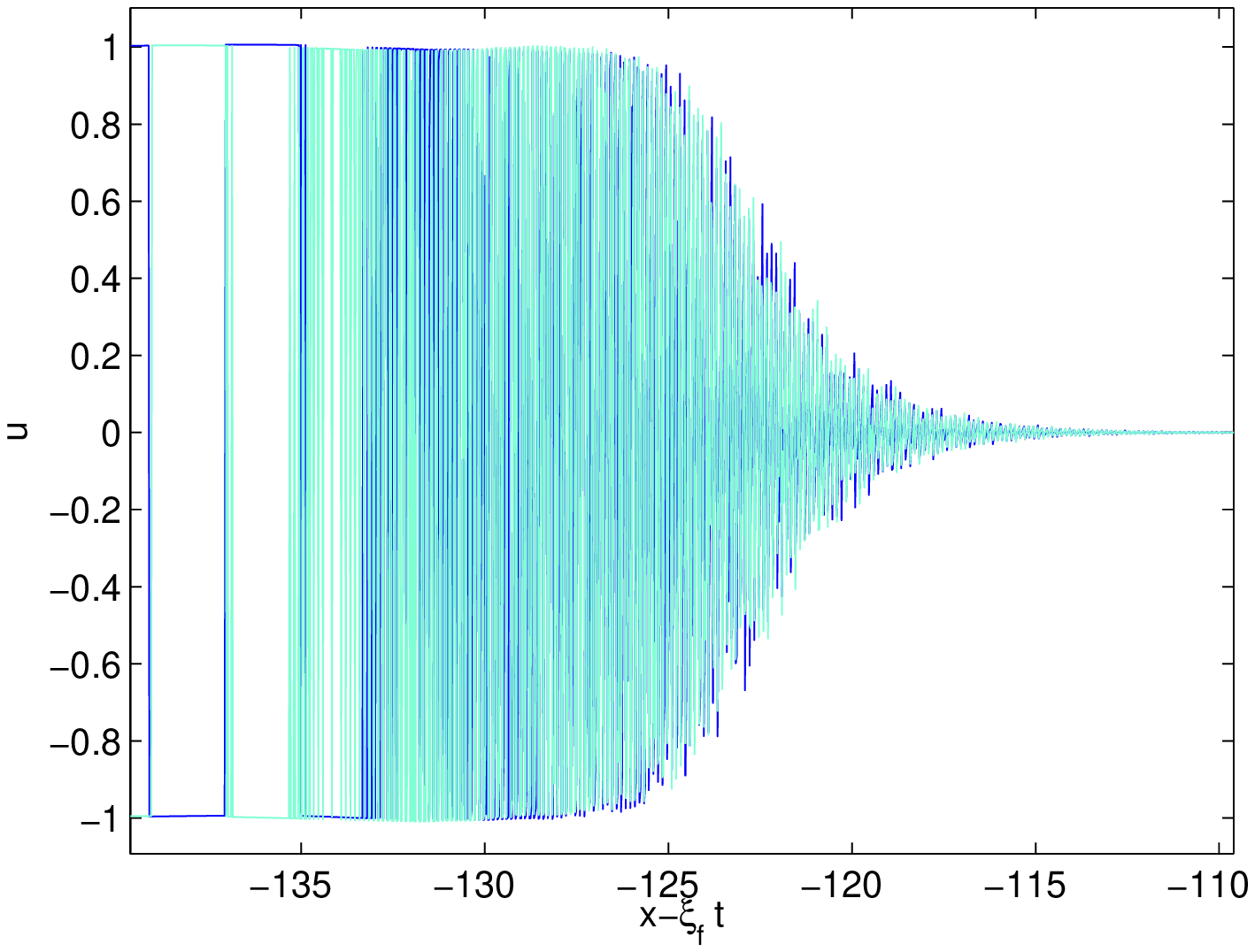}
\label{rosci:speed}}
}
\caption{Numerical computation with $\phi(u)=u^3-u$ for initial data with $\lambda=0.5$, here $\Delta x=0.025$ and $\Delta t=0.01$. \subref{rosci:plain} show the initial development of narrow oscillations around the edge of the front at $t=60$. \subref{rosci:speed} shows the profiles near the front at $t=60$ and at $t=65$ against $x-\xi_f(\alpha +\beta i)\,t$ with $\alpha=0.5$ and $\beta=2\pi/\Delta x$, giving $\xi_f \approx 2$.
}
\label{rapidoscilL05}
\end{figure}

We continue this discussion by mentioning some of the complications that emerge in the numerical simulations. The numerical scheme develops instabilities that seem to be related to the evolution of the leading edge of the front being dominated by a backward heat equation (see Section~\ref{section:2.3new}). This is more apparent for initial perturbations of the form $v_0=e^{-\lambda |x|}$ with $\mbox{Re}(\lambda) <1$. Narrow oscillations of the same order as the spatial step emerge around the leading edge of the front after a number of time iterations. They develop into a faster wave front. We show an example computation with $\phi(u)=u^3-u$ and $\lambda=0.5$ in Figure~\ref{rapidoscilL05}. We have computed numerically the wave speed of such front (by locating the front at several time steps) resulting in approximately $2$ for the spatial steps $\Delta x= 0.1$, $0.05$ and $0.025$. On the other hand, a computation of $\xi_{f}(\alpha+i\beta)$ (see (\ref{poss:fast:speed})) where $\beta=2\pi/\Delta x$ (i.e. the spatial period is given by the grid size) and $\alpha$ being the decay ahead of the front (that in this example is $\alpha\approx 0.5$; remaining close the initially imposed one) yields $\xi_{f}(\alpha+i\beta)\approx 1.9995$ for $\Delta x=0.1$, $1.9998$ for $\Delta x= 0.05$ and $1.9999$ for $\Delta x= 0.025$. We show solutions computed with $\Delta x=0.025$. Figure~\ref{rapidoscilL05}\subref{rosci:plain} shows the appearance of the narrow oscillations and Figure~\ref{rapidoscilL05}\subref{rosci:speed} shows later profiles against the moving coordinate with speed $\xi_f$, showing near overlap. This behaviour can be checked numerically for the linearised problem, leading to the same results for the same choice of parameters, numerical spatial and temporal steps. This suggests that the numerical instability leads to spurious solutions propagating with speed $\xi_f(\l)>\xi^*$.

\begin{figure}
\centering
\mbox{
\subfigure[$u$ against $x$ at $t=140$.]
{
\includegraphics[width=0.45\textwidth,height=.35\textwidth]{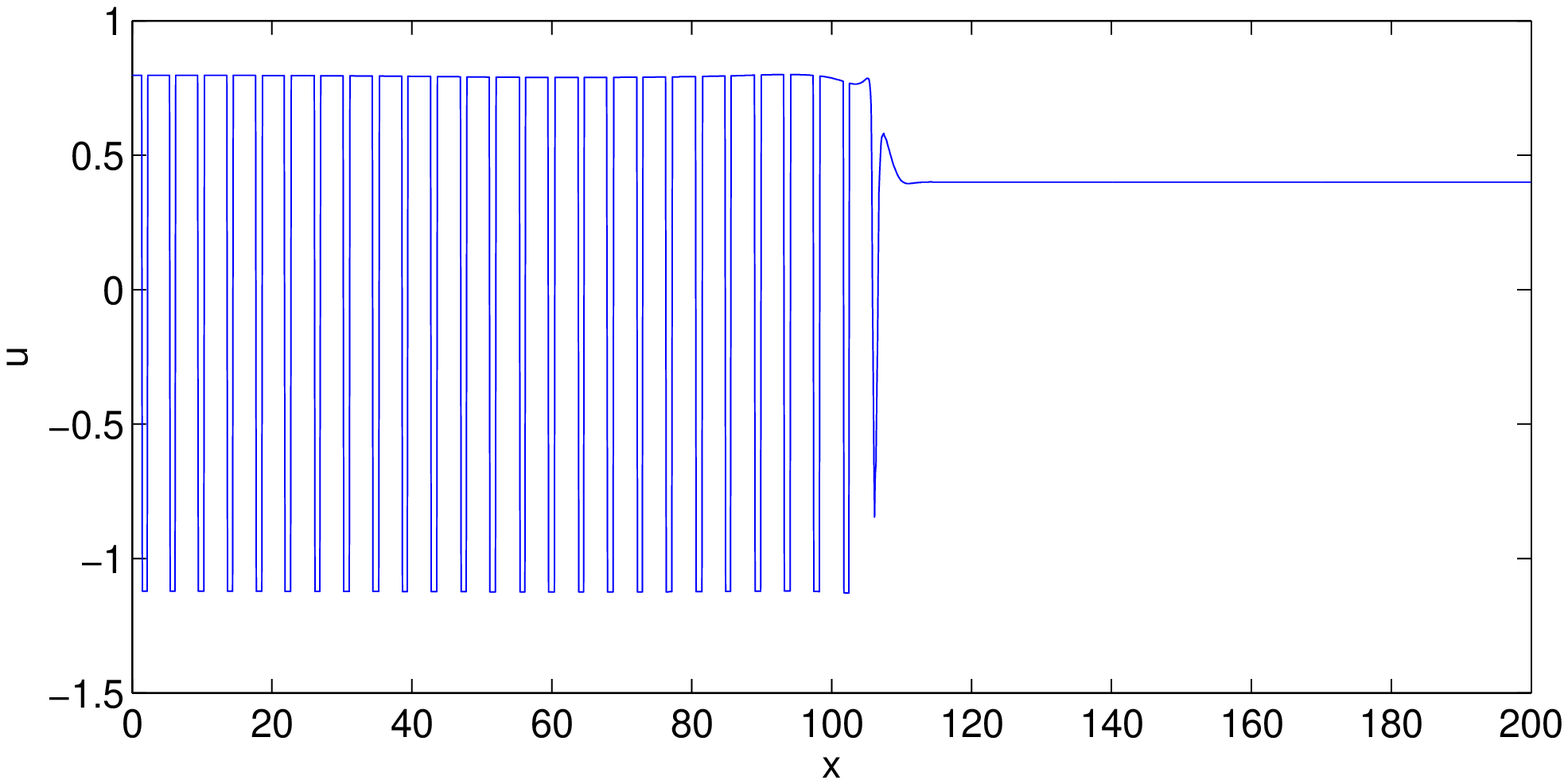}
\label{nonsymsimu}}
}
\quad
\mbox{
\subfigure[$\log|u-0.4|$ against at $t=140$.]
{
\includegraphics[width=0.45\textwidth,height=.35\textwidth]{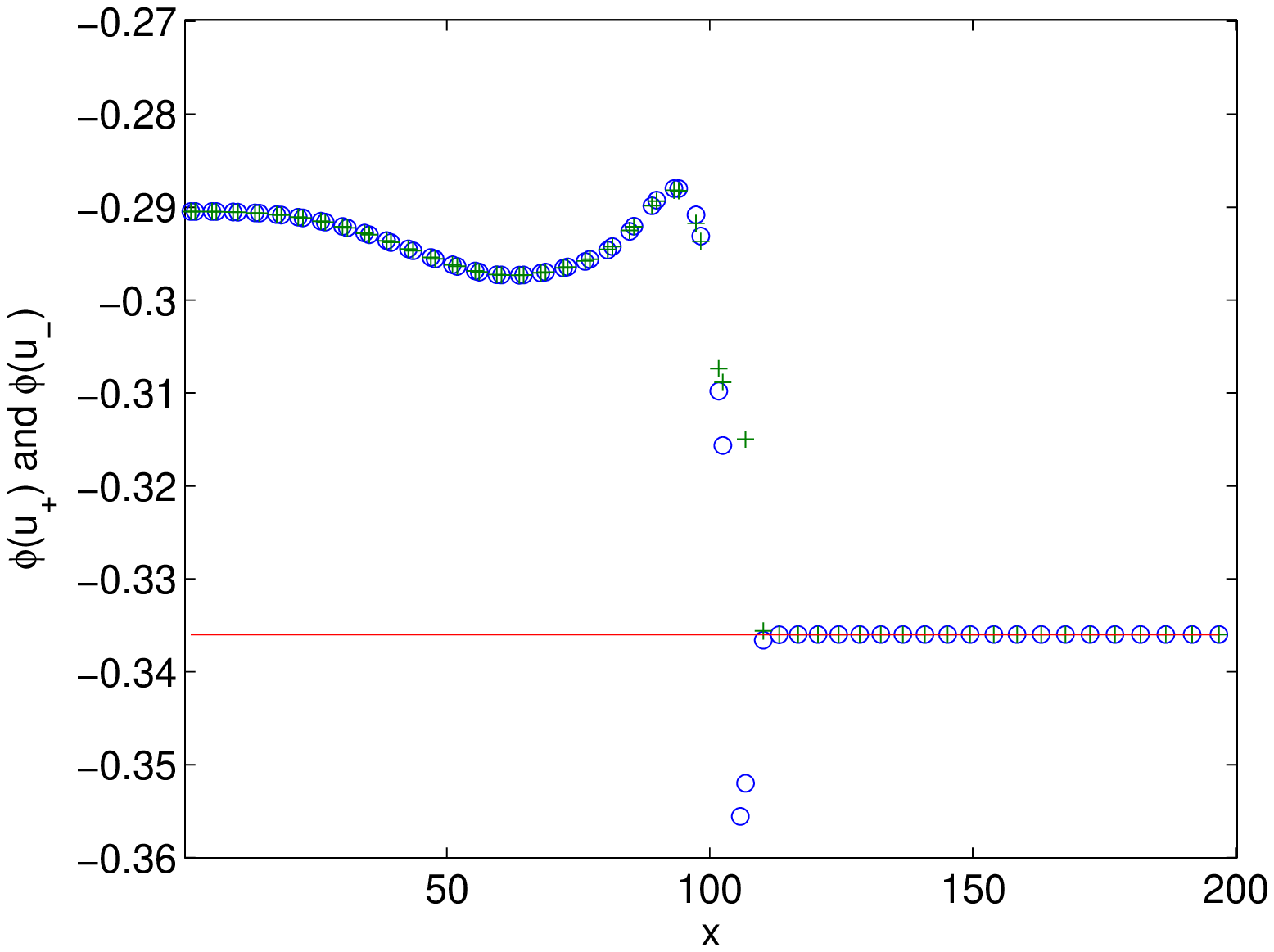}
\label{nonsymphicomp}}
}
\caption{A numerical computation of (\ref{main:eq}) with (\ref{phi:cubic}), $u_u=0.4$ and initial condition $u_0(x)=0.4+0.1 e^{-x^2}$ on the domain 
$[0,200]$. \subref{nonsymsimu} Show the solution at time $t=140$, and \subref{nonsymphicomp} shows a computation of $\phi(u_+)$ (circles) and 
$\phi(u_-)$ (crosses) at $t=140$, shown at the midpoints of $(x_j+x_{j_1})/2$ where $u(x_j)\approx u_+$ or $u_-$. The solid line indicates the 
value of $\phi(0.4)\approx -0.336$.     
}
\label{nonsym:phi:plusminus}
\end{figure}

Finally, we include some remarks about the non-symmetric cases.  
The analysis on front selection performed in Section~\ref{section:2} applies to the 
non-symmetric cases associated to the nonlinearities (\ref{phi:cubic}) and 
(\ref{phi:exposym}). The analysis on the transition region and the pattern also 
applies, the main difference with the symmetric cases being that, in general, 
$A\neq \phi(u_u)$. To illustrate this we include the numerical computation shown 
in Figure~\ref{nonsym:phi:plusminus}. 

The analysis for a $\phi$ of the form (\ref{phi:model2}) would be very different. In this case, there can only be one 
value of $u$ in the stable region satisfying (\ref{phi:const}). The solution might then be expected to approach 
the only available constant stable solution. It is, however, not immediately clear how such a solution arranges
 itself in space as $t\to +\infty$, since the conditions (\ref{mass}) hold; we venture that the solution oscillates 
spatially between a stable value $u_s$ and values that tend to infinity as $t\to +\infty$, presumably approximating a function of the form 
$u(x)=u_s + \sum_{n=1}^N M_n \,\delta(x-x_n)$.

\paragraph{Acknowledgements:} The authors gratefully acknowledge the support of the RTN project `Front-singularities'. C.~M.~Cuesta acknowledges the support of 
the Engineering and Physical Sciences Research Council in the form of a Postdoctoral fellowship 
(while at the University of Nottingham) and that of the MINECO through project MTM2011-24109. J.R. King is grateful to the hospitality of the ICMAT and its support through the MINECO: ICMAT Severo Ochoa project SEV-2011-0087.

\def\cprime{$'$}

\end{document}